\documentclass[11pt]{article}
\usepackage[T1]{fontenc}
\usepackage{latexsym,amssymb,amsmath,amsfonts,amsthm}
\usepackage{graphics}
\usepackage[demo]{graphicx}
\usepackage[section]{placeins}\usepackage{graphicx}
\usepackage{mathrsfs}
\usepackage{subfigure}
\usepackage{float}
\usepackage{color}

\usepackage{tikz}

\usepackage{mathtools}

\usepackage{epstopdf}
\newcommand{\R}{{\mathbb R}}

\newcommand{\N}{{\mathbb N}}
\newcommand{\C}{{\mathbb C}}

\newcommand{\be}{\begin{eqnarray}}
\newcommand{\ben}{\begin{eqnarray*}}
\newcommand{\en}{\end{eqnarray}}
\newcommand{\enn}{\end{eqnarray*}}

\newcommand{\real}{{\rm Re\,}}
\newcommand{\ima}{{\rm Im\,}}

\newcommand{\s}{\mathbb{S}}

\newtheorem{theorem}{Theorem}[section]

\newtheorem{lemma}[theorem]{Lemma}

\newtheorem{definition}[theorem]{Definition}
\newtheorem{remark}[theorem]{Remark}

\definecolor{rot}{rgb}{0,0,0}
\definecolor{hw}{rgb}{0,0,0}

\begin{document}
\renewcommand{\theequation}{\arabic{section}.\arabic{equation}}

\title{Imaging a moving point source from multi-frequency data measured at one and sparse observation directions (part I): far-field case}

	\author{Hongxia Guo\footnotemark[2], Guanghui Hu\footnotemark[2], Guanqiu Ma\footnotemark[1]\; \footnotemark[2]}
	
	\date{}
	\maketitle
	
	\renewcommand{\thefootnote}{\fnsymbol{footnote}}
	\footnotetext[1] {Corresponding author}
	\footnotetext[2]{School of Mathematical Sciences and LPMC, Nankai University, Tianjin, 300071, China(hxguo\_math@163.com, ghhu@nankai.edu.cn, gqma@nankai.edu.cn).}
	\renewcommand{\thefootnote}{\arabic{footnote}}

	\begin{abstract}
		We propose a multi-frequency algorithm for imaging the trajectory of a moving point source from one and sparse far-field observation directions in the frequency domain. The starting and terminal time points of the moving source are both supposed to be known.
 We introduce the concept of observable directions (angles) in the far-field region and derive all observable directions (angles) for straight and circular motions. At an observable direction, it is verified that the smallest trip containing the trajectory and perpendicular to the direction can be imaged, provided the orbit function possesses a certain monotonical property. Without the monotonicity one can only expect to recover a thinner strip.  
The far-field data measured at sparse observable directions can be used to recover the $\Theta$-convex domain of the trajectory. Both two- and three-dimensional numerical examples are implemented to show effectiveness and feasibility of the approach.

\vspace{.2in} {\bf Keywords}: {\bf inverse moving source problem, Helmholtz equation, multi-frequency data, factorization method, uniqueness.}
	\end{abstract}

	\section{Introduction}\label{sec:2}

	\subsection{Time-dependent model and Fourier transform}

We suppose that	the whole space $\mathbb{R}^d$ ($d =2,3$) is filled with a homogeneous and isotropic medium with a unit mass density. Consider a moving point source along the trajectory function $a(t):[t_{\min},t_{\max}]\rightarrow \R^d  \in C^1[t_{\min},t_{\max}]$ with $0 < t_{\min}<t_{\max}$.
The source function $S$ is supposed to radiate wave signals at the beginning time $t_{\min}$ and stop radiating at the time point $t_{\max}$, i.e., it is supported in the interval $[t_{\min}, t_{\max}] $ with respect to the time variable $t>0$. Hence,
the source function takes the form
	\begin{equation}
		S(x,t) = \delta(x-a(t)) \chi(t),
	\end{equation}
	where $\delta$ denotes the Dirac delta function and
		\begin{equation*}
		\chi(t) = \left\{
		\begin{aligned}
			& 1, \quad t \in [t_{\min}, t_{\max}],\\
			& 0, \quad t \notin [t_{\min}, t_{\max}],
		\end{aligned}
		\right.
	\end{equation*}
	is the characteristic function over the interval $[t_{\min}, t_{\max}]$.
Denote the trajectory by $\Gamma \coloneqq \{x: x=a(t),\, t\in [t_{\min}, t_{\max}]\}$.
One can easily find Supp $S(\cdot, t) \subset \Gamma$ for all $t\in [t_{\min}, t_{\max}]$.
The propagation of the radiated wave fields $U(x,t)$ is governed by the initial value problem
	\begin{equation}\label{timeeqn}
		\left\{
		\begin{aligned}
			&\frac{\partial^2 U}{\partial t^2} = \Delta U + S(x,t), \quad &&(x, t) \in \mathbb{R}^d \times \mathbb{R}^+, \mathbb{R}^+ \coloneqq \{t\in \R: t>0\},\\
			&U(x,0)=\partial_t U(x,0) = 0, &&x\in \mathbb{R}^d.
		\end{aligned}
		\right.
	\end{equation}

The solution $U$ can be written explicitly as the convolution of the fundamental solution $G_d (d=2,3)$ to the wave equation with the source term,
	\begin{equation}\label{time-solu}
		U(x,t) = G_d(x;t) * S(x,t) \coloneqq \int_{\mathbb{R}^+} \int_{\mathbb{R}^d}
		G_d(x-y; t-\tau)S(y,\tau)\,dyd\tau
	\end{equation}
	where
	\begin{equation*}
		G_d(x;t) = \left\{
		\begin{aligned}
			& \frac{H(t-|x|)}{2\pi \sqrt{t^2-|x|^2}}, \quad &&\quad\mbox{if}\quad d=2; 
			\\
			& \frac{\delta(t-|x|)}{4\pi |x|}, \quad && \quad\mbox{if}\quad d=3, 
		\end{aligned}\right.
	\end{equation*}
where $H$ denotes the Heaviside function. In this paper the one-dimensional Fourier and inverse Fourier transforms are defined by
	\ben
		(\mathcal{F}g)(k)=\frac{1}{\sqrt{2\pi}}\int_{\R}g(t)e^{-ikt}\,dt,\quad
		(\mathcal{F}^{-1}v)(t)=\frac{1}{\sqrt{2\pi}}\int_{\R}v(k)e^{ikt}\,dk,
	\enn
	respectively.
The Fourier transform of $S$ is thus given by
 	\begin{equation}\label{sourcef}
 		f(x,k) \coloneqq (\mathcal{F} S(x, \cdot))(k)=\frac{1}{\sqrt{2\pi}}\int_{\mathbb{R}} \delta(x-a(t))\chi(t) e^{-ikt}\, dt = \frac{1}{\sqrt{2\pi}}\int_{t_{\min}}^{t_{\max}} \delta(x-a(t)) e^{-ikt}\, dt.
 	\end{equation}
 	It is obvious $f(x,k) = 0$ for all $x\notin\Gamma$ and $k\in [k_{\min},k_{\max}]$.
 	From the expression \eqref{time-solu}, one deduces the Fourier transform of the wave fields $U$,
 	\begin{equation}\label{wFU}
 		\begin{aligned}
 		w(x,k)=(\mathcal{F}U)(x,k) &= \int_{\mathbb{R}^d} (\mathcal{F}G_d)(x-y;k) (\mathcal{F}S)(y,k)\,dy\\
 		&= \int_{\mathbb{R}^d} \Phi_d(x-y;k) f(y,k)\,dy.
 		\end{aligned}
 	\end{equation}
 	Here, $\Phi_d(x-y;k)$ is the fundamental solution to the Helmholtz equation $(\Delta + k^2)w = 0$, given by
\begin{equation*}
		\Phi_d(x-y;k) = \left\{
		\begin{aligned}
			&\frac{i}{4}H_0^{(1)}(k|x-y|), &&d=2,\\
			&\frac{e^{ik|x-y|}}{4 \pi |x-y|}, &&d=3,
		\end{aligned}
		\right.
		\qquad\quad  x \neq y ,\, x,y\in \mathbb{R}^d,
	\end{equation*}
and	$H_0^{(1)}$ is the Hankel function of the first kind of order zero.
 	On the other hand, taking the Fourier transform on the wave equation yields the inhomogeneous Helmholtz equations
 	\begin{equation}\label{eq1}
 		\Delta w(x,k) + k^2 w(x,k) = -f(x,k), \qquad x\in \R^{d}, \;k>0.
 	\end{equation}
From \eqref{wFU} we observe that $w$ satisfies 	
 the Sommerfeld radiation condition
 	\be\label{SRC}
	\lim\limits_{r \to \infty} r^{\frac{d-1}{2}} (\partial_r w - ikw) = 0,\quad r = |x|,\en
which holds uniformly in all directions $\hat{x}=x/r\in \s^{d-1}:=\{x\in \R^{d}: |x|=1\}$.

	\subsection{Formulation in the frequency domain and literature review}
	Denote by $[k_{\min}, k_{\max}]$ an interval of wavenumbers/frequencies on the positive real axis.
From the time-domain setting we see
	$$f(x,k) \neq 0, \, x\in \Gamma,\quad f(x,k) = 0, \, x\notin \Gamma$$
	for all $k\in [k_{\min},k_{\max}]$,
implying supp $f(\cdot, k) = \Gamma$ for all $k \in [k_{\min}, k_{\max}]$.
For every $k > 0$, the unique solution $w \in H^2_{loc}(\mathbb{R}^d)$ to \eqref{eq1}-\eqref{SRC} is given by \eqref{wFU}, i.e.,
	\begin{equation}\label{expression-w}
		w(x, k) = \int_{\mathbb{R}^d} \Phi_d(x-y;k) f(y, k) dy,  \quad x \in \mathbb{R}^d.
	\end{equation}
	The Sommerfeld radiation condition leads to the asymptotic behavior of $w$ at infinity:
	\begin{equation}\label{far-field}
		w(x) = C_d \frac{e^{ik|x|}}{ |x|^\frac{d-1}{2}} \{w^{\infty}(\hat{x}, k) + \mathcal{O}(r^{-\frac{d+1}{2}})\} \quad \text{as} \quad |x| \to \infty, \, d=2,3,
	\end{equation}
	where $C_2 = e^{i\pi /4}/\sqrt{8\pi k}$,  $C_3 = 1/4\pi$, and $w^\infty(\cdot, k)\in C^\infty(\s^{d-1})$ is known as the far-field pattern (or scattering amplitude) of $w$. It is well known that the function  $\hat{x}\mapsto w^\infty(\hat{x}, k)$ is real-analytic on $\s^{d-1}$, where $\hat{x}\in \s^{d-1}$ is usually referred as the observation direction. By \eqref{expression-w}, the far-field pattern $w^\infty$ of $w$ can be expressed as
	\be\nonumber\label{u-infty}
		w^\infty(\hat{x}, k)= \int_{\mathbb{R}^{d}} e^{-ik\hat{x}\cdot y} f(y, k)\,dy
		= \frac{1}{\sqrt{2\pi}} \int_{t_{\min}}^{t_{\max}} e^{-ik \big(a(t)\cdot \hat{x}+t \big)}\,dt
		\quad
	\en
for $\hat{x}\in \s^{d-1}$ and $k>0$.
	Noting that the time-dependent source $S$ is real valued, we have $f(x, -k)= \overline{f(x,k)}$ for all $k>0$ and thus $w^{\infty}(x, -k)=\overline{w^{\infty}(x, k)}$.

\begin{figure}[!ht]
		\centering
		\scalebox{0.5}{
		\begin{tikzpicture}

		\draw [very thick,dotted] (5,0) arc [ start angle = 0, end angle = 360, radius = 5];
		\draw [very thick,smooth] (-1.5,-1.5) .. controls (1,-2) .. (2,0);

		\draw (0,5) node [above] {$|x_j|\to \infty,\,w^{\infty}(\hat{x}_j,k)$};

		\draw (1,-2) node [below] {$\Gamma$};

		\draw (-1.5,-1.5) node [left] {$a(t_{\min})$};
		\draw (2,0) node [right] {$a(t_{\max})$};
		
		\fill (5,0) circle (3pt);
		\fill (0,5) circle (3pt);
		\fill (-5,0) circle (3pt);
		\fill (0,-5) circle (3pt);		
						
		\end{tikzpicture}
		}
		\caption{Imaging the trajectory $\Gamma$ from knowledge of multi-frequency far-field patterns measured at sparse observation directions
		$\hat{x}_j:=(\cos(j\pi/2),\sin(j\pi/2))$, $j=0,1,2,3$.}
		\label{ip}
	\end{figure}
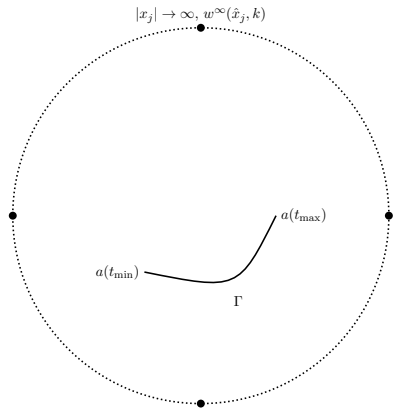

	In this paper we are interested in the following inverse problem (see Fig. \ref{ip}): 
	\begin{description}
	\item[(IP):] Recovery the trajectory $\Gamma$ from knowledge of the multi-frequency far-field patterns
	$$\{w^{\infty}(\hat{x}_j,k): k\in [k_{\min},k_{\max}],\,j=1,2,\cdots,M\}.$$
where $\hat{x}_j\in \s^{d-1}$ are sparse observation directions and $[k_{\min},k_{\max}]$ denotes a broad band of wavenumbers/frequencies.
\end{description}
In particular, we are interested the following question: 
\begin{description}
\item What kind information on $\Gamma$ can be extracted
from the the multi-frequency far-field patterns
	$\{w^{\infty}(\hat{x},k): k\in [k_{\min},k_{\max}]\}$ at a fixed observation direction $\hat{x}\in \s^{d-1}$ ?
\end{description}	
	 The above questions are of great importance in industrial, medical and military applications, because the number the measurement positions is usually quite limited and multi-frequency data are always available by  Fourier transforming the time-dependent measurement data. Although  multi-frequency far-field patterns are taken as the measurement data within this paper, the approach explored here  carries over naturally to the near-field data case at least in three dimensions.

To the best of the authors' knowledge, there are quite few mathematical studies on direct and inverse scattering theory for moving targets, in comparision with vast literatures devoted to scattering by stationary objects (see the monograph \cite{Isa1989}). Cooper \& Strauss \cite{Co79, CoStra} and Stefanov \cite{S1991} contribute rigorous mathematical theory to direct and inverse scattering from moving obstacles. We also refer to \cite{CB2008} for a linearized imaging theory with applications to various radar systems.
Recently there have been growing research interests in detecting the motion of moving point sources governed by inhomogeneous wave equations. Such kind of inverse source problems can be regarded as a linearized inverse obstacles problem. Consequently, various inversion algorithms have been proposed for recovering the orbit, profile and magnitude of a moving point source, such as the algebraic method \cite{NIO2012, T2020, Ohe2011}, the time-reversal method \cite{GF2015}, the method of fundamental solutions \cite{CGMS2020}, matched-filter and correlation-based imaging scheme \cite{FGPT2017},
 the iterative thresholding scheme \cite{Liu2021} and
 the Bayesian inference \cite{LGS2021, WKT2022}. See also \cite{LGY21, HKLZ2019, HKZ2020,HLY20,T2022} for uniqueness and stability results on inverse problems of identifying moving sources.

The purpose of this paper is to establish a factorization method for imaging $\Gamma$ from sparse far-field measurements at multiple frequencies.
The Factorization method was firstly proposed by Kirsch in 1998 \cite{KG08}. It has been successfully applied to various inverse scattering problems with multi-static data at a fixed energy (or equivalently, the Dirichlet-to-Neumann map). Its multi-frequency version was rigorously justified by Griesmaier and Schmiedecke \cite{GS2017} for inverse wave-number-independent source problems. It was verified in \cite{GS2017} that the smallest strip containing the support of a stationary source and perpendicular to a single observation direction in the far field  can be imaged. With sparse far-field observations, the so-called $\Theta$-convex polygon (that is, a convex polygonal whose normals coincide with observation directions) of the support can be recovered. If the dependance of the underlying source on the wavenumber takes the form of a windowed Fourier transform, one can also establish the analogue of the multi-frequency factorization method \cite{GGH2022}. The approach of \cite{GGH2022} also provides inspirations for dealing with other kinds of wave-number-dependent sources (or equivalently, time-dependent sources).   Although preliminary tests are implemented in \cite{GGH2022} for imaging the trajectory of a moving source, a comprehensive mathematical framework still needs to be built, which is the primary task of this work. Extensive numerical tests are implemented in the frequency domain in this paper. The counterpart of our inversion theory for wave equations using time-dependent near-field data deserves to be further investigated, which will be reported in our subsequent publications.

Motivated by earlier studies on sampling-type methods to inverse source problems \cite{GS2017, GGH2022, AHLS, GHZ22}, one can at most expect to recover the smallest strip containing the trajectory and perpendicular to the observation direction through the multi-frequency data measured at a single direction. 
However, our studies show that imaging such a strip turns out to be impossible for inverse moving source problems with a general orbit function.  
The recovery of the motion can be achieved only if the observation direction is {\em observable} in the sense of Definition \ref{obd} and the orbit function possesses a certain monotonicity property; see Theorem \ref{TH4.3} (ii). 
In particular, the monotonicity can be fulfilled if the velocity of the moving source is slower than the wave speed. Otherwise, 
 one can only get a thinner strip $K_{\Gamma}^{(\hat{x})}$ (see \eqref{K} for the definition, whose width is less than the aforementioned smallest strip) at an observable direction. For non-observable directions, the choice of the test function cannot lie in the range of the data-to-pattern operator (see Lemma \ref{lem3.4}). Hence, 
 it is impossible to extract any information on the motion of a moving source by our theory, although numerics still show  
partial information which however remains unclear to us. Using sparse observable directions, we design an indicator function for imaging the $\Theta$-convex domain of the trajectory. The $\Theta$-convex domain is a subset of the $\Theta$-convex polygon introduced in \cite{GS2017} and the $\Theta$-convex scattering support in \cite{SK}, because it 
is defined for observable directions only. Some uniqueness results will be summarized in Theorem \ref{TH4.3}, as a byproduct of the factorization scheme established in Theorems \ref{Th:factorization} and \ref{TH4.2}.




The remaining part is organized as follows. In Section \ref{sec:4}, the  multi-frequency far-field operator $\mathcal{F}^{(\hat{x})}$ for a fixed observation $\hat{x}$ is factorized in terms of the data-to-pattern operator $\mathcal{L}^{(\hat{x})}$, following the spirit of \cite{GGH2022}.
A range identity is given to connect the ranges of 
$\mathcal{F}^{(\hat{x})}$ and $\mathcal{L}^{(\hat{x})}$.
Section \ref{RangeLx} is devoted to the choice of test functions for characterizing the strip $K_{\Gamma}^{(\hat{x})}$ through analysis on the range of the data-to-pattern operator $\mathcal{L}^{(\hat{x})}$. In Section \ref{IdF} we define indicator functions using the far-field data measured at one or several observable directions. 2D and 3D numerical tests will be reported in the final Section \ref{num}.

	\section{Factorization of far-field operator}\label{sec:4}

The aim of this section is to explore the factorization method for recovering the trajectory $\Gamma = \text{Supp} f(\cdot, k)$ from the data measured at a single observation direction  $\hat{x}\in \mathbb{S}^{d-1}$. We shall proceed with the lines of \cite{GGH2022} to derive a factorization of the far-field operator $\mathcal{F}^{(\hat{x})}$.
Following the spirit of \cite{GS2017}, we introduce the central frequency $\kappa$ and half of the bandwidth of the given data $K$ as
	\begin{equation*}
		\kappa \coloneqq \frac{k_{\min} + k_{\max}}{2}, \quad K \coloneqq \frac{k_{\max} - k_{\min}}{2}.
	\end{equation*}
Define the far-field operator $\mathcal{F}^{(\hat{x})} : L^2(0, K) \to L^2(0, K)$ by
	\begin{equation}\label{FarO}
		(\mathcal{F}^{(\hat{x})}\phi)(\tau) \coloneqq \int_0^{K} w^{\infty}(x, \kappa + \tau - s)\,\phi(s)\,ds,\qquad \tau\in(0,K).
	\end{equation}
Recall from \eqref{u-infty} that $w^{\infty}$ is analytic in $k \in \mathbb{R}$. Hence the far-field operator $\mathcal{F}^{(\hat{x})}$ is linear and bounded. Further, it holds that	\begin{equation} \label{def:F}
		\begin{aligned}
			(\mathcal{F}^{(\hat{x})}\phi)(\tau) &= \int_0^{K}\int_{\mathbb{R}^d} e^{-i(\kappa + \tau - s)\hat{x}\cdot y}f(y, \kappa + \tau - s)\,dy\,\phi(s)\,ds \\
			&= \int_0^{K}\int_{\mathbb{R}^d} e^{-i(\kappa + \tau - s)\hat{x}\cdot y} \left(\frac{1}{\sqrt{2\pi}} \int_{t_{\min}}^{t_{\max}} e^{-i(\kappa + \tau - s)t} \delta(y-a(t)) \,dt \right)\,dy\,\phi(s)\,ds \\
			&= \frac{1}{\sqrt{2\pi}} \int_0^{K} \int_{t_{\min}}^{t_{\max}} e^{-i(\kappa + \tau - s)(t+ \hat{x}\cdot a(t))} \,dt \,\phi(s)\,ds
		\end{aligned}
	\end{equation}
	Below we shall prove a factorization of the above far-field operator.

	\begin{theorem}\label{facN}
		We have $\mathcal{F}^{(\hat{x})} =  \mathcal{L}  \mathcal{T}  \mathcal{L}^*$ where $ \mathcal{L} = \mathcal{L}^{(\hat{x})} : L^2(t_{\min},t_{\max}) \to L^2(0, K)$ is defined by
		\begin{equation}\label{tildeL}
			( \mathcal{L}\psi)(\tau) \coloneqq \int_{t_{\min}}^{t_{\max}} e^{-i\tau(t +\hat{x}\cdot a(t))}\psi(t)\, dt, \quad \tau \in (0, K)
		\end{equation}
		for all $\psi \in L^2(t_{\min},t_{\max})$. Here the middle operator $ \mathcal{T} :  L^2(t_{\min},t_{\max}) \to  L^2(t_{\min},t_{\max})$ is a multiplication operator defined by
		\begin{equation}\label{oT}
			( \mathcal{T}\varphi)(t) \coloneqq \frac{1}{\sqrt{2\pi}}e^{-i\kappa (t+\hat{x}\cdot a(t))} \varphi(t).
		\end{equation}
	\end{theorem}
\begin{remark}
In the remaining part of this paper the operator $\mathcal{L}$ will be referred to as the data-to-pattern operator corresponding to the orbit function $a(t)$. It is obvious that the far-field data \eqref{u-infty} can be expressed as $w^\infty(\hat{x},k)=(\mathcal{L}^{(\hat{x})}\, 1)(k)$.
We refer to \cite{KG08} for the analogue of the data-to-pattern operator for multi-static far-field operators at a fixed frequency.
\end{remark}
	\begin{proof}
		We first show that the adjoint operator $ \mathcal{L}^* : L^2(0, K) \to L^2(t_{\min},t_{\max})$ of $ \mathcal{L}$ can be expressed by
		\begin{equation}\label{aj-tildeL}
			( \mathcal{L}^* \phi)(t) \coloneqq \int_{0}^{K} e^{i s (t + \hat{x}\cdot a(t))}\phi(s)\, ds, \quad \phi \in L^2(0, K).
		\end{equation}
		Indeed, for $\psi\in L^2(t_{\min},t_{\max})$ and $\phi \in L^2(0,K)$, it holds that
		\begin{eqnarray*}
			\langle  \mathcal{L}\psi, \phi \rangle_{L^2(0, K)} &=& \int_0^{K} \left(  \int_{t_{\min}}^{t_{\max}} e^{-i\tau (t + \hat{x}\cdot a(t))} \psi(t)\,dt\right) \,\overline{\phi(\tau)}\,d\tau \\
			&=& \int_{t_{\min}}^{t_{\max}} \psi(t) \left(\int_0^{K} \overline{e^{i\tau(t + \hat{x}\cdot a(t))}\phi(\tau)}d\tau\right)\,dt \\
			&=&\langle \psi,  \mathcal{L}^* \phi \rangle_{L^2(t_{\min},t_{\max})}.
		\end{eqnarray*}
		which implies \eqref{aj-tildeL}. By the definition of $\mathcal{T}$, we have
		\begin{equation*}
			( \mathcal{T}  \mathcal{L}^* \phi)(t) = \frac{1}{\sqrt{2\pi}} e^{-i\kappa (t+\hat{x}\cdot a(t))}  \int_{0}^{K} e^{is(t+\hat{x}\cdot a(t))} \phi(s) \,ds,\quad \phi \in L^2(0, K).
		\end{equation*}
		Hence, using \eqref{sourcef} and \eqref{def:F},
		\begin{eqnarray*}
			( \mathcal{L}  \mathcal{T}  \mathcal{L}^* \phi)(\tau) &=& \int_{t_{\min}}^{t_{\max}} e^{-i\tau(t+\hat{x}\cdot a(t))} \left(\frac{1}{\sqrt{2\pi}} e^{-i\kappa (t+\hat{x}\cdot a(t))}   \int_{0}^{K} e^{is(t+\hat{x}\cdot a(t))} \phi(s) \,ds\right)dt\\
			&=& \frac{1}{\sqrt{2\pi}}\int_0^{K}\int_{t_{\min}}^{t_{\max}} e^{-i(\kappa + \tau - s)(t+\hat{x}\cdot a(t))} \,dt\,\phi(s)\,ds\\
			&=& (\mathcal{F}^{(\hat{x})} \phi)(\tau).
		\end{eqnarray*}
 		This proves the factorization $\mathcal{F}^{(\hat{x})} =  \mathcal{L}  \mathcal{T}  \mathcal{L}^*$.
	\end{proof}

	Denote by $\text{Range} ( \mathcal{L}^{(\hat{x})})$ the range of the data-to-pattern operator $ \mathcal{L}=\mathcal{L}^{(\hat{x})}$ (see \eqref{tildeL}) acting on $L^2(t_{\min},t_{\max})$.

	\begin{lemma}\label{com}
		The operator $ \mathcal{L} : L^2(t_{\min},t_{\max}) \to L^2(0, K)$ is compact with dense range.
	\end{lemma}

	\begin{proof}
		For any $\psi \in L^2(t_{\min},t_{\max})$, it holds that $ \mathcal{L}\psi \in H^1(0, K)$, which is compactly embedded into $L^2(0, K)$. This proves the compactness of $ \mathcal{L}$. By \eqref{aj-tildeL},
$(\mathcal{L}^* \phi)(t)$ coincides with the inverse Fourier transform of $\phi$ at the variable $t+\hat{x}\cdot a(t)$. Since the set $\{t+\hat{x}\cdot a(t): t\in [t_{\min}, t_{\max}]\}$ forms an interval of $\R$,
the relation $(\mathcal{L}^* \phi)(t)=0$ implies $\phi=0$ in $L^2(0, K)$. Hence, $ \mathcal{L}^*$ is injective. The denseness of $\text{Range} ( \mathcal{L}^{(\hat{x})})$ in $L^2(0, K)$ follows from the injectivity of $\mathcal{L}^*$.
	\end{proof}

Within the framework of Factorization method, it is essential to connect the ranges of $\mathcal{F}^{(\hat{x})}$ and $\mathcal{L}$. 
	We first recall that, for a bounded operator $F: Y\rightarrow Y$ in a Hilbert space $Y$ the real and imaginary parts of $F$ are defined respectively by
	\begin{equation*}
		\real F=\frac{F+F^*}{2},\quad \ima F=\frac{F-F^*}{2i},
	\end{equation*}
	which are both self-adjoint operators. Furthermore, by spectral representation we define the self-adjoint and positive operator $|\real F|$ as
	\begin{equation*}
		|\real F|=\int_{\R} |\lambda|\, d E_\lambda,\qquad \mbox{if}\quad \real F=\int_{\R} \lambda\, d E_\lambda.
	\end{equation*}
	The selfadjoint and positive operator $|\ima F|$ can be defined analogously.
 Introduce a new operator
	\begin{equation*}
		F_{\#}:=|\real F| +|\ima F|.
	\end{equation*}
	Since $F_{\#}$ is selfadjoint and positive, its square root $F_{\#}^{1/2}$ is defined as
	\begin{equation*}
		F_{\#}^{1/2}:=\int_{\R^+} \sqrt{\lambda}\, d E_\lambda,\qquad \mbox{if}\quad  F_{\#}=\int_{\R^+} \lambda\, d E_\lambda.
	\end{equation*}
	In this paper we need the following result from functional analysis.
	\begin{theorem}(\cite{GGH2022}) \label{range}
		Let $X$ and $Y$ be Hilbert spaces  and let $F: Y\rightarrow Y$, $L: X\rightarrow Y$, $T: X\rightarrow X$ be linear bounded operators such that $F=LTL^*$. We make the following assumptions
		\begin{itemize}
			\item[(i)] $L$ is compact with dense range and thus $L^*$ is compact and one-to-one.
			\item[(ii)] $\real T$ and $\ima T$ are both one-to-one,  and the operator $T_{\#}=|\real T| +|\ima T|: X\rightarrow X$ is coercive, i.e., there exists $c>0$ with
				\begin{equation*}
					\big\langle T_{\#}\, \varphi, \varphi\big\rangle\geq c\,||\varphi||^2\quad\mbox{for all}\quad \varphi\in X.
				\end{equation*}
		\end{itemize}
		Then the operator $F_{\#}$ is positive and  the ranges of $F_{\#}^{1/2}:Y\rightarrow Y$ and  $L: X\rightarrow Y$ coincide.
	\end{theorem}

	To apply Theorem \ref{range} to our inverse problem, we set
	\begin{equation*}
		F=\mathcal{F}^{(\hat{x})}, \quad L=\mathcal{L}, \quad T=\mathcal{T}, \quad X=L^2(t_{\min},t_{\max}),\quad Y=L^2(0, K),
	\end{equation*}
	where $\mathcal{T}$ is the multiplication operator of \eqref{oT}. It is easy to see
	\begin{eqnarray*}
		\left[(\real \mathcal{T})\, \varphi\right](t) &=& \frac{1}{\sqrt{2\pi}}\cos [\kappa (t+\hat{x}\cdot a(t))] \varphi(t),\\
		\left[(\ima \mathcal{T})\, \varphi\right](t) &=& -\frac{1}{\sqrt{2\pi}}\sin [\kappa (t+\hat{x}\cdot a(t))] \varphi(t)
	\end{eqnarray*}
 	are both one-to-one operators from $L^2(t_{\min},t_{\max})$ onto $L^2(t_{\min},t_{\max})$. The coercivity assumption of $\mathcal{F}^{(\hat{x})}$ yields the coercivity of  $\mathcal{T}_{\#}$. As a consequence of Theorem \ref{range}, we obtain
	\begin{equation}\label{RI}
		\mbox{Range}\, [(\mathcal{F}^{(\hat{x})})_{\#}^{1/2}]=\mbox{Range}\,(\mathcal{L}^{(\hat{x})})\quad\mbox{ for any }\, \hat{x}\in \mathbb{S}^{d-1}.
	\end{equation}

	Let $\varphi\in L^2(0, K)$ be a test function. We want to characterize the range of $\mathcal{L}^{(\hat{x})}$ through the choice of $\varphi$. Denote by $(\lambda_n^{(\hat x)}, \psi_n^{(\hat x)})$ an eigensystem of the positive and self-adjoint operator $ (\mathcal{F}^{(\hat{x})})_{\#}$, which is uniquely determined by the multi-frequency far-field patterns  $\{w^{\infty}(\hat{x}, k) : k \in (k_{\min}, k_{\max})\}$. Applying Picard's theorem and Theorem \ref{range}, we obtain
	\begin{equation}\label{Indicator}
		\varphi \in \text{Range}( \mathcal{L}^{(\hat{x})}) \quad \text{if and only if } \quad \sum\limits_{n=1}^{\infty} \frac{|\langle \varphi, \psi_n^{(\hat x)}\rangle|^2}{|\lambda_n^{(\hat x)}|} < +\infty.
	\end{equation}
	To establish the factorization method, we now need to choose a proper class of test functions which usually rely on a sample variable in $\mathbb{R}^d$. 

	\section{Range of $ \mathcal{L}^{(\hat{x})}$ and test functions}\label{RangeLx}

To characterize the range of $ \mathcal{L}^{(\hat{x})}$, we need to investigate monotonicity of the function $h(t):=\hat{x}\cdot a(t)+t\in C^1[t_{\min}, t_{\max}]$. For this purpose we define the division points of a continuous function over a closed interval.

\begin{definition}\label{DIDP}
		Let $f\in C[t_{\min}, t_{\max}]$. The point $t\in (t_{\min}, t_{\max})$ is called a division point if \\
		(1) $f(t)=0$;\\
		(2) There exist an $\epsilon_0>0$ such that either $|f(t+\epsilon)|>0$ or $|f(t-\epsilon)|>0$ for all $0<\epsilon<\epsilon_0$.
	\end{definition} 
	Obviously, the division points constitute a subset of the zero set of a continuous function. However, a division point cannot be an interior point of the zero set.	
Since $a(t)\in C^1[t_{\min}, t_{\max}]$, there are finitely many division points of the function $h'$, which we denote by
 $t_1<t_2<\cdots<t_{n-1}$.
 The interval $[t_{\min},t_{\max}]$ is then divided into $n$ sub-intervals $[t_{j-1},t_j]$, $j=1,2,\cdots,n$, where $t_{\min}=t_0$ and $t_{\max}=t_n$.
 Let $a_j$ and $h_j$ be the restrictions of $a$ and $h$ to $[t_{j-1},t_j]$, respectively.
Set $$\xi^{(\hat{x})}_{j,\min} := \inf\limits_{t\in [t_{j-1},t_j]} \{h_j(t)\},\quad \xi^{(\hat{x})}_{j,\max} := \sup\limits_{t\in [t_{j-1},t_j]} \{h_j(t)\},\quad j=1,2,\cdots n.$$
In each subinterval $(t_{j-1}, t_{j})$, one of
 following cases must hold:
	\begin{itemize}
		\item $h'_j(t)>0$ for all $t\in (t_{j-1},t_j)$. There holds
		$$\xi^{(\hat{x})}_{j,\min} = t_{j-1}+\hat{x}\cdot a_j(t_{j-1}), \quad \xi^{(\hat{x})}_{j,\max} = t_j+\hat{x}\cdot a_j(t_j);$$
		\item $h'_j(t)<0$ for all $t\in (t_{j-1},t_j)$. We have
		$$\xi^{(\hat{x})}_{j,\min} = t_j+\hat{x}\cdot a_j(t_j), \, \quad\xi^{(\hat{x})}_{j,\max} = t_{j-1}+\hat{x}\cdot a_j(t_{j-1});$$
		\item $h'_j(t)=0$ for all $t\in (t_{j-1},t_j)$. Consequently,
		$$\xi^{(\hat{x})}_{j,\min} = \xi^{(\hat{x})}_{j,\max} = t+\hat{x}\cdot a_j(t),\quad t\in[t_{j-1},t_j].$$
	\end{itemize}
Define
	\be\label{def:xi}\xi^{(\hat{x})}_{\min} := \min\limits_j \xi^{(\hat{x})}_{j,\min}=
	\inf\limits_{t\in [t_{\min},t_{\max}]} \{h(t)\}
	,\quad \xi^{(\hat{x})}_{\max} := \max\limits_j \xi^{(\hat{x})}_{j,\max}
	=\sup\limits_{t\in [t_{\min},t_{\max}]} \{h(t)\},
	\en
which denote the minimum and maximum of $h$ over $[t_{\min}, t_{\max}]$, respectively.
	 If $|h'_j(t)|>0$, the monotonicity of the function $\xi=h_j(t)$ for $t\in[t_j, t_{j-1}]$ implies the inverse function $t=h_j^{-1}(\xi)\in C^1[\xi^{(\hat{x})}_{j,\min}, \xi^{\hat{x})}_{j,\max}]$.
Set $$J=\{j\in \N : 1\leq j\leq n, h_j^{\prime}(t)\equiv 0, t\in (t_{j-1},t_j)\}.$$ and assume $h_j(t)\equiv c_j\in \R$ for $j\in J$. Note that it is possible that $J=\emptyset$.

With these notations we can rephrase the operator $\mathcal{L}^{(\hat{x})}$ defined by \eqref{tildeL} as
	\begin{equation}\label{pic-L}
		\begin{aligned}
			(\mathcal{L}^{(\hat{x})} \psi)(\tau) &= \sum\limits_{j=1}^{n} \int_{t_{j-1}}^{t_j} e^{-i\tau h_j(t)} \psi (t) \, dt \\
			&= \sum\limits_{j\notin J} \int_{t_{j-1}}^{t_j} e^{-i\tau h_j(t)} \psi (t) \, dt + \sum\limits_{j\in J} e^{-i\tau c_j} \int_{t_{j-1}}^{t_j}  \psi (t) \, dt.
		\end{aligned}
	\end{equation}
For $j\in J$, using $e^{-i\tau c} = \sqrt{2\pi} \mathcal{F}\delta(t-c)$ we can rewrite each terms in the second sum as
	\begin{equation}\label{delta-L}
		e^{-i\tau c_j} \int_{t_{j-1}}^{t_j}  \psi (t) \, dt = \sqrt{2\pi} \mathcal{F}\delta(t-c_j) \int_{t_{j-1}}^{t_j}  \psi (t) \, dt.
	\end{equation}
For $j\notin J$ and $h'_j(t)>0$, the integral in the first summation on the right hand of \eqref{pic-L} takes the form
	\ben
		\int_{t_{j-1}}^{t_j} e^{-i\tau h_j(t)}\psi(t)\,dt &=& \int_{\xi^{(\hat{x})}_{j,\min}}^{\xi^{(\hat{x})}_{j,\max}} e^{-i\tau \xi} \psi(h_j^{-1}(\xi))\,(h_j^{-1}(\xi))^{\prime}\,d\xi\\
	&=& \int_{\xi^{(\hat{x})}_{j,\min}}^{\xi^{(\hat{x})}_{j,\max}} e^{-i\tau \xi} \psi(h_j^{-1}(\xi))|(h_j^{-1}(\xi))^{\prime}|\,d\xi.	
	\enn
Note that $[h_j^{-1}(\xi)]^{\prime}>0$, due to the relation $h_j^{\prime}(t) [h_j^{-1}(\xi)]^{\prime} =1$. Analogously, if $h'_j(t)<0$ for some $j\notin J$, we have $[h_j^{-1}(\xi)]^{\prime}<0$ and thus	
\begin{equation*}
		\begin{aligned}
			\int_{t_{j-1}}^{t_j} e^{-i\tau h_j(t)}\psi(t)\,dt &= -\int_{\xi^{(\hat{x})}_{j,\min}}^{\xi^{(\hat{x})}_{j,\max}} e^{-i\tau \xi} \psi(h_j^{-1}(\xi))(h_j^{-1}(\xi))^{\prime}\,d\xi \\
			&= \int_{\xi^{(\hat{x})}_{j,\min}}^{\xi^{(\hat{x})}_{j,\max}} e^{-i\tau \xi} \psi(h_j^{-1}(\xi))|(h_j^{-1}(\xi))^{\prime}|\,d\xi.
		\end{aligned}
	\end{equation*}
Now, extending $h_j^{-1}$ by zero from $(\xi^{(\hat{x})}_{j,\min},\xi^{(\hat{x})}_{j,\max})$ to $\mathbb{R}$ and extending $\psi\in L^2(t_{\min},t_{\max})$ by zero to $L^2(\mathbb{R})$, we can write each term for $j\notin J$ as

	\begin{equation}\label{con-L}
		\int_{t_{j-1}}^{t_j} e^{-i\tau h_j(t)}\psi(t)\,dt = \int_{\mathbb{R}} e^{-i\tau \xi} \psi(h_j^{-1}(\xi))|(h_j^{-1}(\xi))^{\prime}|\,d\xi.
	\end{equation}

	Combining \eqref{pic-L}, \eqref{delta-L} and \eqref{con-L}, we get
	\begin{equation}\label{Lxi}
			(\mathcal{L}^{(\hat{x})}\psi)(\tau) = \int_{\mathbb{R}} e^{-i\tau\xi} g(\xi)\,d\xi,	
	\end{equation}
	with $$g(\xi)=\sum\limits_{j\notin J} \psi(h_j^{-1}(\xi))\, |(h_j^{-1}(\xi))'| + \sum\limits_{j\in J} \delta(\xi-c_j)\int_{t_{j-1}}^{t_j} \psi(t)\,dt.$$
	Note that $g$ is a generalized function if $J\neq \emptyset$ and that $g$ coincides with the inverse Fourier transform of $\mathcal{L}^{(\hat{x})} \psi$ up to some constant.
	Since supp $h_j^{-1} \subset [\xi_{\min}^{(\hat{x})},\xi_{\max}^{(\hat{x})}]$ for $j\notin J$ and $c_j\in [\xi_{\min}^{(\hat{x})},\xi_{\max}^{(\hat{x})}]$, we may estimate that the support of $g$ (equivalently, the inverse Fourier transform of $\mathcal{L}^{(\hat{x})} \psi$) as follows:
	$$\mbox{supp}(g(\xi)) \subset [\xi_{\min}^{(\hat{x})},\xi_{\max}^{(\hat{x})}].$$
Summing up the above arguments we arrive at
	\begin{lemma}\label{lem:supp}
		Let $\Gamma = \{y: y=a(t), t\in [t_{\min},t_{\max}]\} \subset \mathbb{R}^d$ be a $C^1$-smooth curve with $t_{\max}>t_{\min}$. 
		Then
		\begin{equation}
			(\mathcal{F}^{-1}\mathcal{L}^{(\hat{x})}\psi)(\xi) = \sqrt{2\pi} \left(\sum\limits_{j\notin J} \psi(h_j^{-1}(\xi))\, |(h_j^{-1}(\xi))'| + \sum\limits_{j\in J} \delta(\xi-c_j)\int_{t_{j-1}}^{t_j} \psi(t)\,dt\right).
\end{equation}	
Moreover,
		\begin{equation*}
			{\rm supp} (\mathcal{F}^{-1}\mathcal{L}^{(\hat{x})}\psi) \subset [\xi^{(\hat{x})}_{\min},\xi^{(\hat{x})}_{\max}].
		\end{equation*}
		\end{lemma}
		Below we provide a sufficient condition to ensure trivial intersections of the ranges of two data-to-pattern operators corresponding to different trajectories. 		
\begin{lemma}\label{lem:interrange}
 		Let $\Gamma_a = \{y: y=a(t), t\in [t_{\min},t_{\max}]\} \subset \mathbb{R}^d$ and $\Gamma_b = \{y: y=b(t), t\in [t_{\min},t_{\max}]\} \subset \mathbb{R}^d$ be $C^1$-smooth curves such that 
 		\be\nonumber
 			&&\left(\inf\limits_{t\in [t_{\min},t_{\max}]}(t + \hat{x}\cdot a(t)),\sup\limits_{t\in [t_{\min},t_{\max}]}(t + \hat{x}\cdot a(t))\right) \\ \label{condition-T}
			&\bigcap& \left(\inf\limits_{t\in [t_{\min},t_{\max}]}(t + \hat{x}\cdot b(t)),\sup\limits_{t\in [t_{\min},t_{\max}]}(t + \hat{x}\cdot b(t))\right) = \emptyset.
		\en
Let $ \mathcal{L}^{(\hat{x})}_a$ and $ \mathcal{L}^{(\hat{x})}_b$ be the data-to-pattern operators associated with  $\Gamma_a$ and $\Gamma_b$, respectively.
		Then $\mbox{Range}( \mathcal{L}^{(\hat{x})}_a) \cap \mbox{Range}( \mathcal{L}^{(\hat{x})}_b) = \{0\}$.
	\end{lemma}

	\begin{proof}
		Let $f_a,f_b \in L^2(t_{\min},t_{\max})$ be such that $ (\mathcal{L}^{(\hat{x})}_a f_a)(\tau)=   (\mathcal{L}^{(\hat{x})}_b f_b)(\tau) \coloneqq  Q(\tau, \hat{x})$. We need to prove $Q(\cdot, \hat{x}) \equiv 0$. By the definition of $ \mathcal{L}$ (see \eqref{tildeL}), the function
		\begin{equation*}
			 \tau \to Q(\tau, \hat{x}) =  \int_{t_{\min}}^{t_{\max}} e^{-i\tau (t+ \hat{x}\cdot a(t))} f_a(t)\,dt
			=  \int_{t_{\min}}^{t_{\max}} e^{-i\tau (t + \hat{x}\cdot b(t))} f_b(t)\,dt
		\end{equation*}
		belongs to $L^2(0, K)$.
		Since $Q(\tau, \hat{x})$ is analytic in $\tau\in \R$, the previous relation is well defined for any $\tau \in \mathbb{R}$.
		By Definition \ref{DIDP}, we suppose that $\{t_j\}_{j=1}^{n-1}$ and $\{\tilde{t}_j\}_{j=1}^{m-1}$ are division points of the functions $h_{a}(t) = t+\hat{x}\cdot a(t)$ and $h_{b}(t) = t+\hat{x}\cdot b(t)$, respectively.
Analogously we define $h_{j,a}(t):= t+\hat{x}\cdot a_j(t)$, $h_{j,b}(t) := t+\hat{x}\cdot b_j(t)$, and $J_a:=\{j\in \N : 1\leq j\leq n, h_{j,a}^{\prime}(t)\equiv 0, t\in (t_{j-1},t_j)\}$, $J_b:=\{j\in \N : 1\leq j\leq m, h_{j,b}^{\prime}(t)\equiv 0, t\in (\tilde{t}_{j-1},\tilde{t}_j)\}$. Denote $h_{j,a}(t)\equiv c_{j,a}$ for $j\in J_a$ and $h_{j,b}(t)\equiv c_{j,b}$ for $j\in J_b$.
Using the formula \eqref{Lxi}, the function $ Q(\cdot, \hat{x})$ can be rewritten as the Fourier transforms:
		\begin{equation}\label{G}
			 {Q}(\tau, \hat{x}) = \int_{\mathbb{R}} e^{-i\tau  \xi}  {g}_a( \xi, \hat{x})\,d \xi =\int_{\mathbb{R}} e^{-i\tau  \xi}  {g}_b( \xi, \hat{x})\,d \xi,
		\end{equation}
		with
		\begin{equation*}
			g_a(\xi,\hat{x})=\sum\limits_{j\notin J_a} f_a(h_{j,a}^{-1}(\xi))\, |(h_{j,a}^{-1}(\xi))'| + \sum\limits_{j\in J_a} \delta(\xi-c_{j,a})\int_{t_{j-1}}^{t_j} f_a(t)\,dt,
		\end{equation*}
		\begin{equation*}
			g_b(\xi,\hat{x})=\sum\limits_{j\notin J_b} f_b(h_{j,b}^{-1}(\xi))\, |(h_{j,b}^{-1}(\xi))'| + \sum\limits_{j\in J_b} \delta(\xi-c_{j,b})\int_{\tilde{t}_{j-1}}^{\tilde{t}_j} f_b(t)\,dt.
		\end{equation*}
This implies $g_a(\xi,\hat{x})=g_b(\xi,\hat{x})$ for all $\xi\in R$.	On the other hand,	the support sets of $g_a$ and $g_b$ satisfy
		\begin{equation*}
			\begin{aligned}
				&{\rm supp}\,  {g}_a(\cdot, \hat{x}) \subset \left(\inf\limits_{t\in [t_{\min},t_{\max}]}(t + \hat{x}\cdot a(t)),\sup\limits_{t\in [t_{\min},t_{\max}]}(t + \hat{x}\cdot a(t))\right) ,\\
				&{\rm supp}\,  {g}_b(\cdot, \hat{x}) \subset \left(\inf\limits_{t\in [t_{\min},t_{\max}]}(t + \hat{x}\cdot b(t)),\sup\limits_{t\in [t_{\min},t_{\max}]}(t + \hat{x}\cdot b(t))\right).
			\end{aligned}
		\end{equation*}
Hence,	by the condition \eqref{condition-T} we obtain $ {g}_a( \xi, \hat{x}) =  {g}_b( \xi, \hat{x}) \equiv 0$ for all $ \xi \in \mathbb{R}$ . In view of \eqref{G}, we get $ Q(\cdot,\hat{x}) \equiv 0$. 	
	\end{proof}


	For any $y\in \R^d$, define the parameter-dependent test functions $\phi^{(\hat{x})}_{y}\in L^2(0, K)$  by
	\begin{equation}\label{testfunc}
		\phi^{(\hat{x})}_{y}(k)=\frac{1}{|t_{\max}-t_{\min}|}\int_{t_{\min}}^{t_{\max}} e^{-ik (\hat{x}\cdot y+t)} dt, \quad k\in(0, K).
	\end{equation}
Here we stress that the test function $\phi^{(\hat{x})}_{y}$ depends on both the observation direction $\hat{x}\in \mathbb{S}^{d-1}$ and the space variable $y\in \R^d$.
The supporting information of the inverse Fourier transform of the above test function is described as follows.
	\begin{lemma}\label{INVF-PHI}
		We have
		\begin{equation}\label{inv-phi}
			[\mathcal{F}^{-1} \phi^{(\hat{x})}_{y}](\tau)=
			\left\{\begin{array}{lll}
				\sqrt{2\pi}/|t_{\max}-t_{\min}| \quad &&\mbox{if}\quad \tau\in [\hat{x}\cdot y+t_{\min},\; \hat{x}\cdot y+t_{\max}], \\
				0 \quad&&\mbox{if otherwise}.
			\end{array}\right.
		\end{equation}
	\end{lemma}
	
	\begin{proof}
		Let $\tau = \hat{x}\cdot y+t$, we can rewrite the function $\phi^{(\hat{x})}_{y}$ as
		\begin{equation*}
			\phi^{(\hat{x})}_{y}(k)=\int_{\R} e^{-ik \tau}g_y(\tau, \hat{x})\,d\tau,
			\end{equation*}
			where\begin{equation*}
			g_y(\tau,\hat{x}):=\left\{
			\begin{aligned}
			&\frac{1}{|t_{\max}-t_{\min}|} &&{\rm if}\, \tau \in [\hat{x}\cdot y +t_{\min}, \hat{x}\cdot y +t_{\max}],\\
			&0 &&\mbox{ if  otherwise}.
			\end{aligned}\right.
		\end{equation*}
		Therefore, $[\mathcal{F}^{-1}\phi^{(\hat{x})}_{y}](\tau) = \sqrt{2\pi}g_y(\tau,\hat{x})$.
	\end{proof}
In the following we present a necessary condition imposed on the observation direction $\hat{x}$ and radiating period $T:=t_{\max}-t_{\min}$ to ensure that the test function $\phi^{(\hat{x})}_y$ lies in the range of the data-to-pattern operator.	
	\begin{lemma}\label{lem3.5}
		If $\phi^{(\hat{x})}_y \in {\rm Range} (\mathcal{L}^{(\hat{x})})$ for some $y\in \R^d$, we have $\xi^{(\hat{x})}_{\max} - \xi^{(\hat{x})}_{\min} \geq T$. Here $\xi^{(\hat{x})}_{\max}$ and $\xi^{(\hat{x})}_{\min}$ are defined by \eqref{def:xi}.
			\end{lemma}

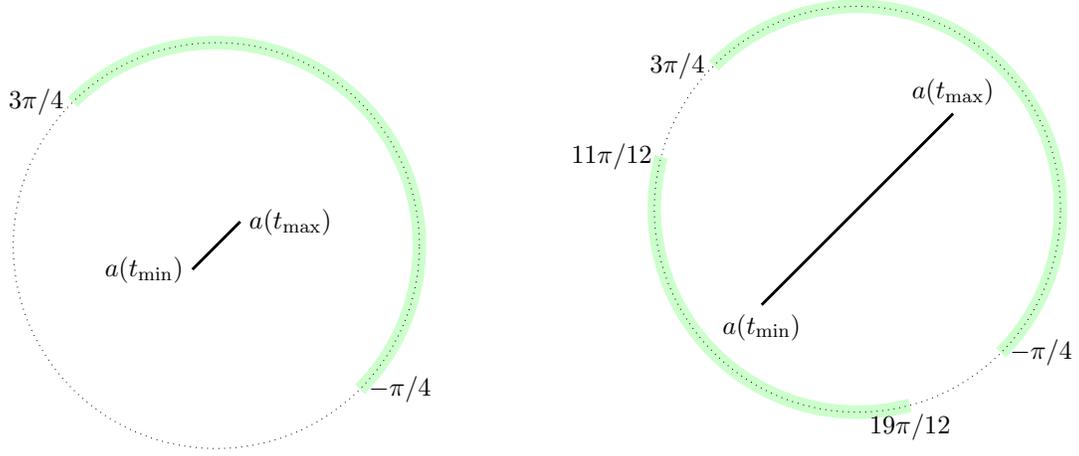
\begin{figure}
\begin{minipage}[b]{.5\textwidth}
		\centering
		\scalebox{0.9}{
		\begin{tikzpicture}
		\fill[green!20] (0:2.9) --(0:3.1) arc (0:135:3.1) --(135:2.9) arc (135:0:2.9);

		\fill[green!20] (315:2.9) --(315:3.1) arc (315:360:3.1) --(360:2.9) arc (360:315:2.9);

		\draw [dotted] (3,0) arc [ start angle = 0, end angle = 360, radius = 3];
		\draw [very thick] (-0.3535, -0.3535) -- (0.3535,0.3535);

		\draw (-0.3535, -0.3535) node [left] {$a(t_{\min})$};

		\draw (0.3535,0.3535) node [right] {$a(t_{\max})$};
		\draw (135:3) node [left] {$3\pi/4$};
		\draw (315:3) node [right] {$-\pi/4$};


		\end{tikzpicture}
		}
		\end{minipage}%
\begin{minipage}[b]{.5\textwidth}
	\centering
		\scalebox{0.9}{
		\begin{tikzpicture}
		\fill[green!20] (0:2.9) --(0:3.1) arc (0:135:3.1) --(135:2.9) arc (135:0:2.9);

		\fill[green!20] (315:2.9) --(315:3.1) arc (315:360:3.1) --(360:2.9) arc (360:315:2.9);

		\fill[green!20] (165:2.9) --(165:3.1) arc (165:180:3.1) --(180:2.9) arc (180:165:2.9);

		\fill[green!20] (180:2.9) --(180:3.1) arc (180:285:3.1) --(285:2.9) arc (285:180:2.9);

		\draw [dotted] (3,0) arc [ start angle = 0, end angle = 360, radius = 3];
		\draw [very thick] (-1.414, -1.414) -- (1.414,1.414);

		\draw (-1.414, -1.414) node [below] {$a(t_{\min})$};

		\draw (1.414,1.414) node [above] {$a(t_{\max})$};
		\draw (135:3) node [left] {$3\pi/4$};
		\draw (315:3) node [right] {$-\pi/4$};
		\draw (165:3) node [left] {$11\pi/12$};
		\draw (285:3) node [below] {$19\pi/12$};


		\end{tikzpicture}
		}
		\end{minipage}
		\caption{Illustration of observable (green arc) and non-observable (dotted arc) directions for the trajectory $a(t) = c\sqrt 2/2(t,t)$ for $t\in [1,2]$ with $c=1$ (left) and $c=4$ (right).		
}\label{S2}
\end{figure}

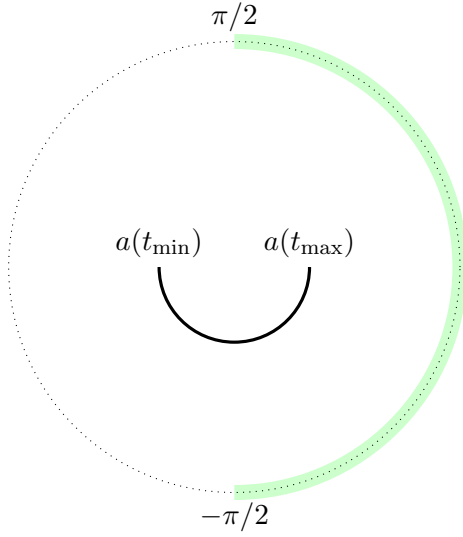
\begin{figure}[ht]
		\centering
		\scalebox{1}{
		\begin{tikzpicture}
		\fill[green!20] (0:2.9) --(0:3.1) arc (0:90:3.1) --(90:2.9) arc (90:0:2.9);

		\fill[green!20] (270:2.9) --(270:3.1) arc (270:360:3.1) --(360:2.9) arc (360:270:2.9);

		\draw [dotted] (3,0) arc [ start angle = 0, end angle = 360, radius = 3];
		\draw [very thick] (-1,0) arc [ start angle =180, end angle =360 ,radius =1];

		\draw (-1, 0) node [above] {$a(t_{\min})$};

		\draw (1, 0) node [above] {$a(t_{\max})$};
		\draw (270:3) node [below] {$-\pi/2$};
		\draw (90:3) node [above] {$\pi/2$};

		\end{tikzpicture}
		}
		\caption{Illustration of the observable (green arc) and non-observable (dotted arc) directions for the trajectory $a(t) = (\cos t+1,\sin t+2)$ with $t\in [\pi,2\pi]$.}\label{A1}
		
	\end{figure}

	\begin{proof}
		If $\phi^{(\hat{x})}_y \in {\rm Range} (\mathcal{L}^{(\hat{x})})$, there exists a function $\psi\in L^2(t_{\min},t_{\max})$ such that $\phi^{(\hat{x})}_y = \mathcal{L}^{(\hat{x})} \psi$ in $L^2(0,K)$. Since both $\phi^{(\hat{x})}_y$ and $\mathcal{L}^{(\hat{x})} \psi$ are analytic functions over $\R$, it holds that $\phi^{(\hat{x})}_y(k) = (\mathcal{L}^{(\hat{x})} \psi) (k)$ for all $k\in \R$. Then their support sets must be identical, i.e., supp$(\mathcal{F}^{-1}\phi^{(\hat{x})}_y)=$ supp$(\mathcal{F}^{-1}\mathcal{L}^{(\hat{x})} \phi) \subset [\xi^{(\hat{x})}_{\min},\xi^{(\hat{x})}_{\max}]$,
where we have used Lemma \ref{lem:supp}. Hence, the length of supp$(\mathcal{F}^{-1}\phi^{(\hat{x})}_y)$, which can be seen from Lemma \ref{INVF-PHI}, must be less than or equal to that of $[\xi^{(\hat{x})}_{\min},\xi^{(\hat{x})}_{\max}]$, i.e.,		
\ben \xi^{(\hat{x})}_{\max} - \xi^{(\hat{x})}_{\min} \geq t_{\max} - t_{\min}=T.\enn		
	\end{proof}	
From the above lemma we conclude that $\phi^{(\hat{x})}_y \notin {\rm Range} (\mathcal{L}^{(\hat{x})})$ for all $y\in \R^d$, if $\xi^{(\hat{x})}_{\max} - \xi^{(\hat{x})}_{\min}< T$. Inspired by this fact we introduce the concept of observable directions.
 \begin{definition}\label{obd}
		Let $\xi^{(\hat{x})}_{\min}$ and $\xi^{(\hat{x})}_{\max}$ be the maximum and minimum of the function $h(t)=\hat{x}\cdot a(t)+t\in C^1[t_{\min}, t_{\max}]$ (see \eqref{def:xi}).
		The unit vector $\hat{x}\in \s^{d-1}$ is called an observable direction if $\xi^{(\hat{x})}_{\max} - \xi^{(\hat{x})}_{\min} \geq T$. The direction $\hat{x}$ is called  non-observable if $\xi^{(\hat{x})}_{\max} - \xi^{(\hat{x})}_{\min} < T$.
	\end{definition}
We remark that the set of observable directions is uniquely determined by the orbit function $a(t)$ together with the starting and terminal time points $t_{\min}$ and $t_{\max}$. For non-observable directions $\hat{x}$, one cannot extract information on the orbit function by our approach, which will be explained in the second assertion of Theorem \ref{Th:factorization} below. If $\hat{x}$ is observable and $h$ is monotonically increasing, the smallest strip containing the trajectory and perpendicular to $\hat{x}$ can be recovered.
If $\hat{x}$ is observable and $h$ is not monotonically increasing, another thinner strip perpendicular to $\hat{x}$  can be imaged. Below we derive the observable directions for orbit functions given by a straight line (see Fig. \ref{S2}) and a semi-circle (see Fig. \ref{A1}) in two dimensions. We refer to Section \ref{num} for further discussions on piecewise linear curves in 2D and a straight line segment in 3D.

\vspace{.1in}	
\textbf{Example 1: A straight line segment in $\R^2$.}

Suppose that an acoustic point source is moving along the straight line
$a(t) = ct( \cos\alpha,\sin\alpha)\in \R^2$ for $t\in[t_{\min},t_{\max}]$, where $c>0$ denotes the velocity and $\alpha\in[0, 2\pi]$ the angle between the trajectory and the $x_1$-axis.

\begin{lemma}\label{line-ob}
		(i) If $c\, \leq 2$, the direction $\hat{x}= (\cos \theta,\sin \theta)$ is observable if $\theta \in [\alpha-\pi/2,\alpha+\pi/2]$.

		(ii) If $c> 2$, the direction $\hat{x}= (\cos \theta,\sin \theta)$ is observable if $\theta \in [\alpha-\pi/2,\alpha+\pi/2]\cup [\alpha+ \arccos(-2/c)), \alpha+2\pi - \arccos(-2/c))]$.
\end{lemma}

		\begin{proof}
		From the expression of the orbit function $a(t)$, we have
		\begin{equation*}
			\begin{aligned}
				h(t) = t+\hat{x}\cdot a(t) &= t + t\, c\, (\cos \theta \cos \alpha +\sin \theta \sin \alpha)
				= t(1+\, c\, \cos (\theta-\alpha)),\\
				h^{\prime}(t) = 1+\hat{x}\cdot a^{\prime}(t) &= 1+\, c\, \cos (\theta-\alpha).
			\end{aligned}
		\end{equation*}
Hence $h'$ is a constant depending on $c, \theta$ and $\alpha$.

	Case (i):	If $h^{\prime}(t)>0$, then $\cos (\theta-\alpha) > -1/c$. If $\hat{x}$ is a non-observable direction, that is,  $\xi^{(\hat{x})}_{\max} - \xi^{(\hat{x})}_{\min} < T$, then it holds
		\begin{equation*}
			(t_{\max}-t_{\min})(1+\, c\, \cos (\theta-\alpha))<T.
		\end{equation*}
	Hence, in this case	$\hat{x}$ is a non-observable direction if $-1/\, c\, < \cos (\theta-\alpha))<0$.

Case (ii):	If $h^{\prime}(t)=0$, one can deduce that each direction $\hat{x}$ is non-observable. Note that $\cos (\theta-\alpha) = -1/c$ in such a case.

Case (iii):		If $h^{\prime}(t)<0$, then $\cos (\theta-\alpha) < -1/c$. Consequently, $\hat{x}$ is a non-observable direction only if
		\begin{equation*}
			(t_{\min}-t_{\max})(1+\, c\, \cos (\theta-\alpha))<t_{\max}-t_{\min}.
		\end{equation*}
Therefore,	the direction $\hat{x}$ is  non-observable for $-2/\, c\, < \cos (\theta-\alpha))<-1/c$.

To sum up, we deduce that the non-observable angles should fulfill the relation
		$$-2/c<\cos (\theta - \alpha)<0.$$
This implies that $\theta \in (\alpha +\pi/2,\alpha+3\pi/2)$ for $\, c\, \leq 2$ and  $\theta \in (\alpha - \arccos(-2/c),\alpha -\pi/2) \cup (\alpha +\pi/2,\alpha+ \arccos(-2/c))$ for $\, c\, >2$.
\end{proof}

\vspace{.1in}\textbf{Example 2: An arc in $\R^2$.}\vspace{.1in}

We suppose that the point source is moving along a semi-circle centered at $z=(z_1, z_2)\in \R^2$.
\begin{lemma}\label{arc-ob}
		Set $a(t) := (\cos t+z_1,\sin t+z_2),\,t\in[t_{\min},t_{\max}]$ for some $z \in \mathbb{R}^2$. Suppose that $T=t_{\max}-t_{\min}<2\pi$.
Then the direction $\hat{x}= (\cos \theta, \sin \theta)$ is observable if  $\theta \in [\frac{t_{\max}+t_{\min}}{2},\frac{t_{\max}+t_{\min}}{2}+\pi]$.
		\end{lemma}

		\begin{proof}
We have
		\begin{equation*}
			\begin{aligned}
				h(t) = t+\hat{x}\cdot a(t) &= t + \cos(\theta-t)+\hat{x}\cdot z,\\
				h^{\prime}(t) = 1+\hat{x}\cdot a^{\prime}(t) &= 1+ \sin (\theta-t).
			\end{aligned}
		\end{equation*}

		\noindent It is obvious that $h^{\prime}(t)>0$ for $t\in[t_{\min},t_{\max}]$ such that $t\neq \theta -\pi/2+2n\pi$, $n\in \mathbb{Z}$. Hence,
		\begin{equation*}
		\xi^{(\hat{x})}_{\min}=t_{\min} + \cos (\theta-t_{\min}) +\hat{x}\cdot z , \quad
	\xi^{(\hat{x})}_{\max}=	
		 t_{\max} + \cos (\theta-t_{\max}) +\hat{x}\cdot z.
		\end{equation*}
	For non-observable directions, we have
		\begin{equation*}			
			\cos (\theta-t_{\max})-\cos (\theta-t_{\min})<0,
		\end{equation*}
that is,
		\begin{equation*}
			\sin(\theta - \frac{t_{\max}+t_{\min}}{2})\;\sin\frac{t_{\max}-t_{\min}}{2}<0.
		\end{equation*}
Recalling from the assumption that $0<T=t_{\max}-t_{\min}<2\pi$, one deduces that
		$$\theta \in (\frac{t_{\max}+t_{\min}}{2}+\pi,\frac{t_{\max}+t_{\min}}{2}+2\pi).$$
		Thus, $\hat{x}\in \mathbb{S}^1$ is an observable direction  if $\theta \in [\frac{t_{\max}+t_{\min}}{2},\frac{t_{\max}+t_{\min}}{2}+\pi].$
		\end{proof}

\begin{figure}[!ht]
		\centering
		\scalebox{0.8}{
		\begin{tikzpicture}
		\draw[line width=0.8584cm,color=blue!20] (0,-3) -- (0,4);
		\draw (-0.5,3.5) node [right] {$K_{\Gamma}^{(\hat{x})}$};
		\draw[->] (-3,0) -- (3,0) node[above] {$x_1$} coordinate(x axis);
		\draw[->] (0,-3) -- (0,3) node[right] {$x_2$} coordinate(y axis);
		\foreach \x/\xtext in {-3,-2,-1, 1, 2, 3}
		\draw[xshift=\x cm] (0pt,1pt) -- (0pt,-1pt) node[below] {$\xtext$};
		\foreach \y/\ytext in {-3,-2,-1, 1, 2, 3}
		\draw[yshift=\y cm] (1pt,0pt) -- (-1pt,0pt) node[left] {$\ytext$};
	
		\draw [<-, very thick,red] (-2,-2) arc [ start angle = 225, end angle = 315, radius = 2.828];
		\draw [very thick, dashed] (-0.4292,-3) -- ( -0.4292,4);
		\draw [very thick, dashed] ( 0.4292,-3) -- ( 0.4292,4);
		\draw [very thick, densely dotted] (-2,-3) -- (-2,4);
		\draw [very thick, densely dotted] (2,-3) -- (2,4);
		\draw (-2,0.5) node [left] {$\inf (\hat{x}\cdot \Gamma)$};
		\draw (2,-0.8) node [right] {$\sup (\hat{x}\cdot \Gamma)$};

		\draw (1,-2) node [below] {$\Gamma$};

		\draw [thick,->] (2.1,1) -- (2.5,1);
		\draw (2.5,1) node[right] {$\hat{x} = (1,0)$};
		
		\end{tikzpicture}
		}
		\caption{Illustration of the strip $K_\Gamma^{(\hat{x})}$ (blue area) with $\hat{x} = (1, 0)$. Here the curve $a(t)=2\sqrt{2}(\cos t, -\sin t),\, t\in [\pi/4,3\pi/4]$ denotes the orbit (the red arc) of a point source moving from right to left. There holds $\inf(\hat{x}\cdot\, \Gamma) = -2, \,\sup(\hat{x}\cdot\, \Gamma) = 2,\,\xi_{\min}^{(\hat{x})}=3\pi/4-2,\,\xi_{\max}^{(\hat{x})}=\pi/4+2$. In this case the strip $K_\Gamma^{(\hat{x})}$ is a subset of $\{y\in \R^2: \hat{x}\cdot y\in \hat{x}\cdot \Gamma\}$.   }
		\label{kgamma}
	\end{figure}
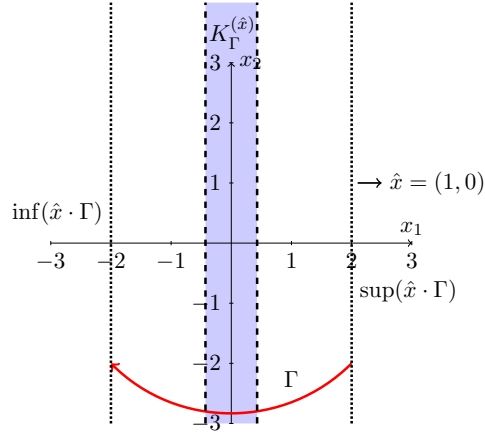

Given the trajectory $\Gamma = \{y: y=a(t),\, t\in [t_{\min}, t_{\max}]\}$, we define
		\begin{equation*}
			\hat{x}\cdot \Gamma \coloneqq \{\tau \in \mathbb{R} : \tau =\hat{x}\cdot y \mbox{ for some }y\in \Gamma\},
		\end{equation*}
which is an interval of $\R$.	
Obviously, the set $\{y\in \R^d: \hat{x}\cdot y\in \hat{x}\cdot \Gamma\}$ denotes the smallest strip containing $\Gamma$ and perpendicular to the direction $\hat{x}$. One can at most expect to recover this strip from the multi-frequency data taken at a single observation direction. 	
	If $\hat{x}$ is an observable direction, we define the strip (see Fig. \ref{kgamma})
	\begin{equation}\label{K}
		K_\Gamma^{(\hat{x})}:=\{y\in \mathbb{R}^d: \xi^{(\hat{x})}_{\min} - t_{\min} \leq \hat{x}\cdot y \leq \xi^{(\hat{x})}_{\max} - t_{\max} \} \subset \mathbb{R}^d.
	\end{equation}
If $h^{\prime}(t)>0$ for $t \in (t_{\min},t_{\max})$, we have
	\begin{equation*}
		K_\Gamma^{(\hat{x})}=\{y\in \mathbb{R}^d: \hat{x}\cdot a(t_{\min}) \leq \hat{x}\cdot y \leq \hat{x}\cdot a(t_{\max}) \}
	\end{equation*}
	which coincides with the strip $\{y\in \R^d: \hat{x}\cdot y\in \hat{x}\cdot \Gamma\}$;
	
	If $h^{\prime}(t)<0$ for $t \in (t_{\min},t_{\max})$, there holds
	\begin{equation*}
		K_\Gamma^{(\hat{x})}=\{y\in \mathbb{R}^d: \hat{x}\cdot a(t_{\max})+T \leq \hat{x}\cdot y \leq \hat{x}\cdot a(t_{\min})-T \},
	\end{equation*}
which is a subset of
$\{y\in \R^d: \hat{x}\cdot y\in \hat{x}\cdot \Gamma\}$; see Lemma \ref{lem3.9} below.

	\begin{lemma}\label{lem3.9}
Let $\hat{x}\in \s^{d-1}$ be an observable direction. We have  $$\hat{x}\cdot y\in \hat{x}\cdot \Gamma\quad\mbox{for all}\quad y\in K_\Gamma^{(\hat{x})}.$$
\end{lemma}

	\begin{proof}
		Suppose that $$\xi^{(\hat{x})}_{\min} = \hat{x}\cdot a(t_1) +t_1 , \quad \xi^{(\hat{x})}_{\max} = \hat{x}\cdot a(t_2) +t_2,\quad\mbox{for some} \quad t_1,t_2\in [t_{\min}, t_{\max}].$$ Therefore,		\ben
			&&\xi^{(\hat{x})}_{\min} - t_{\min} = \hat{x}\cdot a(t_1) +t_1 - t_{\min} \geq \hat{x}\cdot a(t_1) \geq \inf(\hat{x}\cdot \Gamma),\\
			&&	\xi^{(\hat{x})}_{\max} - t_{\max} = \hat{x}\cdot a(t_2) +t_2 - t_{\max} \leq \hat{x}\cdot a(t_2) \leq \sup(\hat{x}\cdot \Gamma).
		\enn
		This implies that for $ y\in K_\Gamma^{(\hat{x})}$,
		\ben
		\hat{x}\cdot y\geq \xi^{(\hat{x})}_{\min} - t_{\min} \geq \inf(\hat{x}\cdot \Gamma),\quad
	\hat{x}\cdot y\leq \xi^{(\hat{x})}_{\max} - t_{\max} \geq \sup(\hat{x}\cdot \Gamma),	
				\enn
which proves $\hat{x}\cdot y\in [\inf(\hat{x}\cdot \Gamma), \sup(\hat{x}\cdot \Gamma)]=\hat{x}\cdot \Gamma$.
	\end{proof}

If $\hat{x}\in \s^{d-1}$ is observable,	we shall prove that the test function $\phi^{(\hat{x})}_y$ lies in the range of $\mathcal{L}^{(\hat{x})}$ if and only if $y \in K_{\Gamma}^{(\hat{x})}$. This together with \eqref{RI}
establishes a computational criterion for imaging $K_\Gamma^{(\hat{x})}$ from the multi-frequency far-field data $u^{\infty}(\hat{x},k)$ with $k \in [k_{\min}, k_{\max}]$.
We also need to discuss non-observable directions.

	\begin{lemma}\label{lem3.4}
		(i) If $\hat{x}$ is non-observable, we have $\phi^{(\hat{x})}_y \notin {\rm Range} (\mathcal{L}^{(\hat{x})})$ for all $y \in \mathbb{R}^d$.\\
		(ii) If $\hat{x}$ is an observable direction, we have $\phi^{(\hat{x})}_y \in {\rm Range} (\mathcal{L}^{(\hat{x})})$ if and only if $y \in K_{\Gamma}^{(\hat{x})}$.
	\end{lemma}

	\begin{proof}
	(i) The first assertion follows directly from Lemma \ref{lem3.5} and the Definition \ref{obd} for non-observable directions.
		
(ii) If $\hat{x}$ is an observable direction, we have $\xi^{(\hat{x})}_{\max} - \xi^{(\hat{x})}_{\min} \geq T$. If $\phi^{(\hat{x})}_y \in {\rm Range} (\mathcal{L}^{(\hat{x})})$, one can find a function $\phi$ satisfying $\phi^{(\hat{x})}_y = \mathcal{L}^{(\hat{x})} \phi$. Then their support sets  must fulfill the relation supp$(\mathcal{F}^{-1}\phi^{(\hat{x})}_y)=$ supp$(\mathcal{F}^{-1}\mathcal{L}^{(\hat{x})} \phi) \subset [\xi^{(\hat{x})}_{\min}, \xi^{(\hat{x})}_{\max}]$ by Lemma \ref{lem:interrange}. Using Lemma \ref{INVF-PHI} yields
		$$[\hat{x}\cdot y + t_{\min}, \hat{x}\cdot y + t_{\max}] \subset [\xi^{(\hat{x})}_{\min}, \xi^{(\hat{x})}_{\max}].$$
		Hence, $\hat{x}\cdot y + t_{\min} \geq \xi^{(\hat{x})}_{\min}$ and $\hat{x}\cdot y + t_{\max} \leq \xi^{(\hat{x})}_{\max}$, leading to		$$\xi^{(\hat{x})}_{\min} - t_{\min} \leq \hat{x}\cdot y \leq \xi^{(\hat{x})}_{\max} - t_{\max}.$$
This proves $y \in K_{\Gamma}^{(\hat{x})}$.

		On the other hand, if $y \in K_{\Gamma}^{(\hat{x})}$, we have
		\begin{equation*}
			[\hat{x}\cdot y +t_{\min}, \hat{x}\cdot y +t_{\max}] \subset [\xi^{(\hat{x})}_{\min},\xi^{(\hat{x})}_{\max}].
		\end{equation*}
Setting $$\psi(t) := \frac{e^{ik\hat{x}\cdot (a(t)-y)}}{|t_{\max}-t_{\min}|}\in L^2(t_{\min},t_{\max}),$$ we find $\phi^{(\hat{x})}_y(k) = (\mathcal{L}^{(\hat{x})}\psi) (k)$. Therefore, $\phi^{(\hat{x})}_y(k) \in \mbox{ Range}(\mathcal{L}^{(\hat{x})})$.
		
	\end{proof}

	\section{Indicator functions and uniqueness}\label{IdF}
If $\hat{x}$ is an observable direction, we know from
 Lemma \ref{lem3.4} that the test functions $\phi^{(\hat{x})}_y$ can be utilized to characterize $K_{\Gamma}^{(\hat{x})}$ through \eqref{RI}. Hence, we define the indicator function
	\begin{equation}\label{indicator}
		W^{(\hat{x})}(y) \coloneqq \left[\sum_{n=1}^\infty\frac{|\langle \phi^{(\hat{x})}_{y}, \psi_n^{(\hat{x})} \rangle|_{L^2(0, K)}^2}{ |\lambda_n^{(\hat{x})}|}\right]^{-1}, \qquad y\in \mathbb{R}^d.
	\end{equation}
 Combining Theorem \ref{range}, Lemma \ref{lem3.4} and Picard theorem, we obtain.
	
	\begin{theorem}\label{Th:factorization}
	
		If $\hat{x}$ is an observable direction, it holds that
		\ben
		W^{(\hat{x})}(y)=\left\{\begin{array}{lll}
		0 &&\quad\mbox{if}\quad y\notin K_\Gamma^{(\hat{x})},\\
		\mbox{finite positive number} &&\quad\mbox{if}\quad y\in K_\Gamma^{(\hat{x})}.		
		\end{array}\right.
		\enn
		If $\hat{x}$ is non-observable, we have $W^{(\hat{x})}(y)=0$ for all $y\in \mathbb{R}^d$.

	\end{theorem}
	
	Hence, for observable directions the values of $W^{(\hat{x})}$ in the strip $K_\Gamma^{(\hat{x})}$ should be relatively bigger than those elsewhere. The values of $W^{(\hat{x})}$ vanished identically in $\R^d$ if $\hat{x}$ is non-observable.
In the case of sparse observable directions $\{\hat{x}_j: j=1,2,\cdots, M\}$, we shall make use of the following indicator function:
	\begin{equation}\label{W}
		W(y)= \left[\sum_{j=1}^M \frac{1}{W^{(\hat{x}_j)}(y)}\right]^{-1}=		
		\left[\sum_{j=1}^M\sum_{n=1}^\infty\frac{|\langle \phi^{(\hat{x}_j)}_{y}, \psi_n^{(\hat{x}_j)} \rangle|_{L^2(0, K)}^2}{ |\lambda_n^{(\hat{x}_j)}|}\right]^{-1}, \qquad y\in \mathbb{R}^d.
	\end{equation}
	Define the $\Theta$-convex domain of $\Gamma$ associated with the observable directions $\{\hat{x}_j: j=1,2,\cdots, M\}$ as
	\be\label{Theta}
		\Theta_\Gamma \coloneqq \bigcap\limits_{j=1,2,\cdots, M} K_\Gamma^{(\hat{x}_j)}.
	\en
We can reconstruct $\Theta_\Gamma$ from 	
the multi-frequency far-field data measured at sparse observable directions.
\begin{theorem}\label{TH4.2}
		It holds that $0<W(y) <+\infty$ if $y \in \Theta_\Gamma$ and
		$W(y)=0$ if $y \notin \Theta_\Gamma$.
\end{theorem}

	\begin{proof}
		If $y\in \Theta_\Gamma$, it means that $y\in K_\Gamma^{(\hat{x}_j)}$ for $j=1,2,\cdots,M$. By Theorem \ref{Th:factorization},
		\begin{equation}\label{wjy}
			\sum\limits_{n=1}^{\infty} \frac{|\langle \phi^{(\hat{x}_j)}_{y}, \psi_n^{(\hat{x}_j)} \rangle|_{L^2(0, K)}^2}{ |\lambda_n^{(\hat{x}_j)}|} < +\infty \quad\mbox{for all}\quad j = 1,2,\cdots,M.
		\end{equation}
		Then the finite sum over the index $j$ must fulfill the relation $0<W(y)<+\infty$.

		If $y\notin \Theta_\Gamma$, we may suppose without loss of generality that $y\notin K_\Gamma^{(\hat{x}_1)}$.
 By Theorem \ref{Th:factorization},
 		\begin{equation*}
		[W^{(\hat{x}_1)}(y)]^{-1}=	\sum\limits_{j=1}^{M}\sum\limits_{n=1}^{\infty} \frac{|\langle \phi^{(\hat{x}_1)}_{y}, \psi_n^{(\hat{x}_j)} \rangle|_{L^2(0, K)}^2}{ |\lambda_n^{(\hat{x}_1)}|} = \infty.
		\end{equation*}
Together with the definition of $W$, this gives
\ben
W(y)<\left[\sum\limits_{n=1}^{\infty} \frac{|\langle \phi^{(\hat{x}_1)}_{y}, \psi_n^{(\hat{x}_j)} \rangle|_{L^2(0, K)}^2}{ |\lambda_n^{(\hat{x}_1)}|}   \right]^{-1}=0.
\enn

	\end{proof}
Consequently, we arrive at the following uniqueness results, which seem unknown in the literature.
\begin{theorem}\label{TH4.3} Denote by $\Gamma=\{a(t): t\in[t_{\min}, t_{\max}]\}$ the trajectory of a moving point source where $a\in C^1[t_{\min}, t_{\max}]$.

(i) The $\Theta$-convex domain of $\Gamma$ associated with all observable directions $\hat{x}\in \s^{d-1}$ (see \eqref{Theta}) can be uniquely determined by the multi-frequency data $\{u^{\infty}(\hat{x}, k): \hat{x}\in \s^{d-1}, k\in(k_{\min}, k_{\max})\}$.

(ii) Let $\hat{x}\in\s^{d-1}$ be an arbitrarily fixed observable direction. Then the strip $K_\Gamma^{(\hat{x})}$ (see \eqref{K}) can be uniquely determined by the multi-frequency data $\{u^{\infty}(\hat{x}, k):  k\in(k_{\min}, k_{\max})\}$. In particular, the strip $\{y\in \R^{d}: \hat{x}\cdot y\in \hat{x}\cdot \Gamma\}$ can be uniquely recovered if $1+\hat{x}\cdot a'(t)>0$ in $[t_{\min}, t_{\max}]$.
\end{theorem}	

\begin{remark} Physically,
the condition $1+\hat{x}\cdot a'(t)>0$ in the second assertion of Theorem \ref{TH4.3} means that the function $h(t)=t+\hat{x}\cdot a(t)$ is monotonically increasing in $[t_{\min}, t_{\max}]$. It can be fulfilled if the velocity of the moving source is less than the propagating speed of waves, i.e., $|a'(t)|<1$. Note that the acoustical speed in the background medium has been normalized to be one.
\end{remark}

The second assertion of Theorem \ref{TH4.3} answers the question what kind of information can be extracted from the multi-frequency data measured at a single observable direction. Unfortunately, we do not know whether an observation direction is observable or not, if there is no a priori information on the orbit function.

	\section{Numerical experements in $\R^d$ ($d=2,3$)}\label{num}
    In this section, we carry out a couple of numerical experiments to validate our algorithm in both two and three dimensions.  In practice, the time-domain data
 should be Fourier transformed to the multi-frequency data and the near-field version of our algorithm should be implemented. To simply the numerical procedures for simulating,  we shall carry out computational tests in the frequency domain only. Our aim is to get information of the trajectory of a moving point source from multi-frequency far-field data taken at a single or multiple observation directions.

Suppose that the wave-number-dependent source term $f(x, k)$ is given by (\ref{sourcef}). Then the far-field pattern can be synthetized by (\ref{u-infty}), i.e.,

\begin{equation}
\begin{aligned}
w^\infty(\hat{x}, k)
=
\int_{t_{\min}}^{t_{\max}} e^{-ik(\hat{x}\cdot a(t)+t)} \, dt,\quad \hat{x}\in \R^d\, (d=2,3),\; k\in(k_{\min}, k_{\max}).
\end{aligned}
\end{equation}
In all our numerical examples below, we set $k_{\min}=0$ for simplicity. The bandwidth can be extended from $(0,k_{max})$ to $(-k_{\max}, k_{\max})$ by $w^{\infty}(\hat{x}, -k)=\overline{w^{\infty}(\hat{x}, k)}$. Then, one deduces from these new measurement data with $k_{min}=-k_{max}$ that $\kappa=0$ and $K=k_{max}$. Thus, the far field operator \eqref{FarO} becomes
    \begin{equation} \label{far-field0}
			(\mathcal{F}^{(\hat{x})}\phi)(\tau) = \int_0^{k_{max}} w^{\infty}(\hat{x},  \tau - s)\,\phi(s)\,ds, \quad \tau\in(0, k_{\max}).
    \end{equation}
    Discretize the frequency interval $(0,k_{\max})$ with $$k_n=(n-0.5)\Delta k, \quad \Delta k:=\frac{k_{\max}}{N}, \quad n=1,2,\cdots,N.$$
    We adopt {$2N-1$ samples $w^{\infty}(\hat x, k_n), n=1,2,\cdots,N$ and $w^{\infty}(\hat x, -k_n), n=1,2,\cdots,N-1$}, of the far field and apply the midpoint rule to approximate the integral in (\ref{far-field0}). Then it follows that
    \begin{equation}
    (\mathcal{F}^{(\hat x)}\phi)(\tau_n) \approx \sum_{m=1}^{N} w^{\infty}(\hat x, \tau_n-s_m)\phi(s_m)\Delta k,
    \end{equation}
    where $\tau_n:=n\Delta k$ and $s_m:=(m-0.5)\Delta k$, $n,m=1,2,\cdots,N$.
    Consequently, a discrete approximation of the far field operator $\mathcal{F}^{(\hat x)}$ is given by the Toeplitz matrix
\be \label{matF}
F^{(\hat x)}:=\begin{pmatrix}
w^{\infty}(\hat{x},k_1) & \overline{w^{\infty}(\hat{x},k_1)} & \cdots & \overline{w^{\infty}(\hat{x},k_{N-2})}  & \overline{w^{\infty}(\hat{x},k_{N-1})}  \\
w^{\infty}(\hat{x},k_2) & w^{\infty}(\hat{x},k_1) & \cdots & \overline{w^{\infty}(\hat{x},k_{N-3})} &\overline{w^{\infty}(\hat{x},k_{N-2})}   \\
\vdots & \vdots  &  &\vdots &\vdots \\
w^{\infty}(\hat{x},k_{N-1}) & w^{\infty}(\hat{x},k_{N-2}) &  \cdots & w^{\infty}(\hat{x},k_1) & \overline{w^{\infty}(\hat{x},k_1)}\\

w^{\infty}(\hat{x},k_N) & w^{\infty}(\hat{x},k_{N-1}) &  \cdots & w^{\infty}(\hat{x},k_2) & w^{\infty}(\hat{x},k_1)\\
\end{pmatrix} \Delta k \in \C^{N}\times \C^{N}
\en
where $\overline{w^{\infty}(\hat{x},k_n)}=w^{\infty}(\hat{x},-k_n)$, $n=1,\cdots,N-1$.

\noindent Similarly, we discretize the test function $\phi_{y}^{(\hat x)}$ from (\ref{testfunc}) by the vector
\be \label{testn}
\phi_y^{(\hat{x})}:= \left(\frac{i}{T\tau_1} (e^{-i \tau_1 t_{\max}}-e^{-i \tau_1 t_{\min}})\,e^{-i\tau_1 \hat{x}\cdot y}, \;\cdots,\; \frac{i}{T\tau_n} (e^{-i \tau_n t_{\max}}-e^{-i \tau_n t_{\min}})\,e^{-i\tau_n \hat{x}\cdot y}\right) \in \C^N,
\en
where $T=t_{\max}-t_{\min}$.
\noindent Denoting by  $\left\{ ( {\tilde \lambda^{(\hat x)}_n}, \psi^{(\hat x)}_n): n=1,2,\cdots,N \right\}$ an eigen-system of the matrix $F^{(\hat x)}$ (\ref{matF}), {then one deduces that  an eigen-system of the matrix $(F^{(\hat x)})_\#:= |Re(F^{(\hat x)})|+|Im(F^{(\hat x)})|$ is $\left\{ ( \lambda^{(\hat x)}_n, \psi^{(\hat x)}_n): n=1,2,\cdots,N \right\}$ , where $ \lambda^{(\hat x)}_n:=|Re (\tilde \lambda^{(\hat x)}_n)| +|Im (\tilde \lambda^{(\hat x)}_n)|$. We } approximate the indicator function $W^{(\hat x)}$ of (\ref{indicator}) by
\ben \label{ind-n}
W^{(\hat{x})}(y):=\left[\sum_{n=1}^N\frac{\left| \phi^{(\hat{x})}_{y} \cdot \overline{\psi_n^{(\hat{x})} }\right|^2}{ |\lambda_n^{(\hat{x})}|}\right]^{-1}, \quad y\in \R^d,
\enn
where $\cdot$ denotes the inner product in $\R^N$.
Accordingly, a plot of $W^{(\hat{x})}(y)$, should yield a visualization of the strip $K_\Gamma^{(\hat{x})}$, which contains information on the source trajectory if $\hat{x}\in \s^{d-1}$ is an observable direction. In the following numerical examples, the frequency band is taken as $(0, 3\pi)$ with $k_{\max}=3\pi$, $N=18$ and $\Delta k=\pi/6$.

In the following figures, the exact trajectory of a moving source is plotted with yellow sold lines. In two dimensions we shall image the trajectory of moving point sources represented by a straight line, an arc or a piecewise linear curve, using the far-field data of one and sparse observation directions. In three dimensions, the recovery of a straight line segment is examined with the data measured at a single direction only.


\subsection{A single observation direction}

\textbf{Example 1: A straight line segment in $\R^2$} \vspace{.1in}

We consider  the same straight line segment from  Example 1 in Section \ref{RangeLx}. The following two cases are studied.

\vspace{.1in}\textbf{Case 1} $c=1$, $\alpha=\pi/2$ and $t\in[1,3]$.

\begin{figure}[htb]
\centering
\subfigure[$\theta=0$ ]{
\includegraphics[scale=0.22]{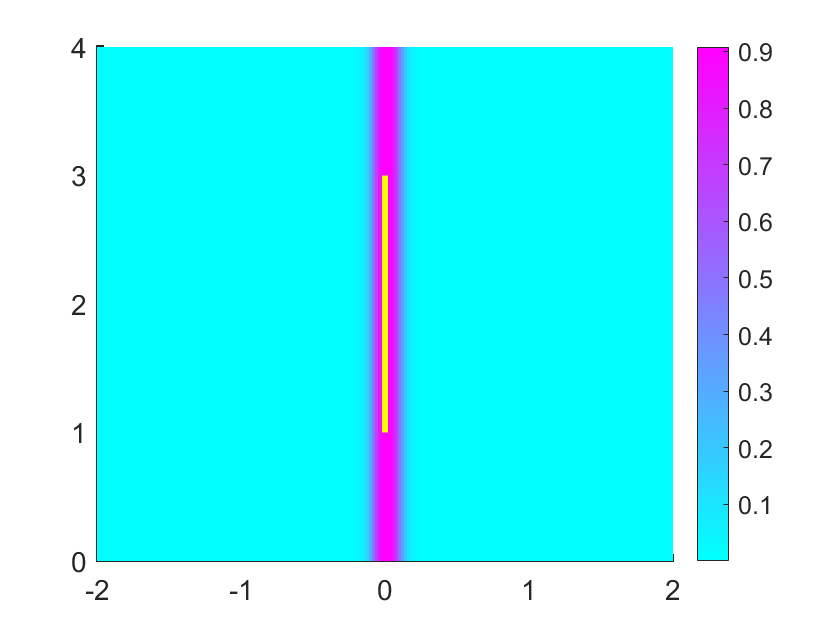}

}
\subfigure[$\theta=\pi/6$ ]{
\includegraphics[scale=0.22]{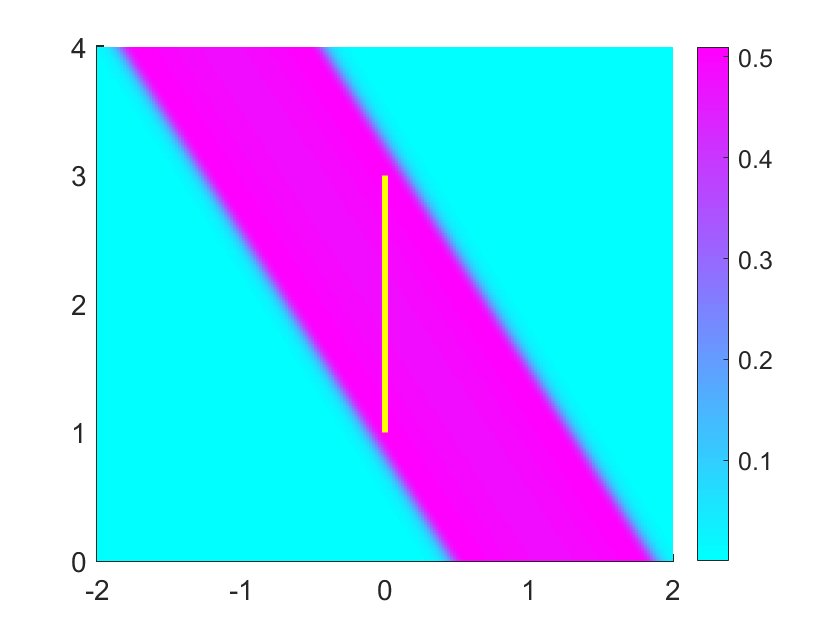}

}
\subfigure[$\theta=\pi/3$]{
\includegraphics[scale=0.22]{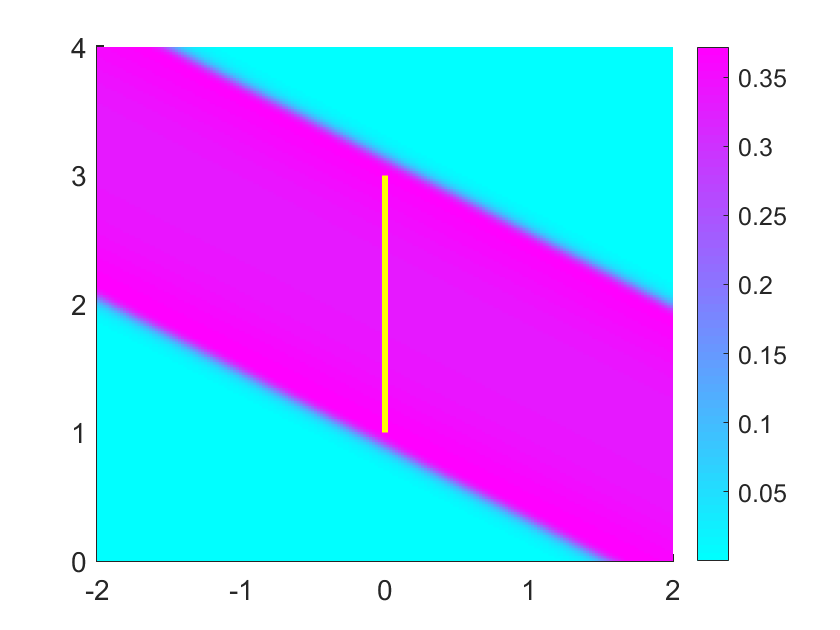}

}

\subfigure[$\theta=\pi/2$]{

\includegraphics[scale=0.22]{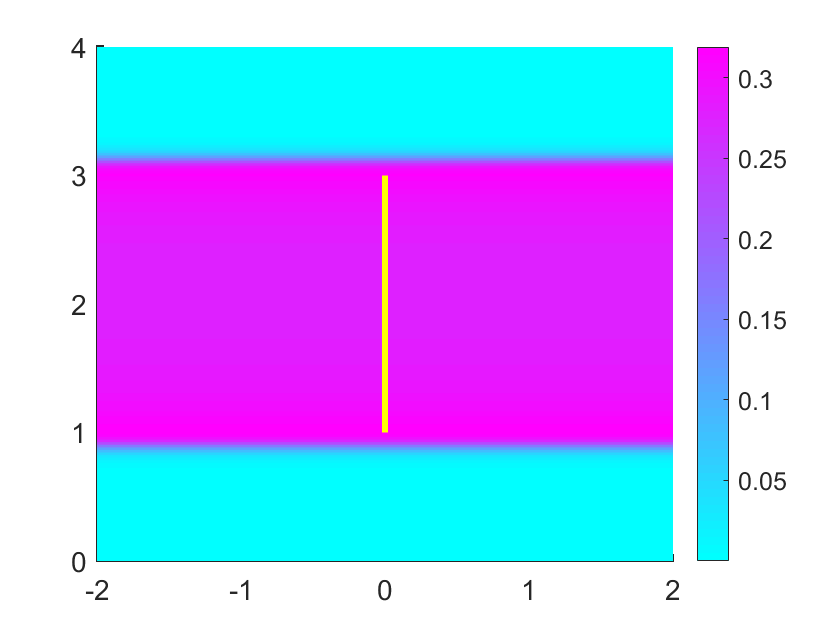}

}
\subfigure[$\theta=2\pi/3$ ]{
\includegraphics[scale=0.22]{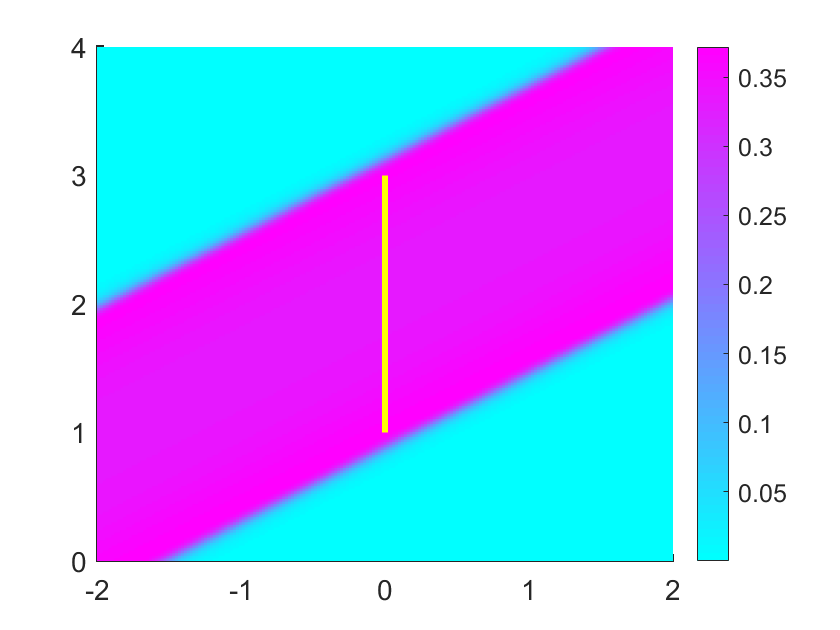}

}
\subfigure[$\theta=5\pi/6$]{
\includegraphics[scale=0.22]{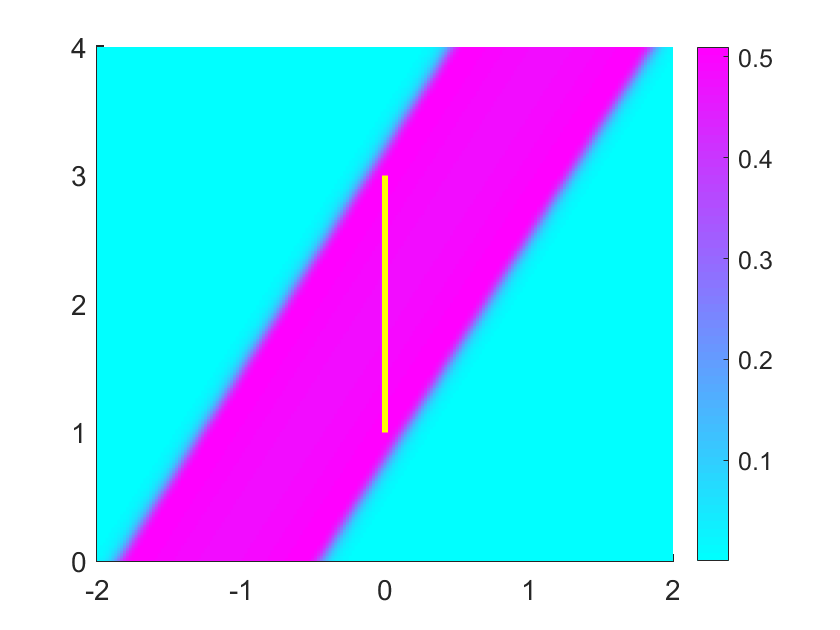}

}
\caption{Reconstruction from a single observable direction $\hat{x}=(\cos \theta, \sin \theta)$ with $\theta\in[0,\pi]$  for a straight line segment $a(t)=(0,t)$ with $t\in[1,3]$.
} \label{fig:line1}
\end{figure}

\begin{figure}[htb] \label{str-2-unob}
\centering
\subfigure[$\theta=9\pi/8$]{

\includegraphics[scale=0.22]{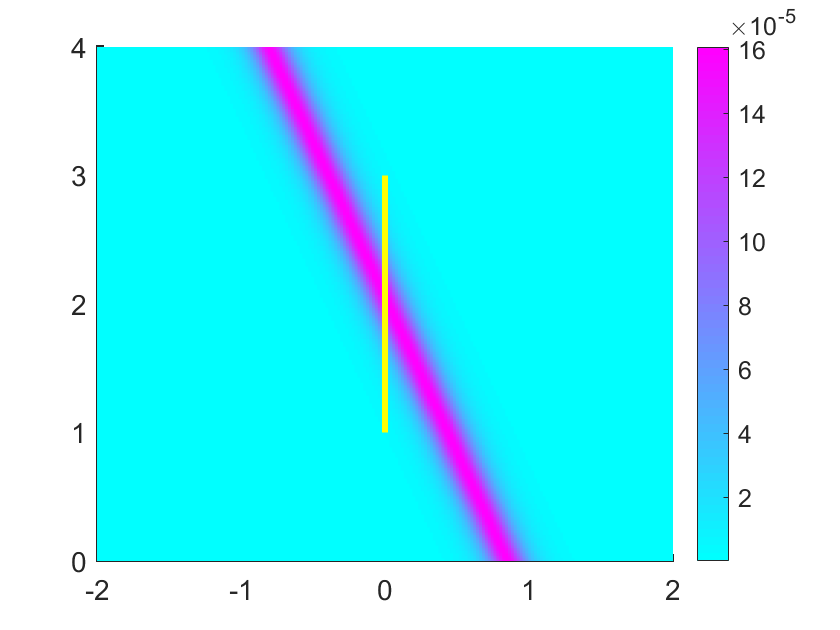}

}
\subfigure[$\theta=10\pi/8$ ]{
\includegraphics[scale=0.22]{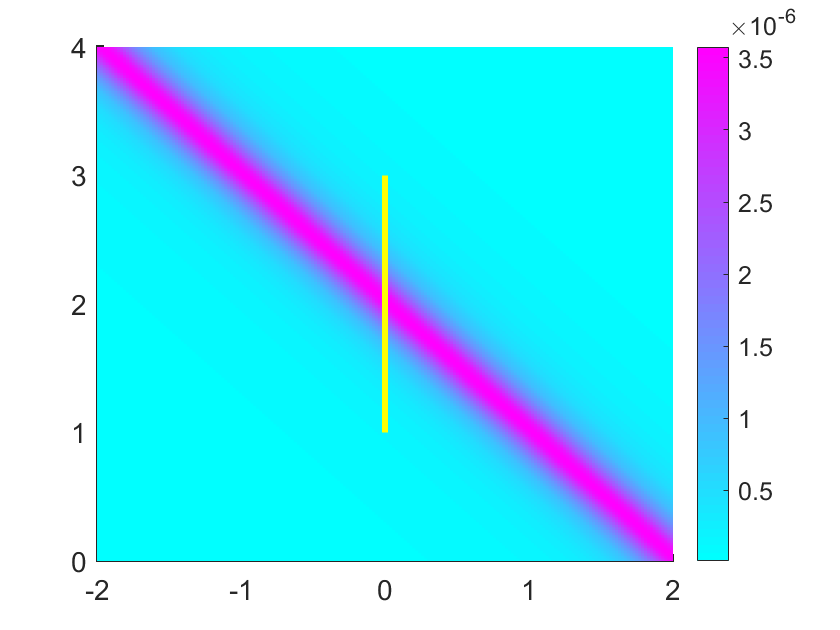}

}
\subfigure[$\theta=11\pi/8$]{
\includegraphics[scale=0.22]{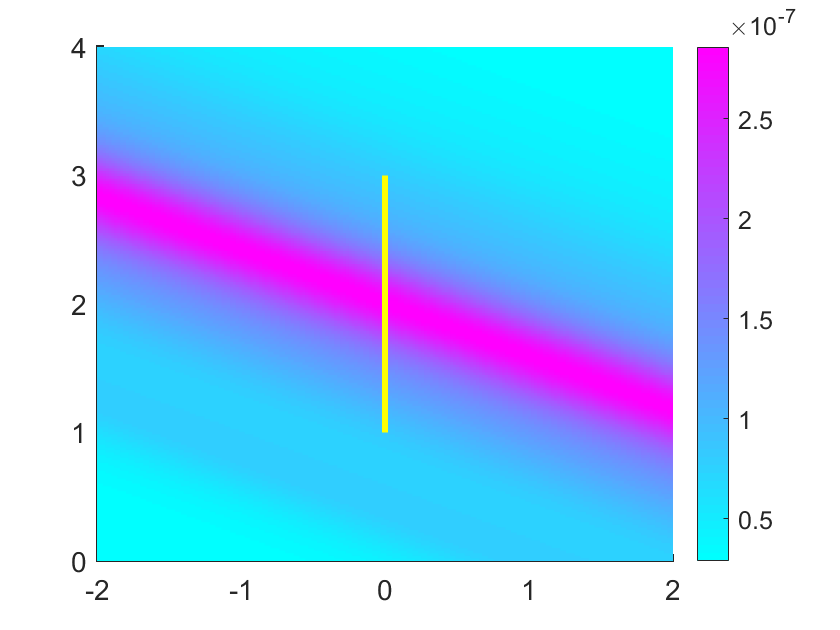}

}

\subfigure[$\theta=13\pi/8$]{

\includegraphics[scale=0.22]{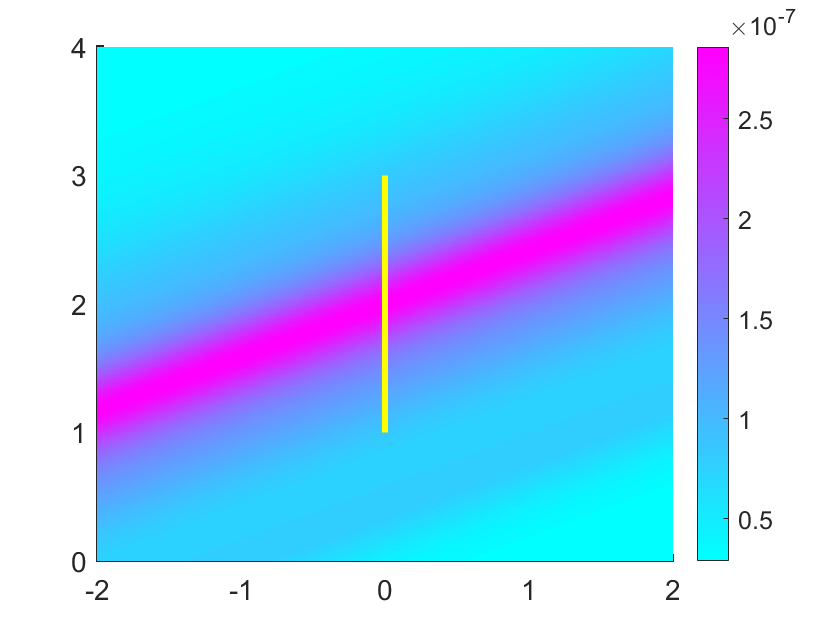}

}
\subfigure[$\theta=14\pi/8$ ]{
\includegraphics[scale=0.22]{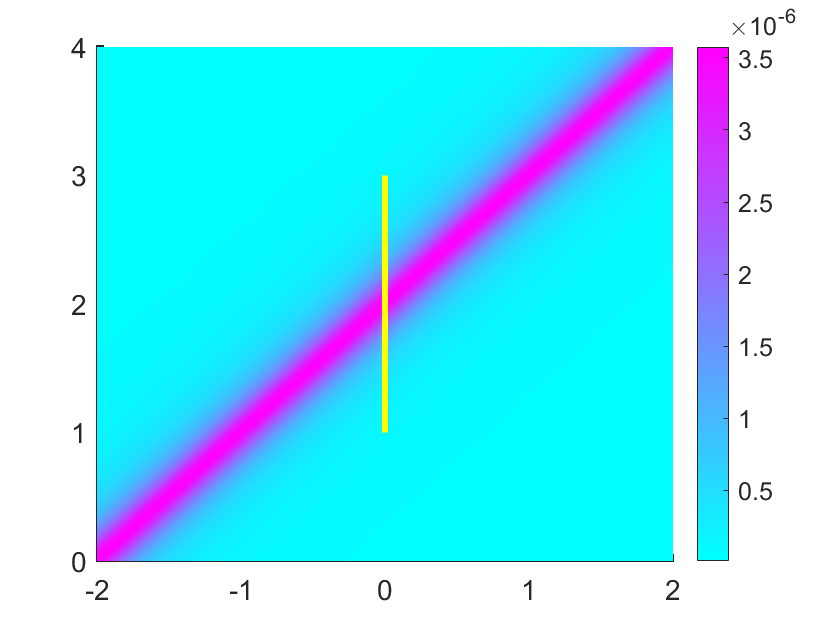}

}
\subfigure[$\theta=15\pi/8$]{
\includegraphics[scale=0.22]{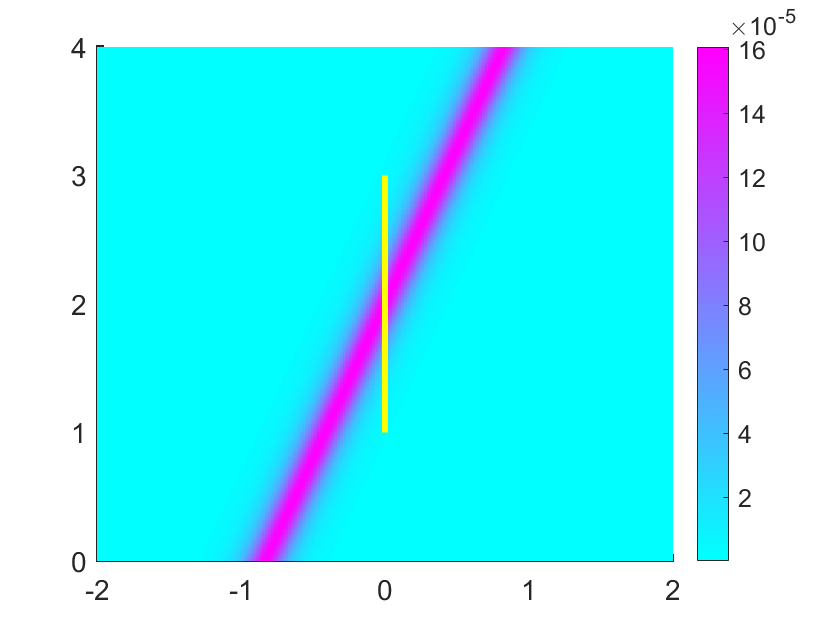}

}
\caption{Reconstruction from a single non-observable direction $\hat{x}=(\cos\theta, \sin\theta)$ with $\theta\in(\pi,2\pi)$ for a straight line segment $a(t)=(0,t)$ with $t\in[1,3]$. } \label{fig:line2}
\end{figure}

 In this case the trajectory of the moving source is  $a(t)=(0,t)$ for $t\in[1,3]$. Choose the search domain as a square of the form $[-2,2]\times[0,4]$. By Lemma \ref{line-ob}, the non-observable directions are $\hat{x}=(\cos \theta, \sin\theta)$ with $\theta\in(\pi,2\pi)$ and the observable directions $\hat{x}=(\cos \theta, \sin \theta)$ with $\theta\in[0,\pi]$. Numerical results are presented in Figs. \ref{fig:line1} and \ref{fig:line2}.

 The observable angles are taken as $\theta=0$, $\pi/6$, $\pi/3$, $\pi/2$, $2\pi/3$ and  $5\pi/6$ in Fig. \ref{fig:line1}. By Lemma \ref{line-ob}, we know  $K_{\Gamma}^{(\hat x)}= \{y\in\R^2:\inf (\hat x \cdot \Gamma) \leq \hat x\cdot y\leq \sup( \hat x \cdot \Gamma)\}$  for  all observable directions, because $h'(t)=1+\cos(\theta-\pi/2)>0$  with $t\in [1,3]$.
 In Fig. \ref{fig:line1}, the trajectory of the moving source is nicely located in the smallest strip $K_{\Gamma}^{(\hat x)}$ perpendicular to the observation direction just as our theoretical results predict.  The numerical results match well with our theoretical analysis.

Observation directions at the angles $\theta=9\pi/8$, $10\pi/8$, $11\pi/8$, $13\pi/8$, $14\pi/8$ and  $15\pi/8$ are non-observable. The numerical results in Fig.\ref{fig:line2} show that the indicator values are all much smaller than $10^{-4}$, which are in good consistent with the results of Theorem \ref{Th:factorization}. Hence, we can not reconstruct the smallest strip containing the trajectory of the moving source. It is very interesting to conclude from Fig.\ref{fig:line2} that, even at a non-observable direction, partial information on the trajectory can still be recovered by our indicator function: the maximum points of $W^{(\hat{x})}$ are degenerated into a straight line perpendicular to $\hat{x}$ and passing through the middle point of the trajectory. However, this phenomenon needs to be further investigated.

\vspace{.1in}\textbf{Case 2 }
$c=4$,  $\alpha=\pi/4$ and $t\in[1,2]$.

In this case $a(t)=(2\sqrt 2 t, 2\sqrt2 t)$ represents a diagonal line segment. The search domain is taken as $[-2,5]\times[-2,5]$.  The observable directions are  $\hat{x}$ with $\theta\in[0,3\pi/4]\cup [7\pi/4, 8\pi/4]\cup[11\pi/12, 19\pi/12]$
 and non-observable directions are $\hat{x}$ with $\theta\in (3\pi/4,11\pi/12)\cup(19\pi/12, 7\pi/4)$. By the proof of Lemma \ref{line-ob},
 $h'(t)=1+4\cos(\theta-\pi/4)>0$ for
 observable angles $\theta \in [0,3\pi/4]\cup[7\pi/4,2\pi)$ and
 $h'(t)<0$ for
 $\theta \in [11\pi/12, 19\pi/12]$.

\begin{figure}[htb]
\centering
\subfigure[$\theta=0\pi/8$]{

\includegraphics[scale=0.22]{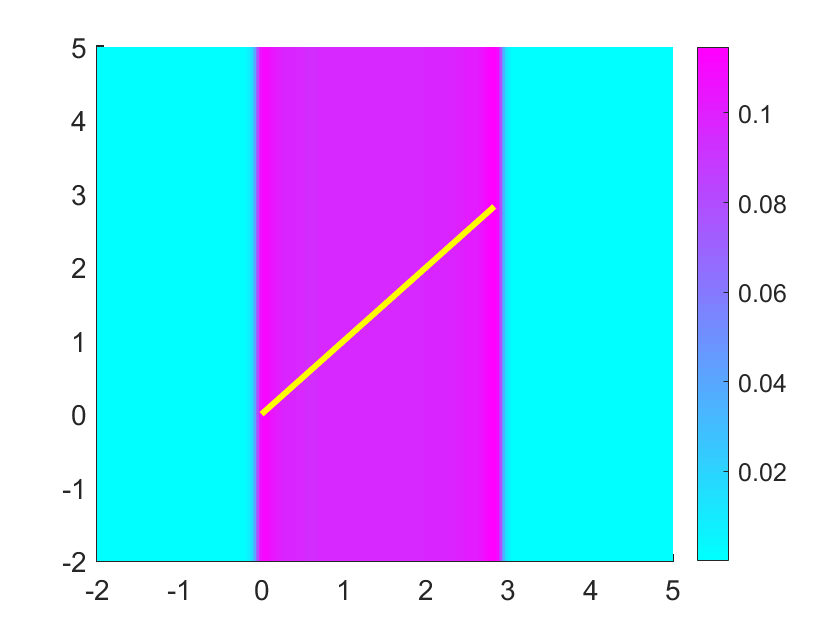}

}
\subfigure[$\theta=3\pi/8$ ]{
\includegraphics[scale=0.22]{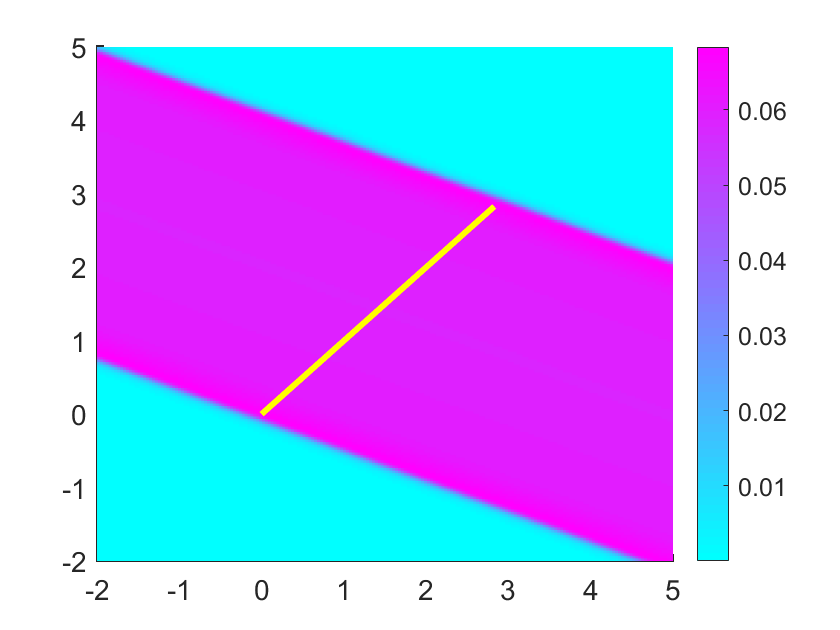}

}
\subfigure[$\theta=4\pi/8$]{
\includegraphics[scale=0.22]{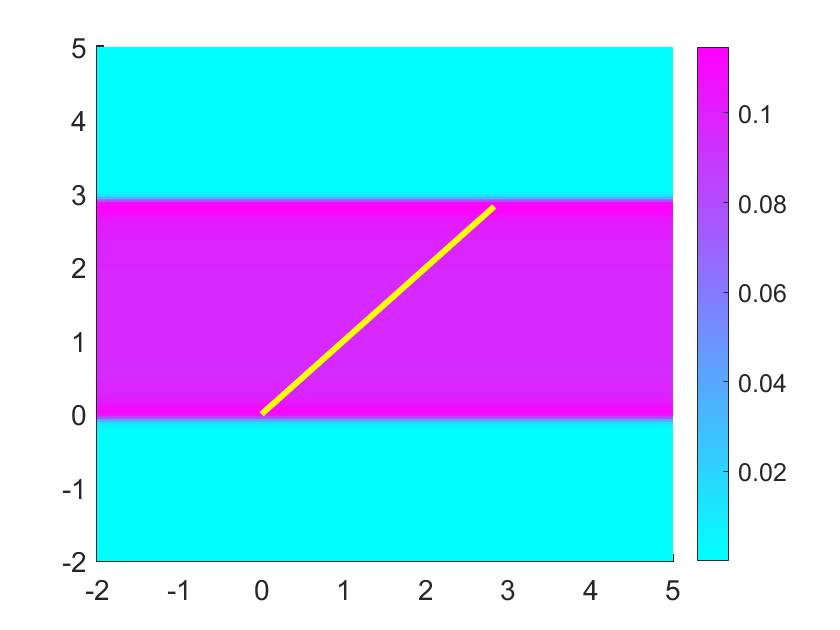}

}
\subfigure[$\theta=5\pi/8$]{

\includegraphics[scale=0.22]{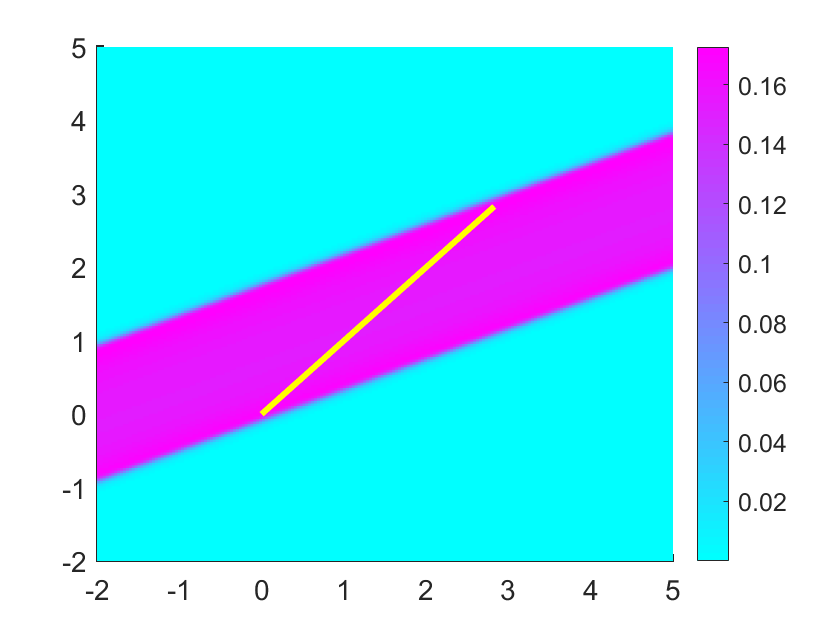}

}
\subfigure[$\theta=6\pi/8$ ]{
\includegraphics[scale=0.22]{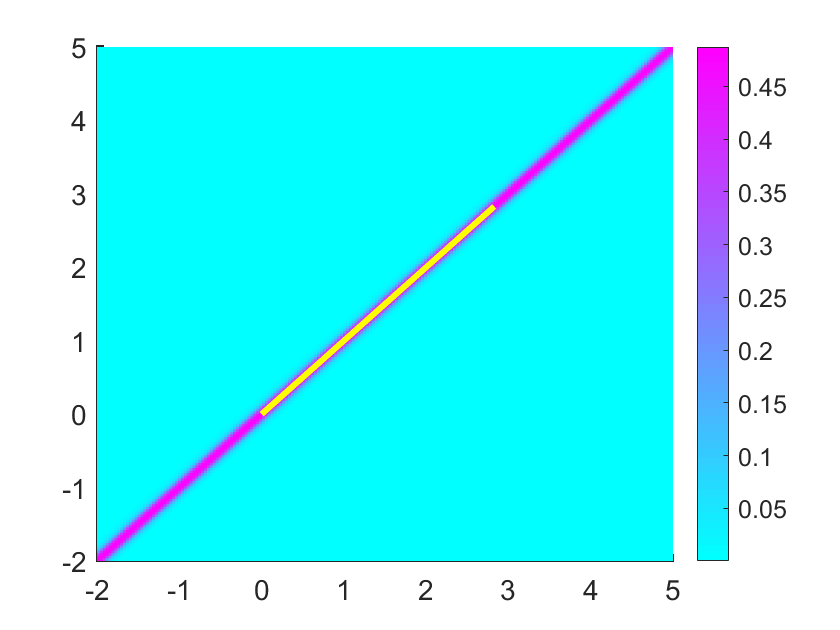}

}
\subfigure[$\theta=15\pi/8$]{
\includegraphics[scale=0.22]{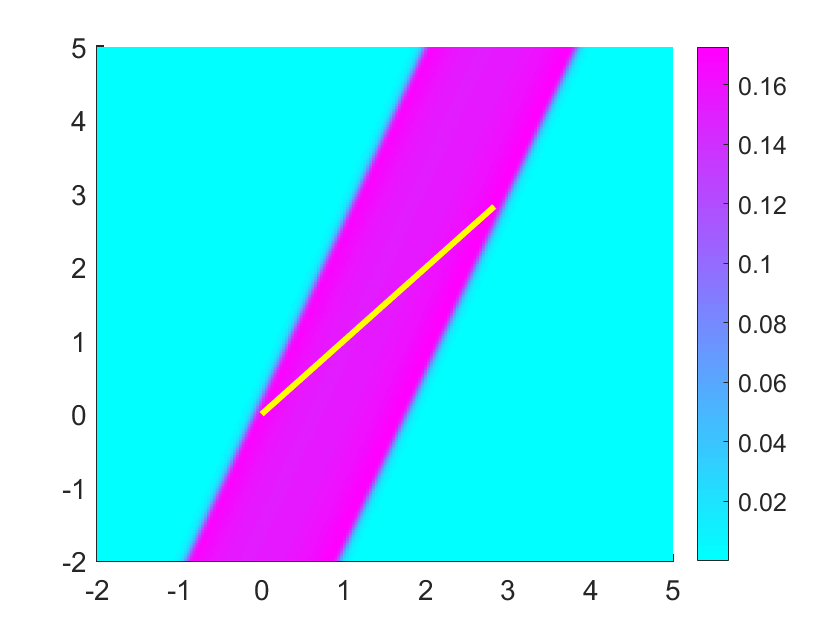}

}
\caption{Reconstruction from a single observable direction $\hat x=(\cos\theta, \sin\theta)$  with $\theta\in[0,3\pi/4]\cup [7\pi/4, 2\pi)$ for a straight line segment $a(t)=(2\sqrt 2t, 2\sqrt 2t)$ with $t\in[1,2]$. Since $h'(t)>0$, the strip $K_\Gamma^{(\hat{x})}$ coincides with $\{y\in \R^2: \hat{x}\cdot y\in \hat{x}\cdot \Gamma\}$.} \label{fig:xie-line1}
\end{figure}

In Fig.\ref{fig:xie-line1}, we take different observable angles $\theta \in [0,3\pi/4]\cup[7\pi/4,2\pi)$. Since $h'(t)>0$,  the trajectory of the moving source can be completely covered by the smallest strip perpendicular the observation direction. The numerical examples indeed show that $K_{\Gamma}^{(\hat x)}= \{y\in\R^2:\inf (\hat x \cdot \Gamma) \leq \hat x\cdot y\leq \sup( \hat x \cdot \Gamma)\}$.

\begin{figure}[htb]
\centering
\subfigure[$\theta=11\pi/12$]{

\includegraphics[scale=0.22]{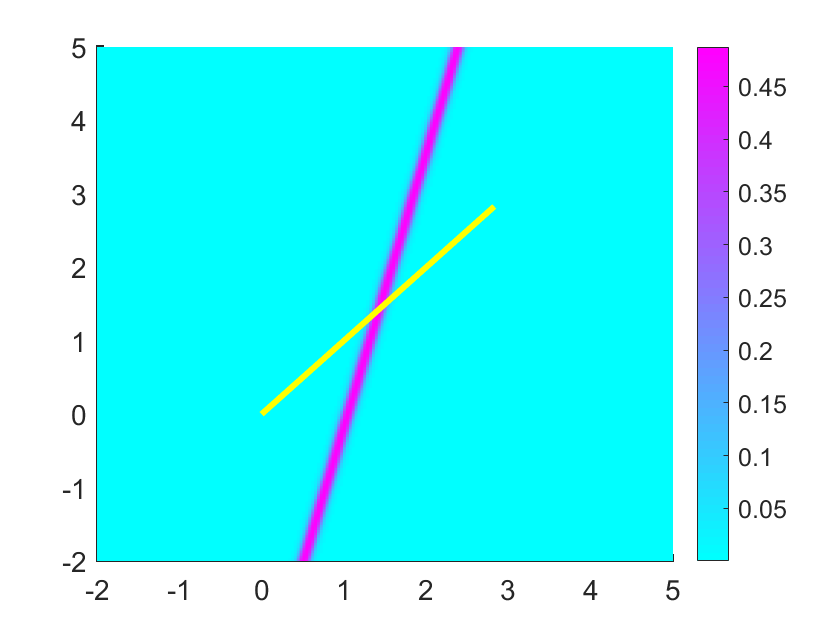}

}
\subfigure[$\theta=8\pi/8$ ]{
\includegraphics[scale=0.22]{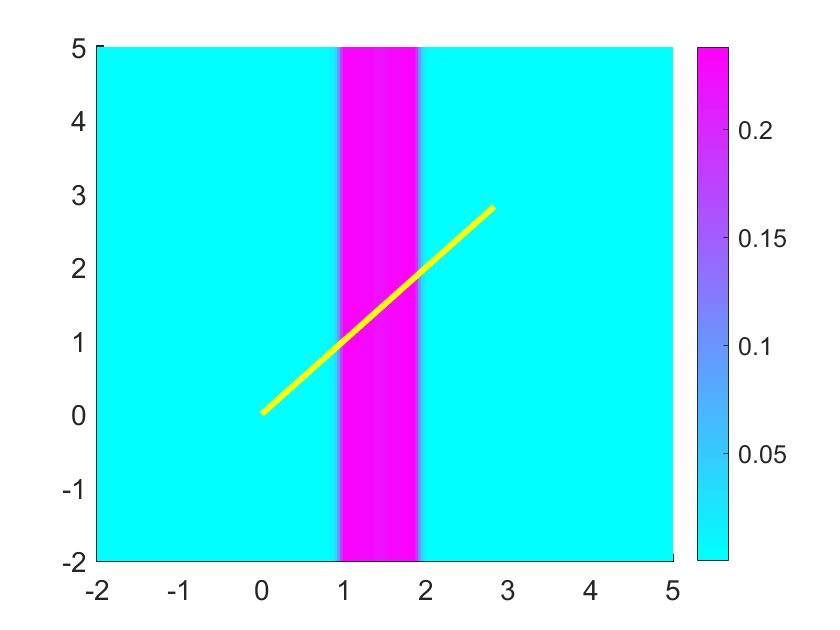}

}
\subfigure[$\theta=9\pi/8$]{
\includegraphics[scale=0.22]{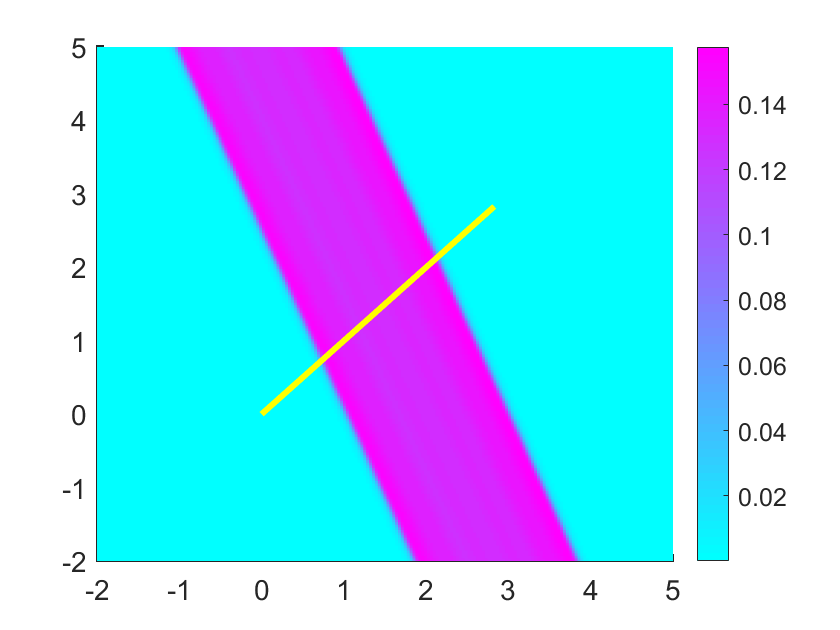}

}
\subfigure[$\theta=10\pi/8$]{

\includegraphics[scale=0.22]{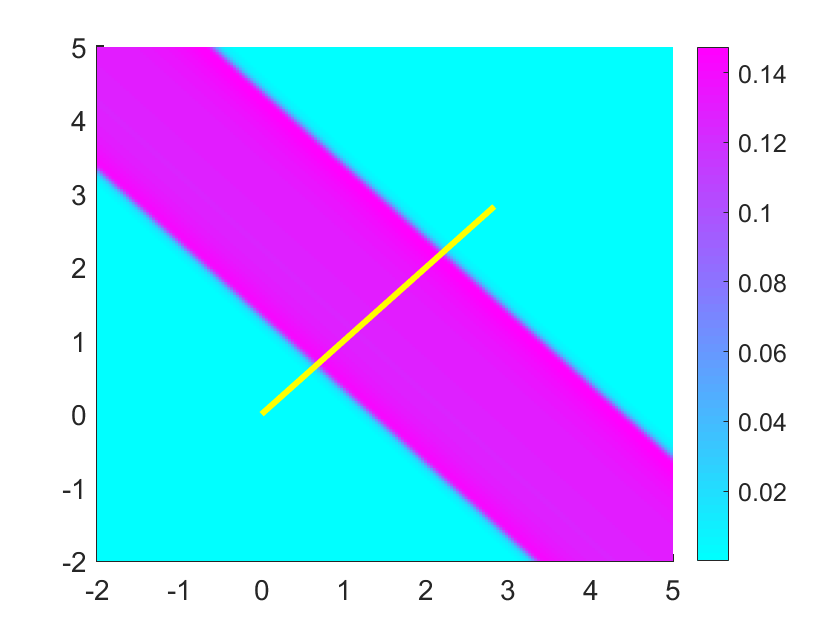}

}
\subfigure[$\theta=11\pi/8$ ]{
\includegraphics[scale=0.22]{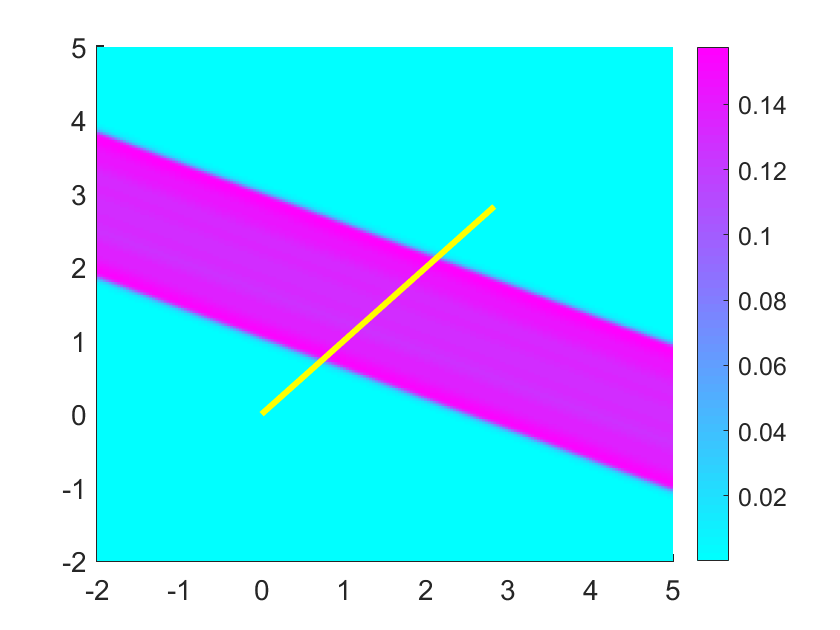}

}
\subfigure[$\theta=12\pi/8$]{
\includegraphics[scale=0.22]{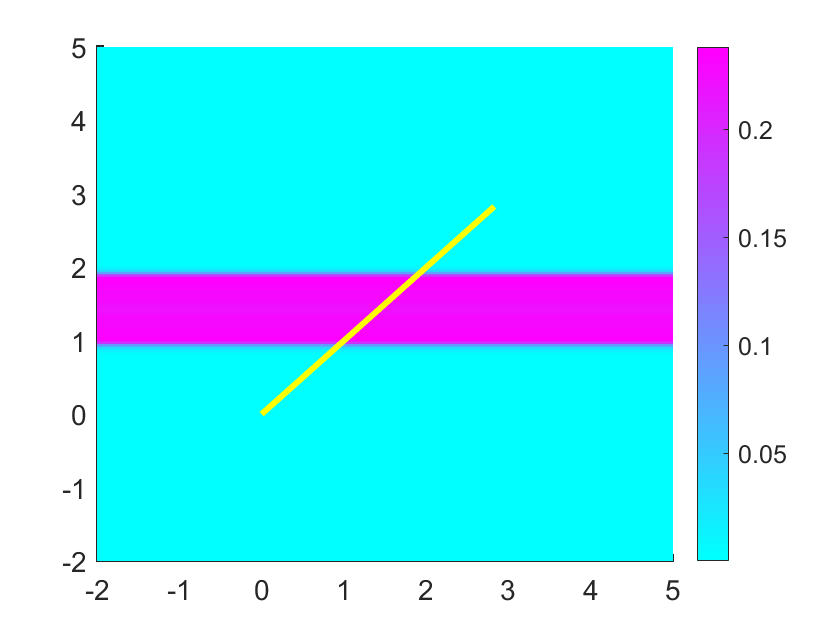}

}

\caption{Reconstruction from a single observable direction $\hat x=(\cos\theta, \sin\theta)$  with $\theta \in [11\pi/12, 19\pi/12]$ for a straight line segment $a(t)=(2\sqrt 2t, 2\sqrt 2t)$ with $t\in[1,2]$. Since $h'(t)<0$, the strip $K_\Gamma^{(\hat{x})}$ is a subset of $\{y\in \R^2: \hat{x}\cdot y\in \hat{x}\cdot \Gamma\}$.} \label{fig:xie-line22}
\end{figure}

In Fig.\ref{fig:xie-line22}, we measure the data at the observable angle $\theta \in [11\pi/12, 19\pi/12]$ so that $ h'(t)<0$. Although these observation directions $\theta$ belong to the class of the observable set, the recovered strips $K_{\Gamma}^{(\hat x)}$ are thinner than the smallest strips  containing the trajectory of the moving source, because  $K_{\Gamma}^{(\hat x)}\subset \{y\in\R^2:\inf (\hat x \cdot \Gamma) \leq \hat x\cdot y\leq \sup( \hat x \cdot \Gamma)\}$ by Lemma \ref{lem3.9}.

\begin{figure}[htb]
\centering
\subfigure[$\theta=19\pi/24$]{

\includegraphics[scale=0.22]{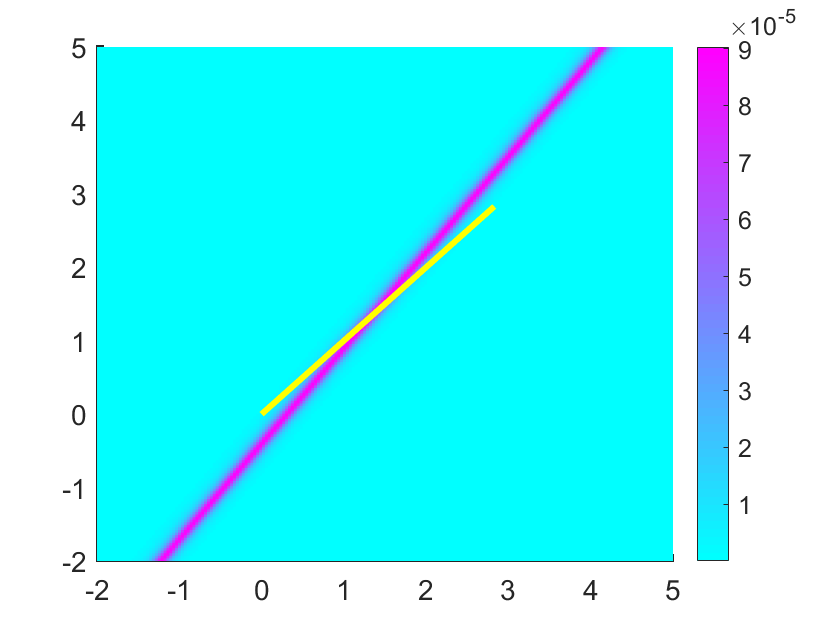}

}
\subfigure[$\theta=20\pi/24$ ]{
\includegraphics[scale=0.22]{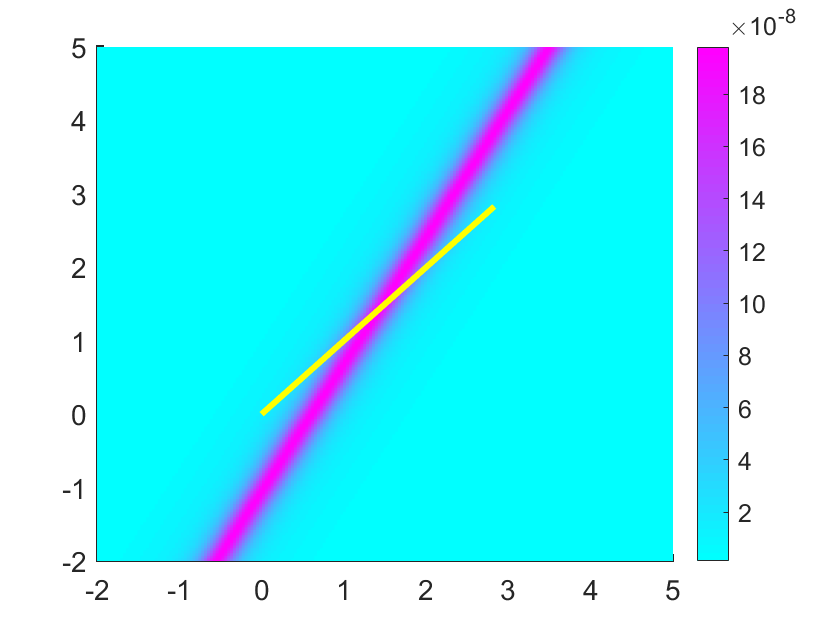}

}
\subfigure[$\theta=21\pi/24$]{
\includegraphics[scale=0.22]{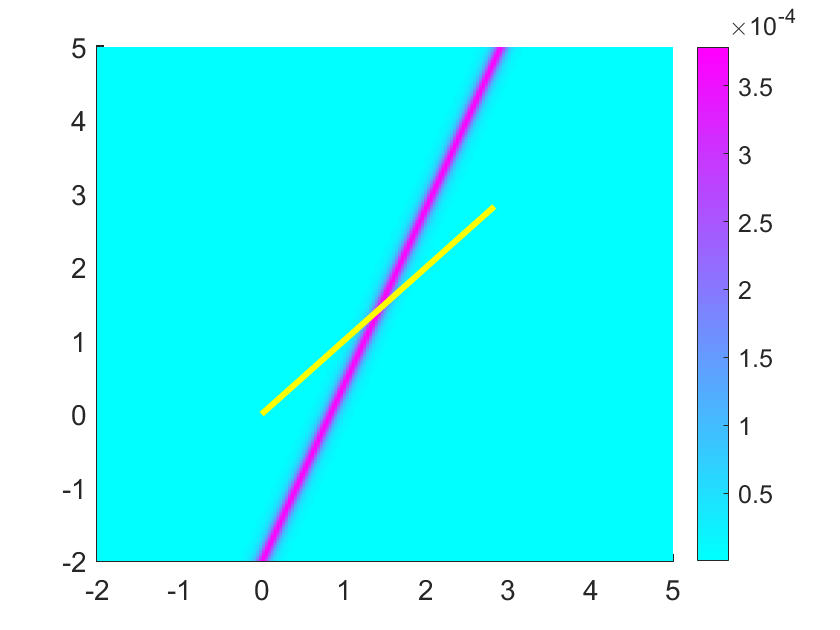}

}
\subfigure[$\theta=39\pi/24$]{

\includegraphics[scale=0.22]{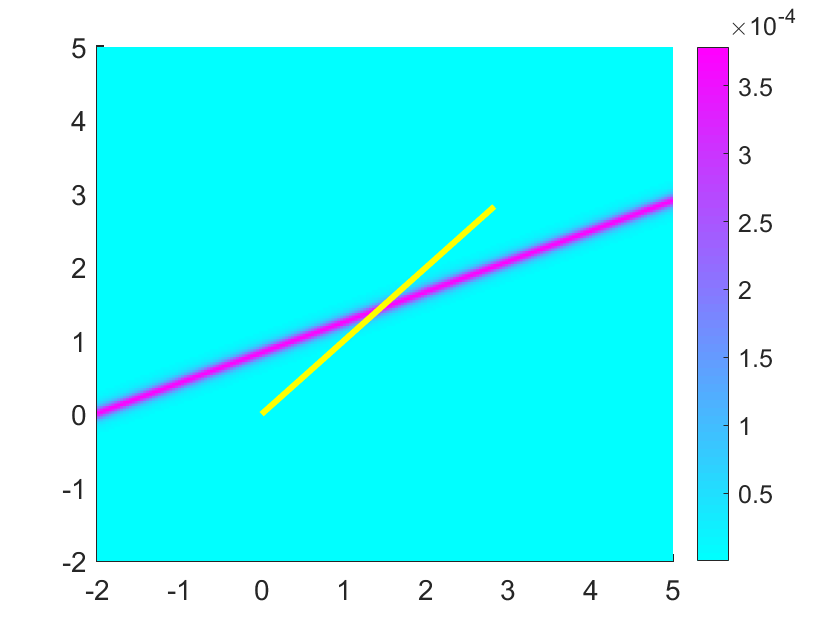}

}
\subfigure[$\theta=40\pi/24$ ]{
\includegraphics[scale=0.22]{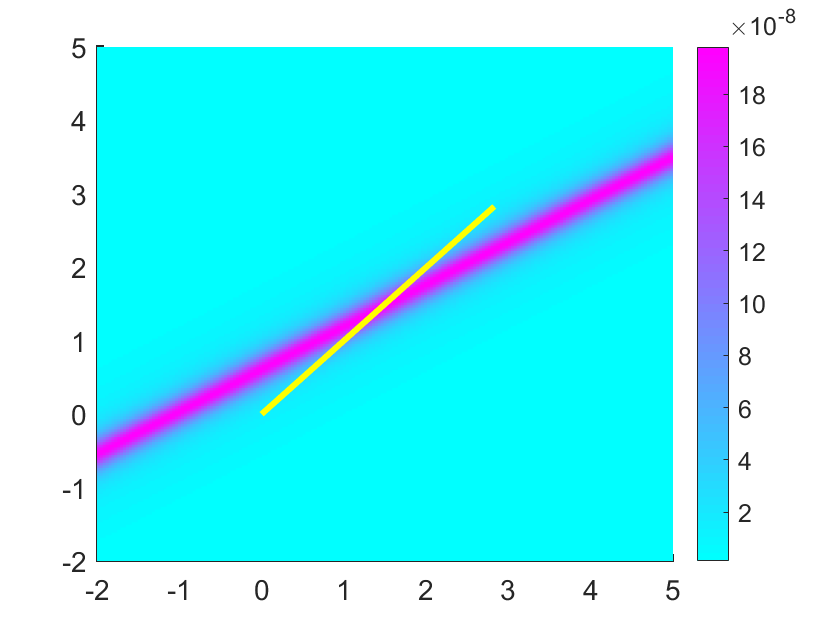}

}
\subfigure[$\theta=41\pi/24$]{
\includegraphics[scale=0.22]{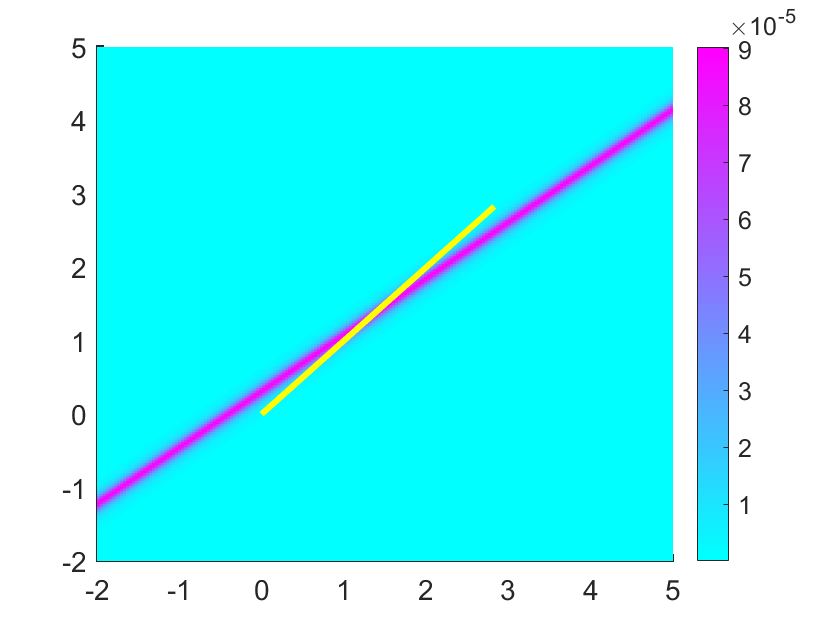}

}
\caption{Reconstruction from a single non-observable direction $\hat x=(\cos\theta, \sin\theta)$  with $\theta \in (3\pi/4,11\pi/12)\cup(19\pi/12, 7\pi/4)$ for a straight line segment $a(t)=(2\sqrt 2t, 2\sqrt 2t)$ with $t\in[1,2]$. } \label{fig:xie-line3}
\end{figure}

In Fig.\ref{fig:xie-line3}, we make use of non-observable angles. The numerical results illustrate that the indicator values are indeed much smaller. Hence, one cannot expect to reconstruct the smallest strip containing the trajectory of the moving source.



\vspace{.1in}\textbf{Example 2: An arc in $\R^2$}\vspace{.1in}

\begin{figure}[htb]
\centering
\subfigure[$\theta=5\pi/8$]{

\includegraphics[scale=0.22]{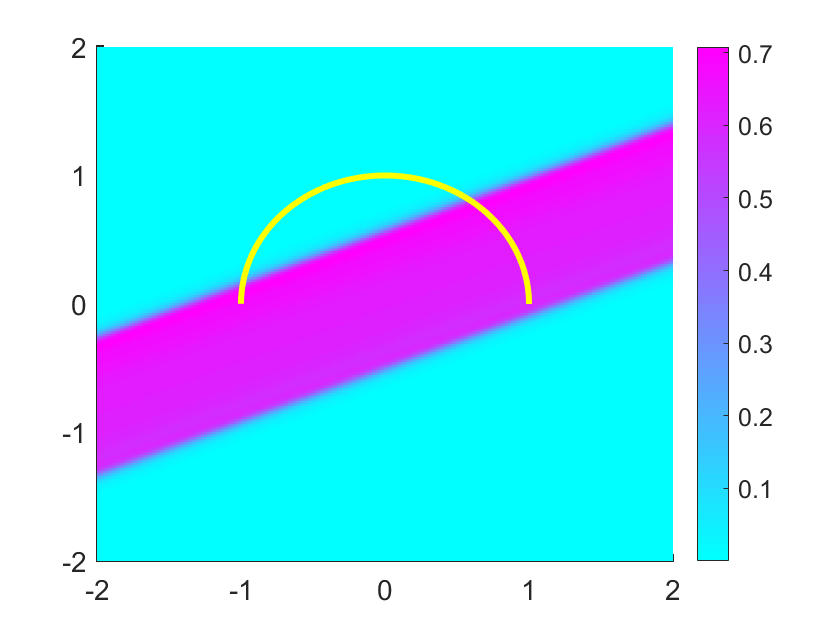}

}
\subfigure[$\theta=6\pi/8$ ]{
\includegraphics[scale=0.22]{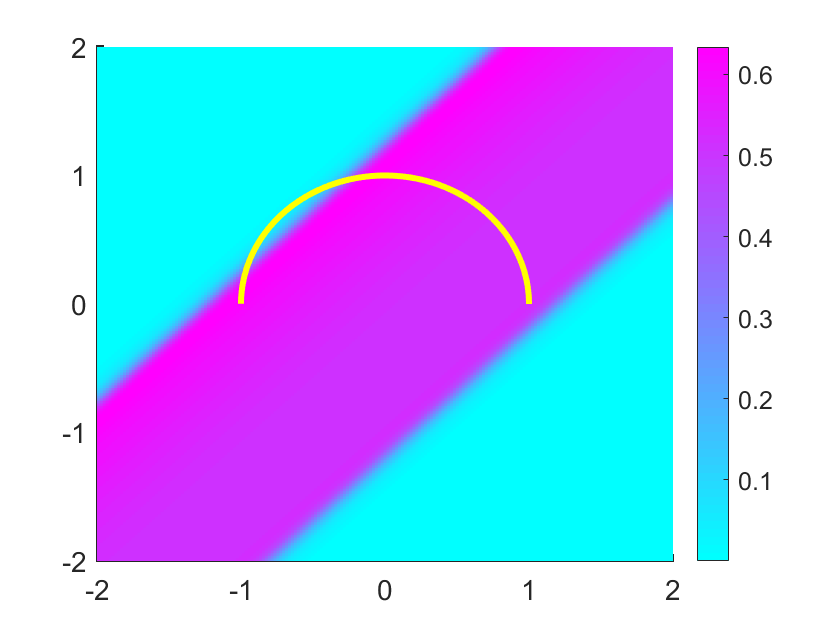}

}
\subfigure[$\theta=7\pi/8$]{
\includegraphics[scale=0.22]{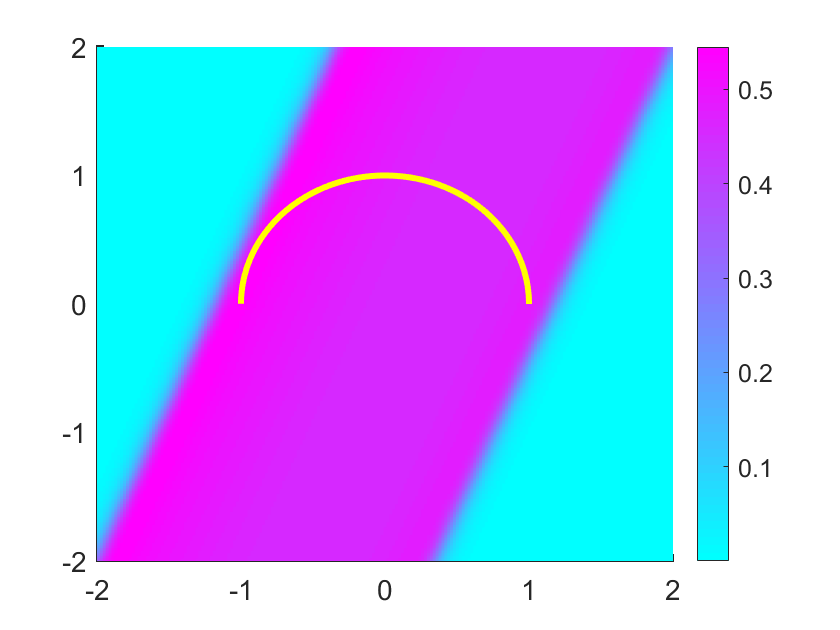}

}
\subfigure[$\theta=8\pi/8$]{

\includegraphics[scale=0.22]{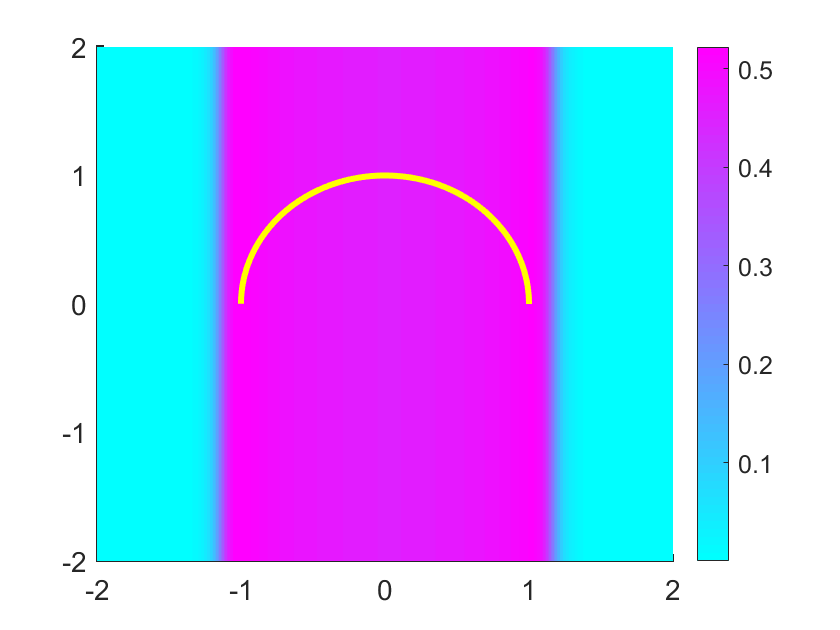}

}
\subfigure[$\theta=9\pi/8$ ]{
\includegraphics[scale=0.22]{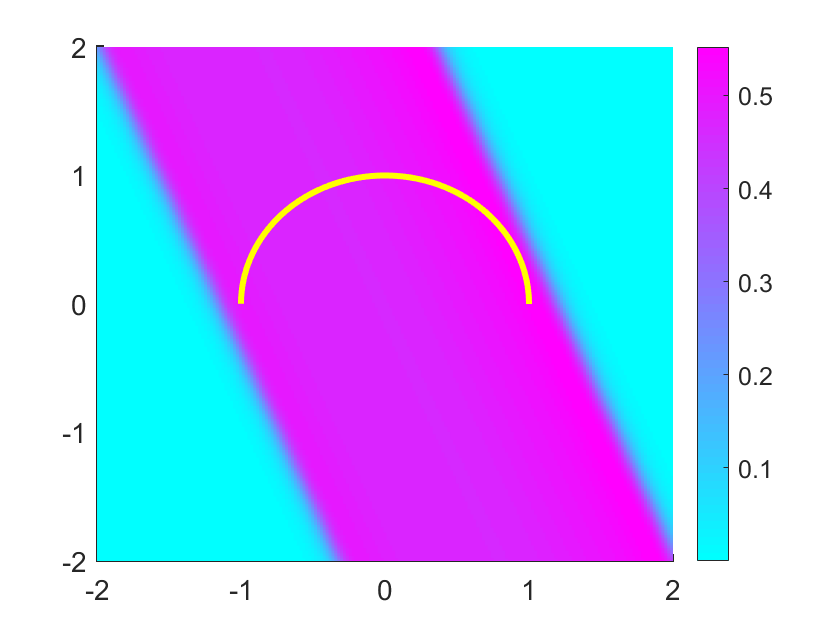}

}
\subfigure[$\theta=10\pi/8$]{
\includegraphics[scale=0.22]{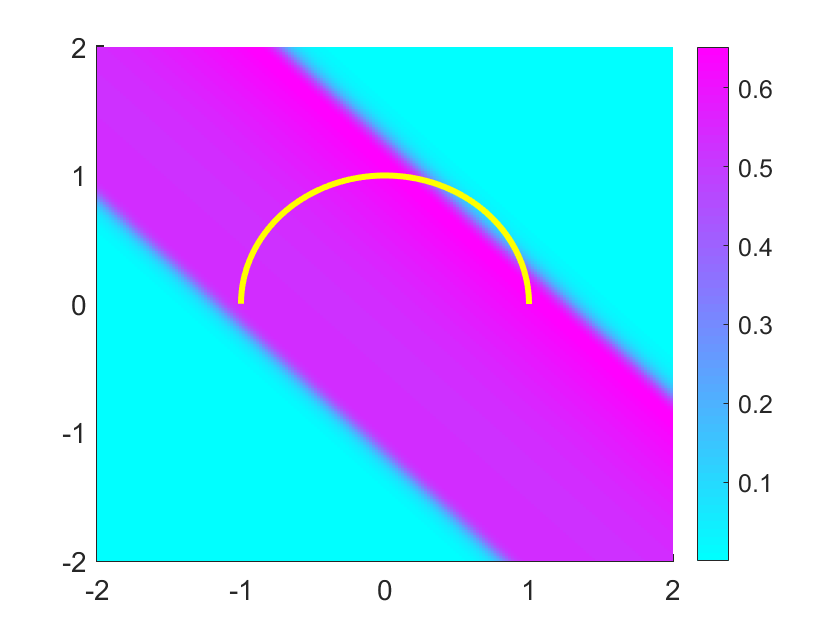}

}
\caption{Reconstruction from a single observable direction $\hat{x}=(\cos\theta, \sin\theta)$ with $\theta\in[\pi/2, 3\pi/2]$  for an arc  $a(t)=(\cos t, \sin t)$ with $t\in[0,\pi]$.} \label{fig:arc1}
\end{figure}

\begin{figure}[htb]
\centering
\subfigure[$\theta=1\pi/8$]{

\includegraphics[scale=0.22]{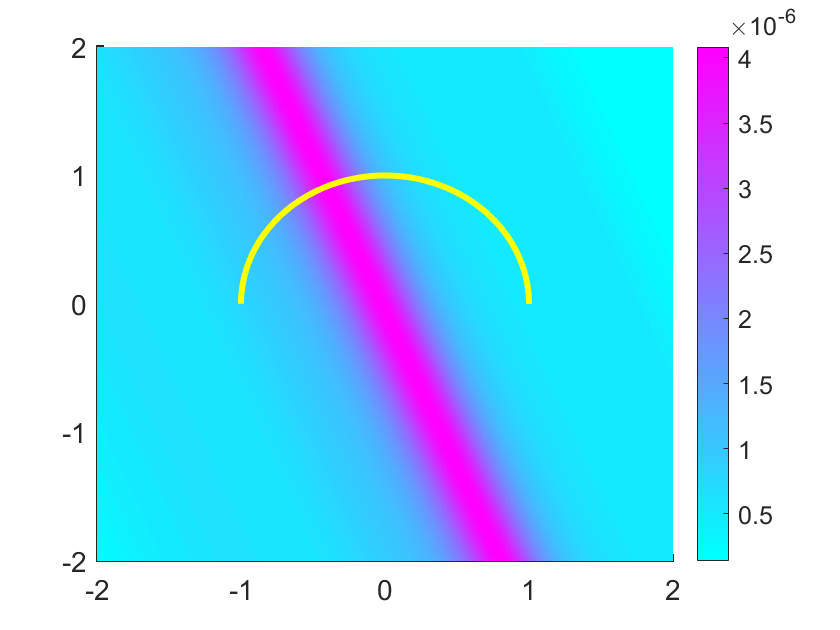}

}
\subfigure[$\theta=2\pi/8$ ]{
\includegraphics[scale=0.22]{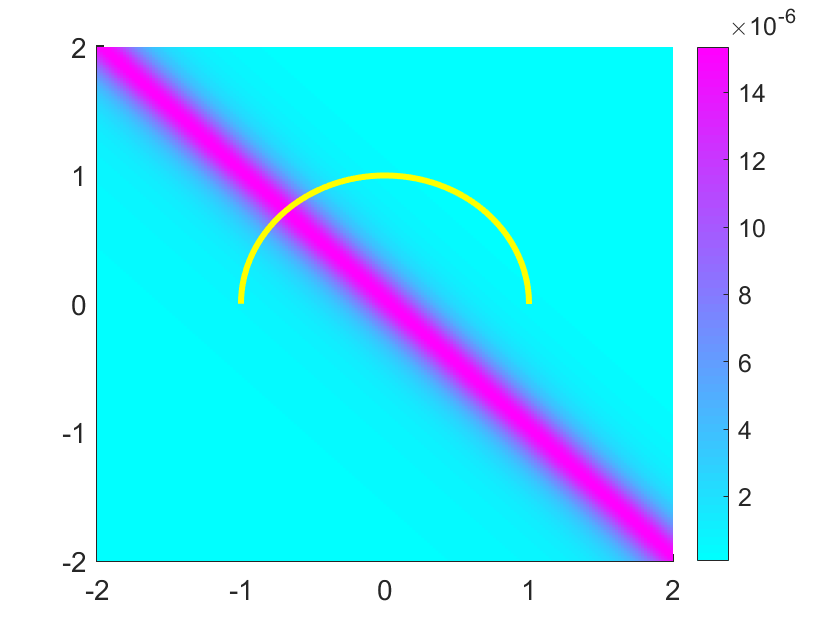}

}
\subfigure[$\theta=3\pi/8$]{
\includegraphics[scale=0.22]{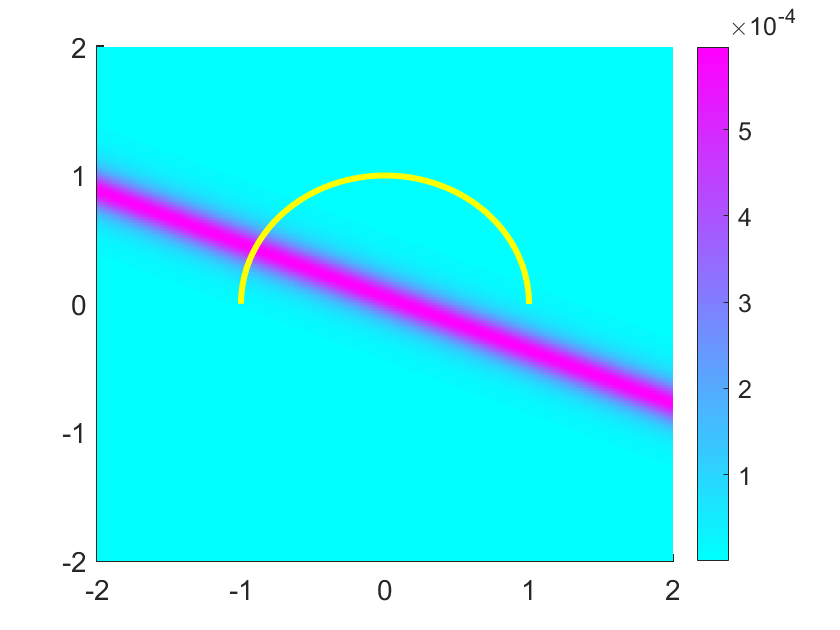}

}
\subfigure[$\theta=13\pi/8$]{

\includegraphics[scale=0.22]{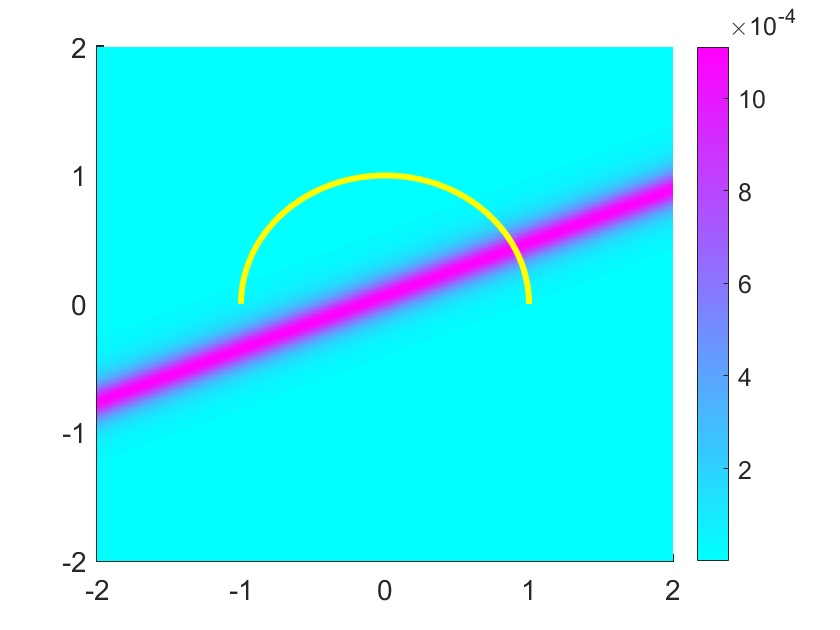}

}
\subfigure[$\theta=15\pi/8$ ]{
\includegraphics[scale=0.22]{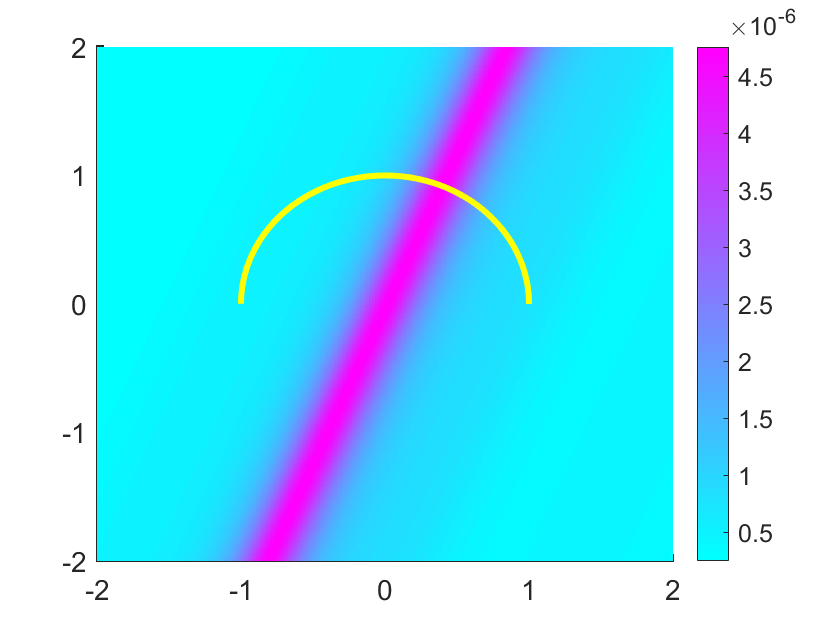}

}
\subfigure[$\theta=16\pi/8$]{
\includegraphics[scale=0.22]{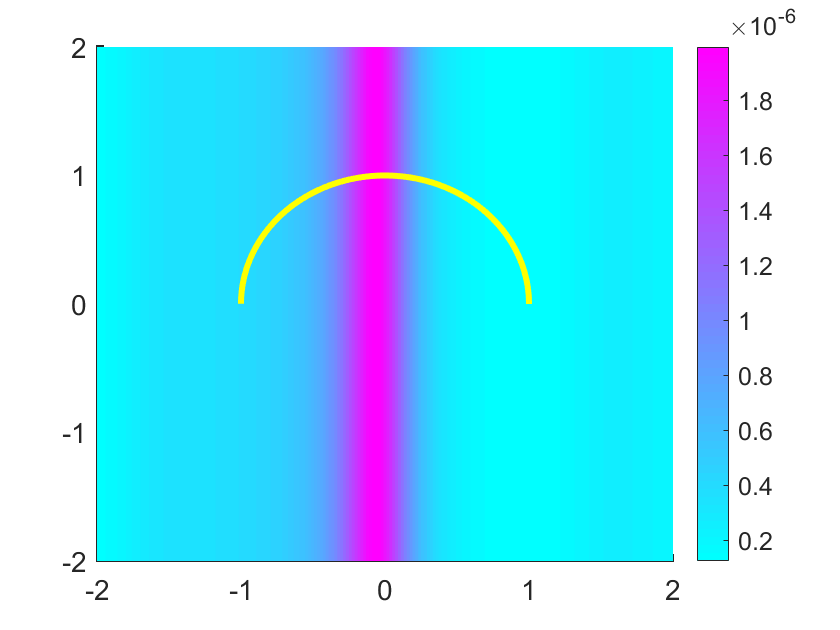}

}
\caption{Reconstruction from a single non-observable direction $\hat{x}=(\cos\theta, \sin\theta)$ with $\theta \in (0, \pi/2) \cup  (3\pi/2, 2\pi)$  for an arc  $a(t)=(\cos t, \sin t)$ with $t\in[0,\pi]$.} \label{fig:arc2}
\end{figure}

As shown in Example 2 of Section \ref{RangeLx}, we take $a(t)=(\cos(t), \sin(t))$ with $t\in [0,\pi]$. The search domain is $[-2,2]^2$. From Lemma \ref{arc-ob}, we know that observable directions are  $\hat{x}$ with $\theta \in [\pi/2, 3\pi/2]$ and non-observable directions are  $\hat{x}$ with $\theta \in (0, \pi/2) \cup  (3\pi/2, 2\pi)$.
 Fig.\ref{fig:arc1} shows the reconstructions using the data from
 observable directions, where the subfigures (c), (d) and (e) nicely give us the the smallest strip $K_{\Gamma}^{(\hat{x})}$ containing the trajectory of the moving source that is perpendicular to the observable direction. Note that $K_{\Gamma}^{(\hat{x})}= \{y\in \R^2: \sup(\hat x \cdot \Gamma)\leq \hat x \cdot y\leq \inf (\hat x \cdot \Gamma) \}$ for $\theta=7\pi/8$,  $8\pi/8$,  and $9\pi/8$, because  $h'(t)>0$ for all $t\in[0,\pi]$ at these angles. However, the strips $K_{\Gamma}^{(\hat{x})}$ in subfigures (a), (b) and (f) do not provide sufficient information on the trajectory. This is due to the reason that $h'(t)<0$ for $\theta=5\pi/8$, $6\pi/8$, $10\pi/8$ and $t\in[0,\pi]$, implying that
 $K_{\Gamma}^{(\hat{x})}\subset \{y\in \R^2: \sup(\hat x \cdot \Gamma)\leq \hat x \cdot y\leq \inf (\hat x \cdot \Gamma) \}$.
  Reconstructions from non-observable angles  $\theta \in(0, \pi/2) \cup  (3\pi/2, 2\pi)$ are illustrated in Fig.\ref{fig:arc2}. The values are still very small and can not reconstruct the strip $\{y\in \R^2: \sup(\hat x \cdot \Gamma)\leq \hat x \cdot y\leq \inf (\hat x \cdot \Gamma) \}$.

	\vspace{.1in}\textbf{Example 3: A piecewise linear curve in $\R^2$} \vspace{.1in}	



	We first remark that the analysis performed in Sections 2-4 carry over to piecewisely $C^1$-smooth orbit functions. Complexity arises only from the definition of the division points made in Def. \ref{DIDP}, where the discontinuity points of $h'(t)$ should be taken into account.	
Assume that the trajectory of the moving source $a(t)$ is given by
		\begin{equation*}
			a(t) = \left\{
			\begin{aligned}
				&(-t+3,-t+3),\quad &t\in[0,1],\\
				&(t+1,-t+3), &t\in[1,2].
			\end{aligned}
			\right.
		\end{equation*}
Let $\hat{x} = (\cos \theta, \sin \theta),\,\theta\in[0,2\pi)$ be the observation direction. We first calculate the observable and non-observable directions. Note that $T=2$. Evidently,
		\begin{equation*}
			h(t) = t+\hat{x}\cdot a(t) = \left\{
			\begin{aligned}
				& t (1- \cos \theta - \sin \theta) +3(\cos \theta + \sin \theta), \, &t\in [0,1],\\
				& t (1+ \cos \theta - \sin \theta) +\cos \theta + 3\sin \theta, \, &t\in [1,2],
			\end{aligned}
			\right.
		\end{equation*}
		and thus
		\begin{equation*}
			h^{\prime}(t) = 1+\hat{x}\cdot a^{\prime}(t) = \left\{
			\begin{aligned}
				& 1- \cos \theta - \sin \theta, \, &t\in (0,1),\\
				& 1+ \cos \theta - \sin \theta, \, &t\in (1,2).
			\end{aligned}
			\right.
		\end{equation*}
Since $h^{\prime}(t) \equiv 0$ in some interval when $\theta = 0,\pi/2$ and $\pi$,  we need to consider the following six cases separately.
		
		\textbf{(1) $\theta =0$.}
We have $h^{\prime}(t)=0$ for $t\in (0,1)$ and $h^{\prime}(t)>0$ for $t\in (1,2)$. This gives $3\leq h(t)\leq 5$ for $t\in(0,2)$, implying $\xi^{(\hat{x})}_{\max} - \xi^{(\hat{x})}_{\min} = T$. Thus, $\hat{x}=(1,0)$ is an observable direction.

        \textbf{(2) $\theta \in (0,\pi/2)$}.	We have $h^{\prime}(t)<0$ for $t\in (0,1)$ and $h^{\prime}(t)>0$ for $t\in (1,2)$. Therefore, for $\xi=h(t), t\in(0, 2)$ it holds that
		\begin{equation*}
			\begin{aligned}
				\xi &\in [1+2(\cos \theta+\sin \theta),3(\cos \theta+\sin \theta)] \cup [1+2(\cos \theta+\sin \theta),2+3\cos \theta -\sin \theta] \\
				&= [1+2(\cos \theta+\sin \theta),\max\{3(\cos \theta+\sin \theta),2+3\cos \theta -\sin \theta\}].
			\end{aligned}
		\end{equation*}
Consequently, $\xi^{(\hat{x})}_{\max} - \xi^{(\hat{x})}_{\min} = \max\{\cos \theta +\sin \theta -1,1+\cos \theta -3\sin \theta\} < T=2$. Thus, each $\hat{x}$ with $\theta \in (0,\pi/2)$.is non-observable.

		\textbf{(3) $\theta =\pi/2$}. We have $h^{\prime}(t)=0$ for $t\in (0,1)$ and $h^{\prime}(t)=0$ for $t\in (1,2)$, implying that $h(t)\equiv 3$.
		Hence,  $\xi^{(\hat{x})}_{\max} - \xi^{(\hat{x})}_{\min} = 0 < T$. Thus, $\hat{x}=(0,1)$ is an non-observable.

		\textbf{(4) $\theta \in (\pi/2,\pi).$} We have $h^{\prime}(t)>0$ for $t\in (0,1)$ and $h^{\prime}(t)<0$ for $t\in (1,2)$. Hence, if $\xi=h(t)$ for some $t\in(0,2)$, then
		\begin{equation*}
			\begin{aligned}
				\xi &\in [3(\cos \theta+\sin \theta),1+2(\cos \theta+\sin \theta)] \cup [2+3\cos \theta -\sin \theta,1+2(\cos \theta+\sin \theta)] \\
				&= [\min\{3(\cos \theta+\sin \theta),2+3\cos \theta -\sin \theta\},1+2(\cos \theta+\sin \theta)].
			\end{aligned}
		\end{equation*}
In this case, we get $\xi^{(\hat{x})}_{\max} - \xi^{(\hat{x})}_{\min} = \max\{1-\cos \theta -\sin \theta,3\sin \theta - \cos \theta -1\} < T$. Thus, the direction $\hat{x}$ with $\theta \in (\pi/2,\pi)$
is non-observable.

		\textbf{(5) $\theta = \pi.$} We have $h^{\prime}(t)>0$ in $(0,1)$ and $h^{\prime}(t)=0$ in $(1,2)$, implying $-3\leq h(t)\leq -1$ for $t\in(0,2)$.
Thus $\xi^{(\hat{x})}_{\max} - \xi^{(\hat{x})}_{\min} = 2$ and $\hat{x}=(-1, 0)$ is an observable direction.

		\textbf{(6) $\theta \in (\pi,2\pi)$}. We have $h^{\prime}(t)>0$ for $t\in (0,1)$ and $h^{\prime}(t)>0$ for $t\in (1,2)$. For $\xi=h(t)$ we have
		\begin{equation*}
			\begin{aligned}
				\xi &\in [3(\cos \theta+\sin \theta),1+2(\cos \theta+\sin \theta)] \cup [1+2(\cos \theta+\sin \theta),2+3\cos \theta -\sin \theta] \\
				&= [3(\cos \theta+\sin \theta),2+3\cos \theta -\sin \theta],
			\end{aligned}
		\end{equation*}
implying that $\xi^{(\hat{x})}_{\max} - \xi^{(\hat{x})}_{\min} = 2-4\sin \theta > 2$. Therefore, each direction $\hat{x}$ with $\theta\in(0,\pi)$ is observable.

\begin{figure}[htb]
\centering
\subfigure[$\theta=11\pi/10$]{
\includegraphics[scale=0.22]{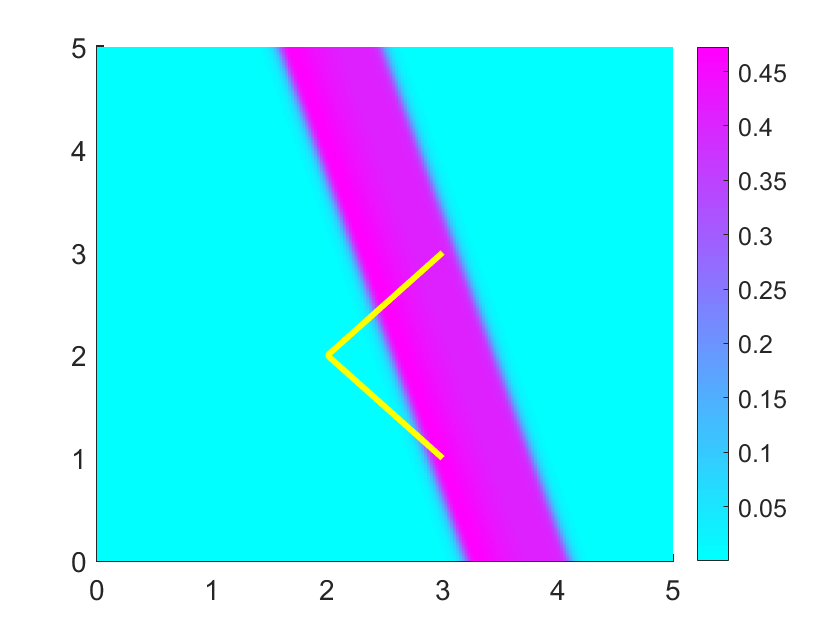}

}
\subfigure[$\theta=5\pi/4$ ]{
\includegraphics[scale=0.22]{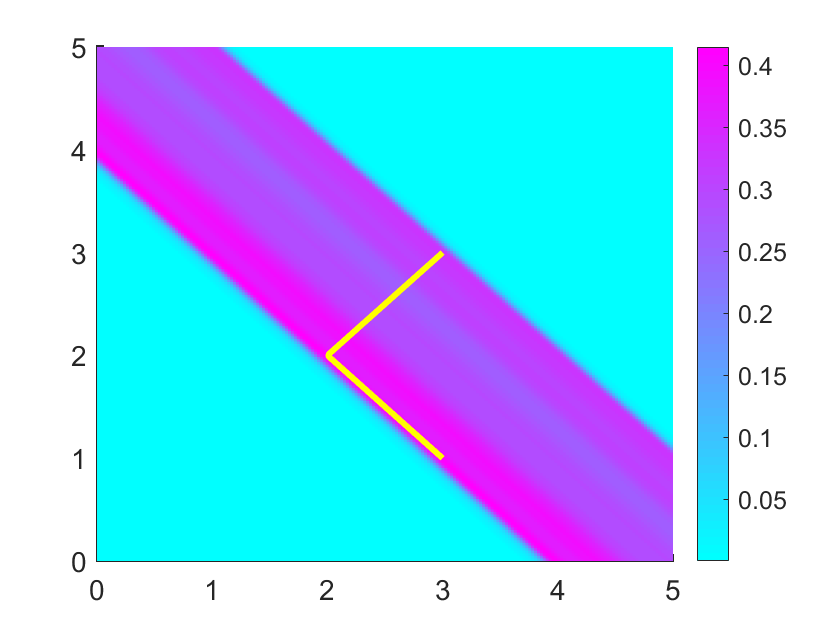}

}
\subfigure[$\theta=3\pi/2$]{
\includegraphics[scale=0.22]{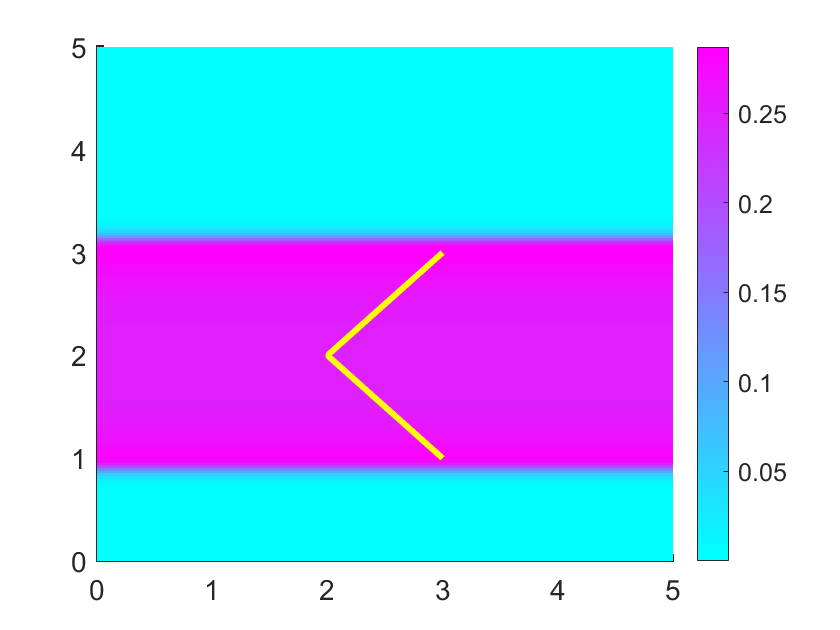}

}
\subfigure[$\theta=5\pi/3$]{

\includegraphics[scale=0.22]{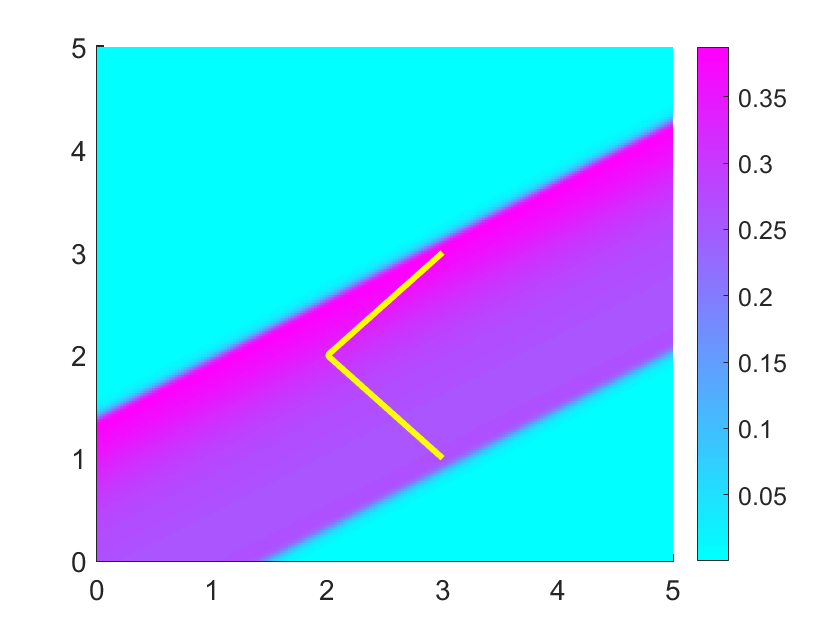}

}
\subfigure[$\theta=7\pi/4$ ]{
\includegraphics[scale=0.22]{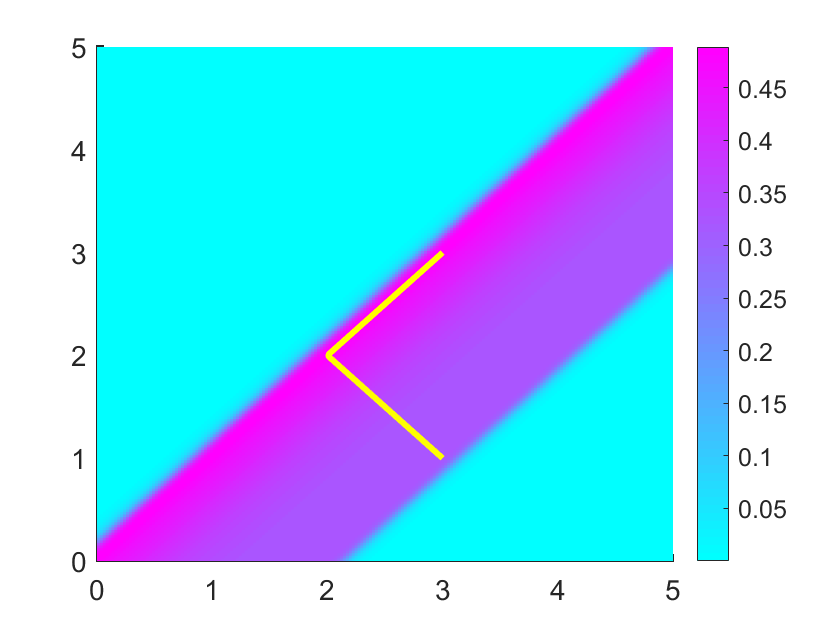}

}
\subfigure[$\theta=15\pi/8$]{
\includegraphics[scale=0.22]{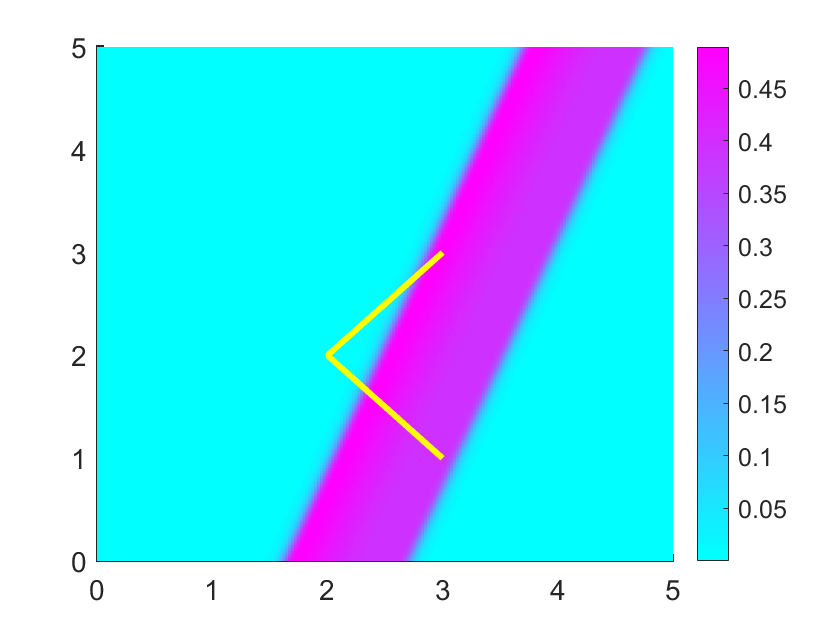}

}
%
%
%
%
\caption{Reconstruction from a single observable direction $\hat{x}=(\cos\theta, \sin\theta)$ with $\theta\in[\pi,2\pi]$  for a broken line segment $a(t)=(-t+3,-t+3)$ with $t\in[0,1]$ and $a(t)=(t+1,-t+3)$ with $t\in[1,2]$ in $\R^2$. } \label{fig:broken1}
\end{figure}

\begin{figure}[htb]
\centering
\subfigure[$\theta=\pi/7$]{

\includegraphics[scale=0.22]{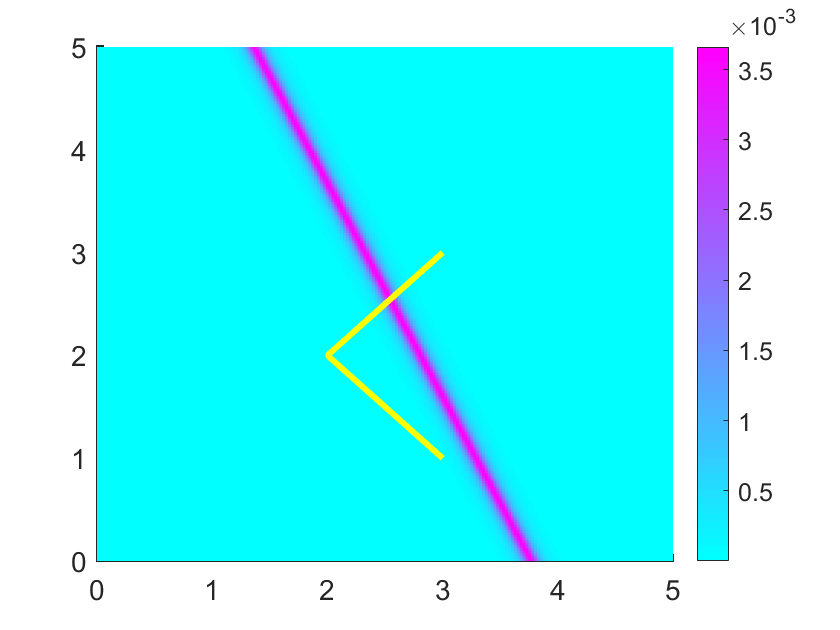}

}
\subfigure[$\theta=2\pi/7$ ]{
\includegraphics[scale=0.22]{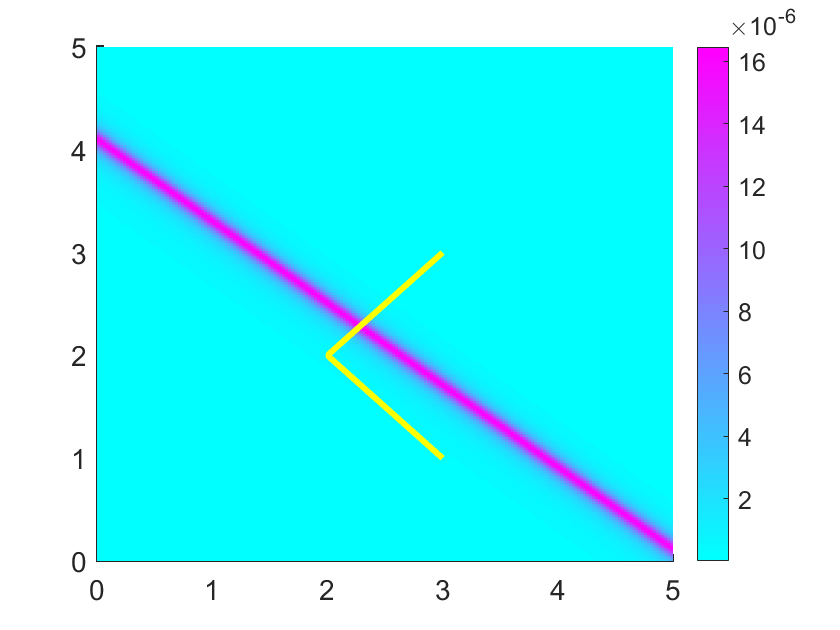}

}
\subfigure[$\theta=3\pi/7$]{
\includegraphics[scale=0.22]{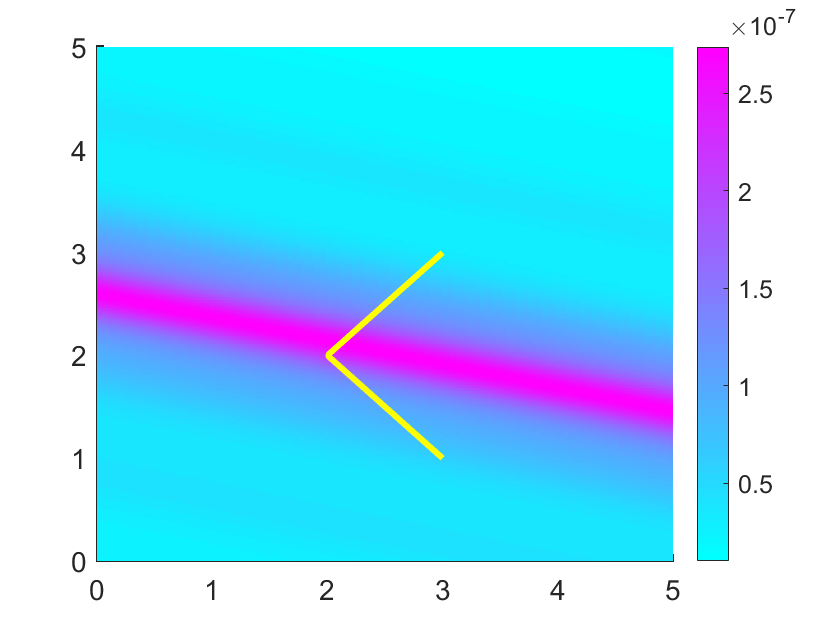}

}
\subfigure[$\theta=4\pi/7$]{

\includegraphics[scale=0.22]{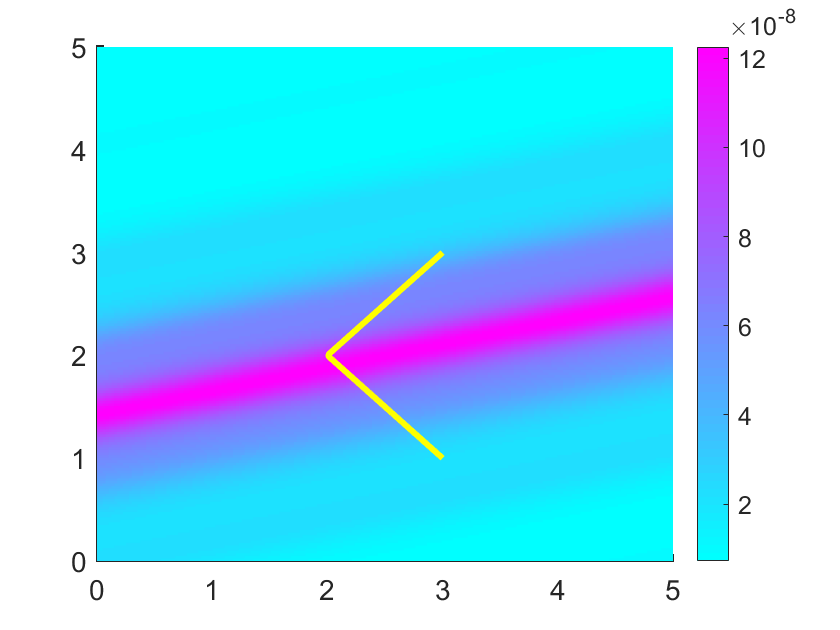}

}
\subfigure[$\theta=5\pi/7$ ]{
\includegraphics[scale=0.22]{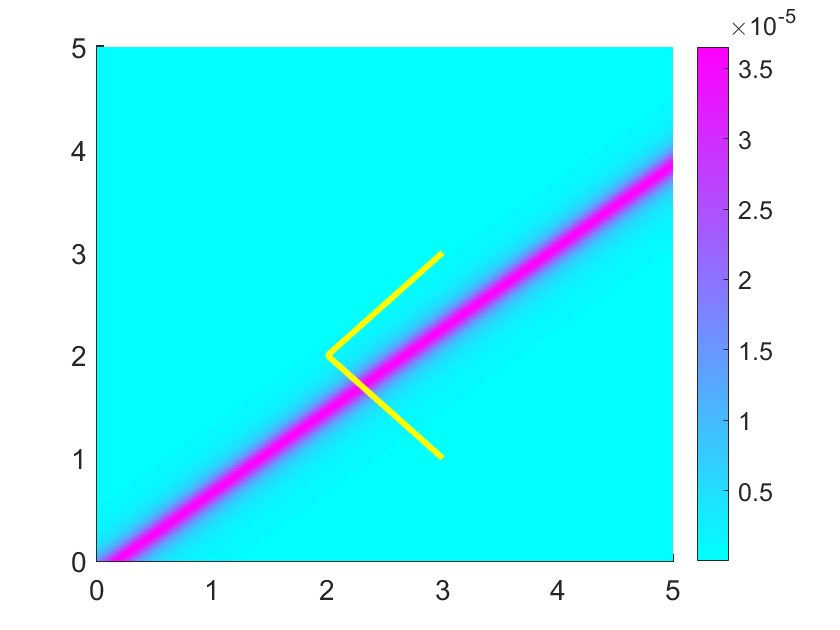}

}
\subfigure[$\theta=6\pi/7$]{
\includegraphics[scale=0.22]{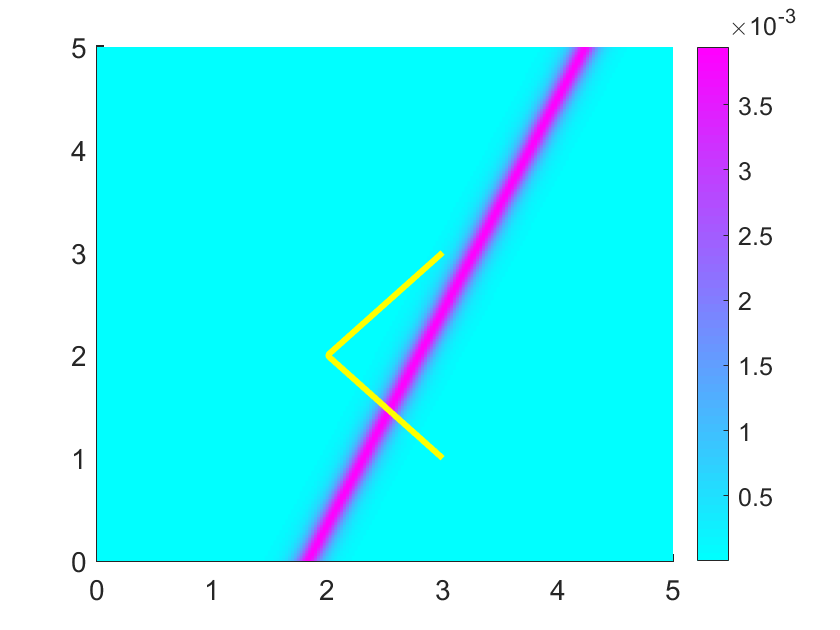}

}
\caption{Reconstruction from a single non-observable direction $\hat{x}=(\cos\theta, \sin\theta)$ with $\theta\in(0,\pi)$  for a piecewise linear curve $a(t)=(-t+3,-t+3)$ with $t\in[0,1]$ and $a(t)=(t+1,-t+3)$ with $t\in[1,2]$ in $\R^2$. } \label{fig:broken2}
\end{figure}

Summing up,  we conclude that $[\pi,2\pi]$ consists of observable angles and $(0, \pi)$ the non-observable ones.  In Fig.\ref{fig:broken1}, we plot the indicator functions for different observable angles  in $[\pi,2\pi]$. 
In subfigures (b), (c), (d) and (e),
the reconstructed strip $K_{\Gamma}^{(\hat{x})}$ coincides with $\{y\in \R^2: \sup(\hat x \cdot \Gamma)\leq \hat x \cdot y\leq \inf (\hat x \cdot \Gamma) \}$  for all $t\in[0,2]$ and
$\theta=5\pi/4$, $3\pi/2$, $5\pi/3$  and  $7\pi/4$.
The strips $K_{\Gamma}^{(\hat{x})}$ in subfigures (a) and (f) are subsets of $\{y\in \R^2: \sup(\hat x \cdot \Gamma)\leq \hat x \cdot y\leq \inf (\hat x \cdot \Gamma) \}$ for all $t\in[0,2]$ and
 $\theta =11\pi/10$ and $15\pi/8$.
 In Fig.\ref{fig:broken2} we show reconstructions from non-observable angles $\theta\in(0,\pi)$.

\vspace{.1in}\textbf{Example 4: A straight line segment in $\R^3$}\vspace{.1in}

\begin{figure}[htb]
\centering
\subfigure[$\phi=\pi/6, \theta=\pi/8$]{
\includegraphics[scale=0.3]{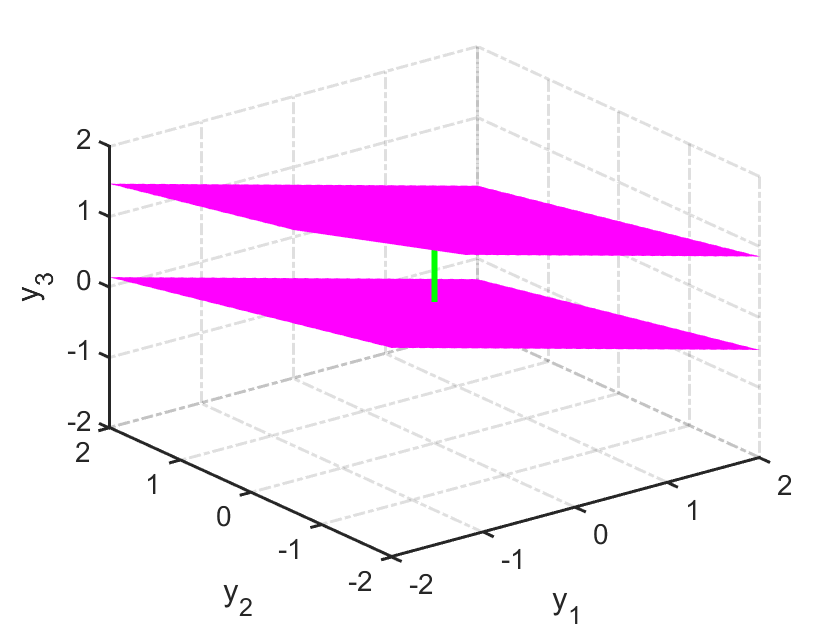}}
\subfigure[$\phi=\pi/4, \theta=\pi/6$ ]{
\includegraphics[scale=0.3]{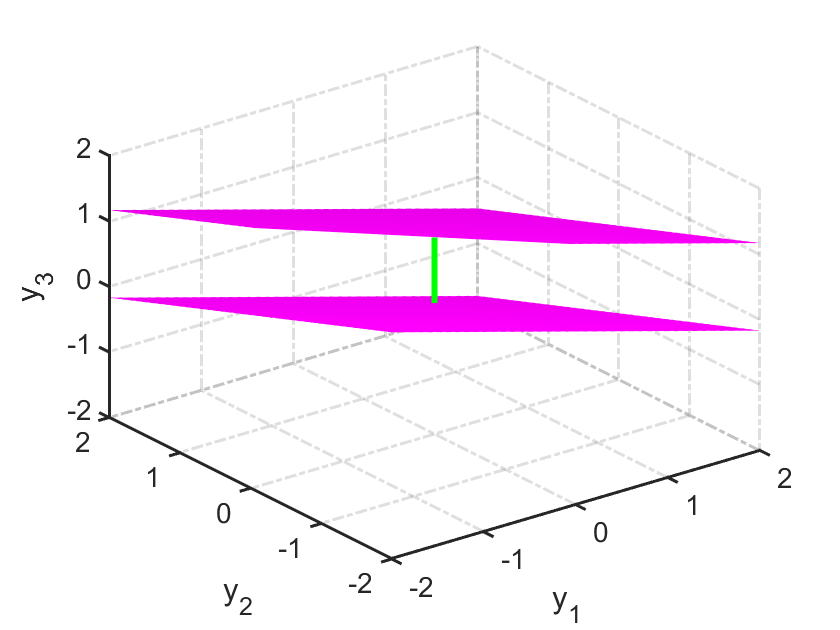}}
\caption{Illustration of the hyperplanes perpendicular to the observable direction $\hat x= (\sin\theta \cos \phi,  \sin\theta\sin\phi,  \cos \phi)$. Here we take isosurface level $=0.02$.} \label{fig:hyperplane}
\end{figure}

\begin{figure}[htb]
\centering
\subfigure[$\phi=\pi/4, \theta=\pi/8$]{

\includegraphics[scale=0.22]{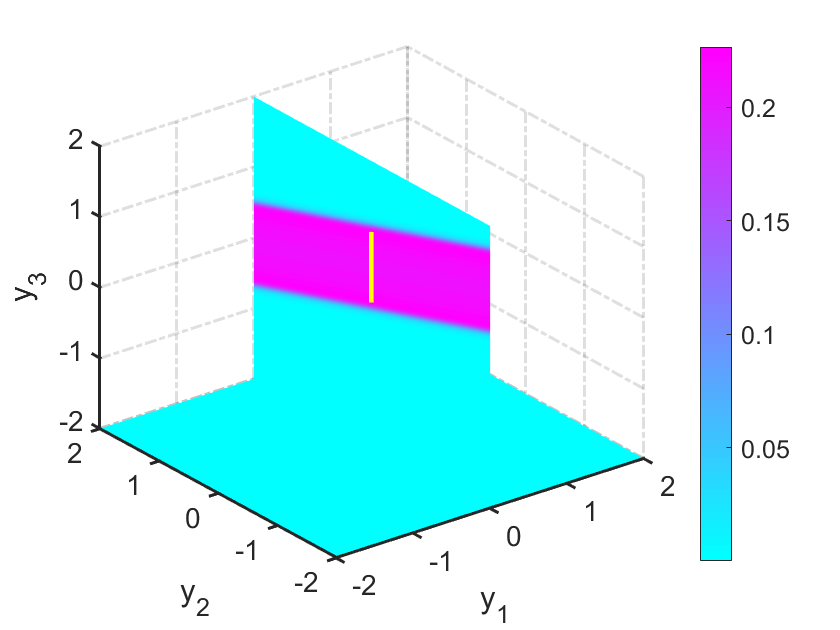}

}
\subfigure[$\phi=\pi/4, \theta=2\pi/8$ ]{
\includegraphics[scale=0.22]{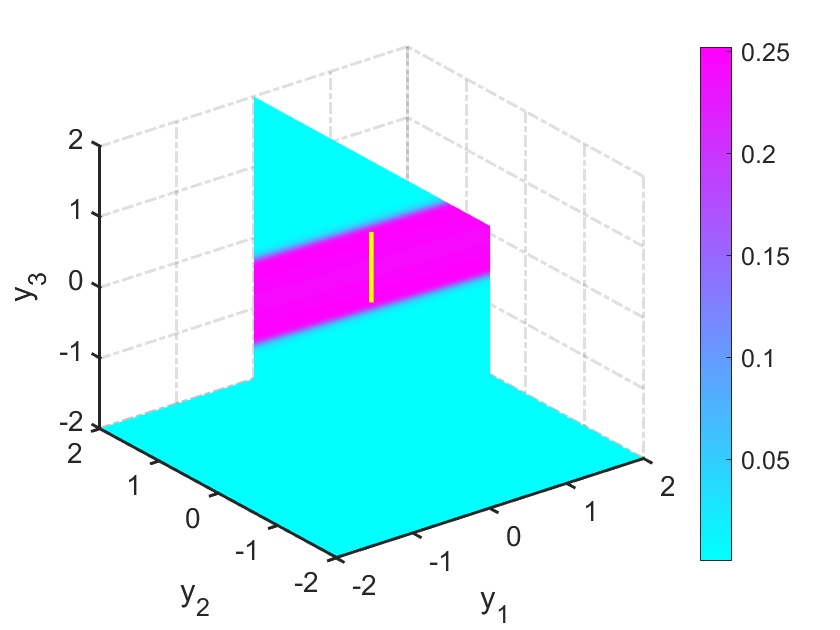}

}
\subfigure[$\phi=\pi/4, \theta=3\pi/8$]{
\includegraphics[scale=0.22]{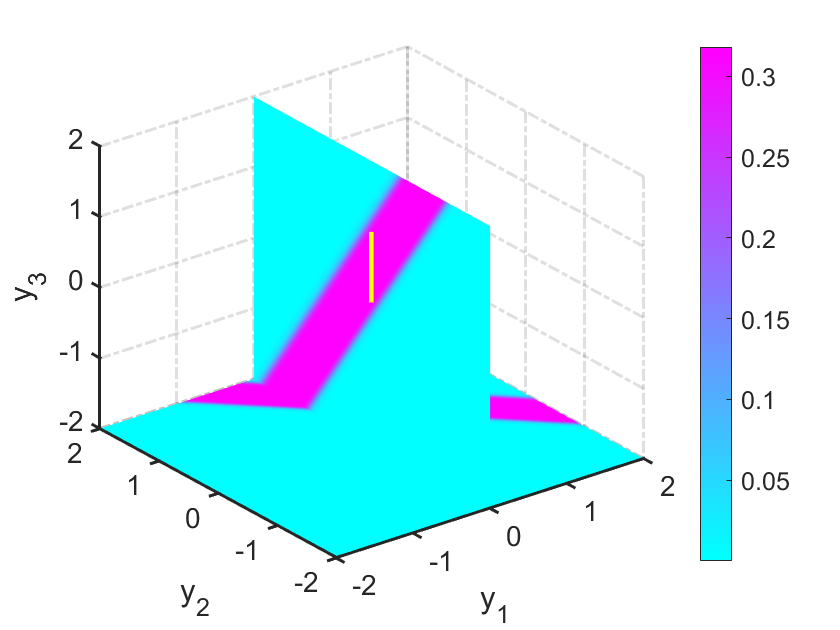}

}

\subfigure[$\phi=\pi/2, \theta=\pi/8$]{

\includegraphics[scale=0.22]{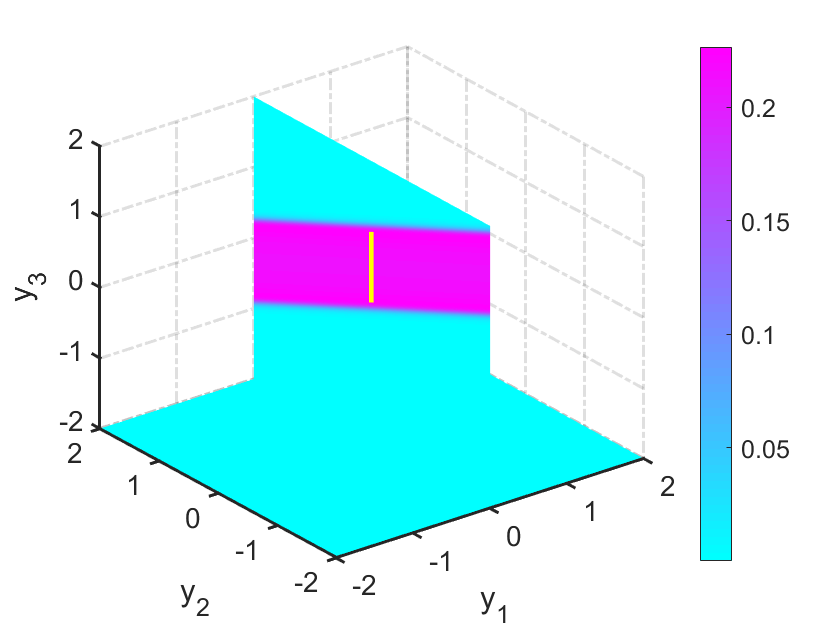}

}
\subfigure[$\phi=\pi/2, \theta=2\pi/8$ ]{
\includegraphics[scale=0.22]{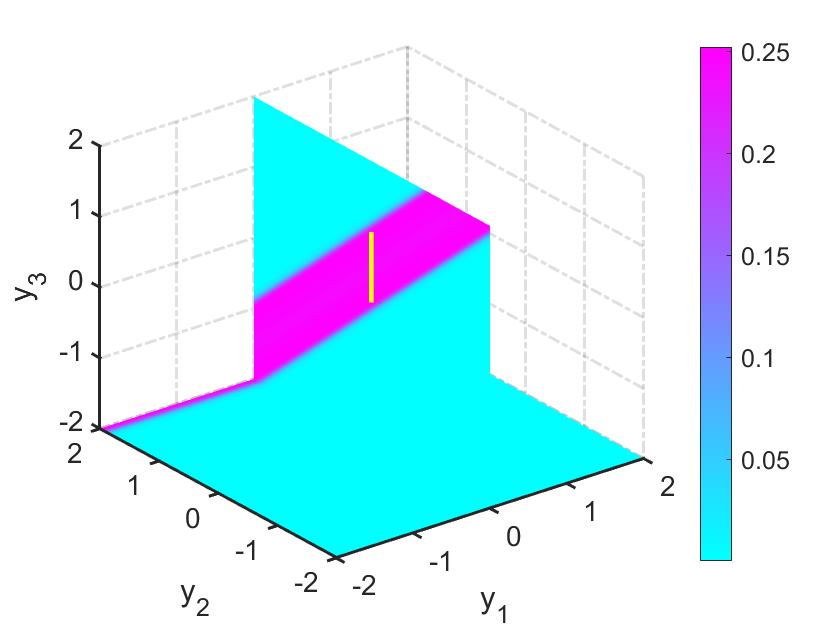}

}
\subfigure[$\phi=\pi/2, \theta=3\pi/8$]{
\includegraphics[scale=0.22]{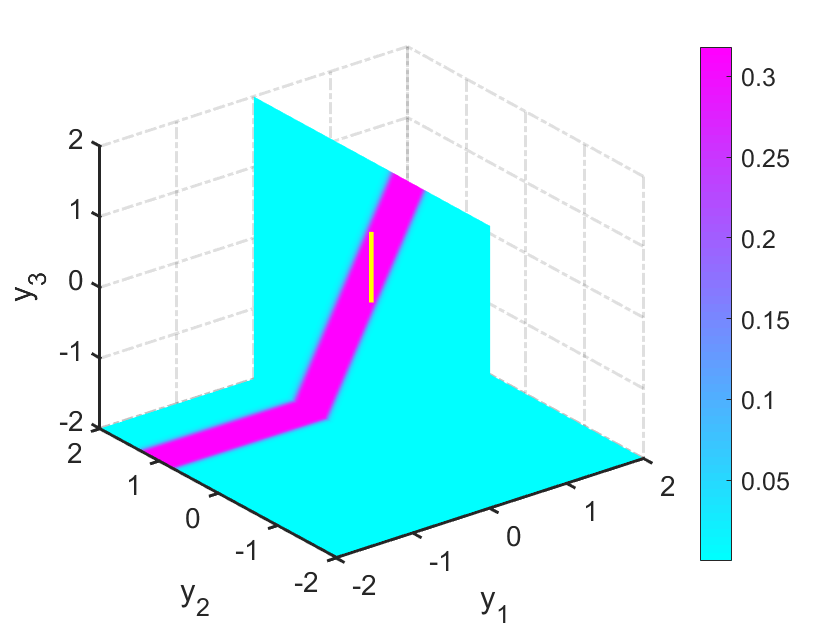}

}
\subfigure[$\phi=\pi, \theta=\pi/8$]{

\includegraphics[scale=0.22]{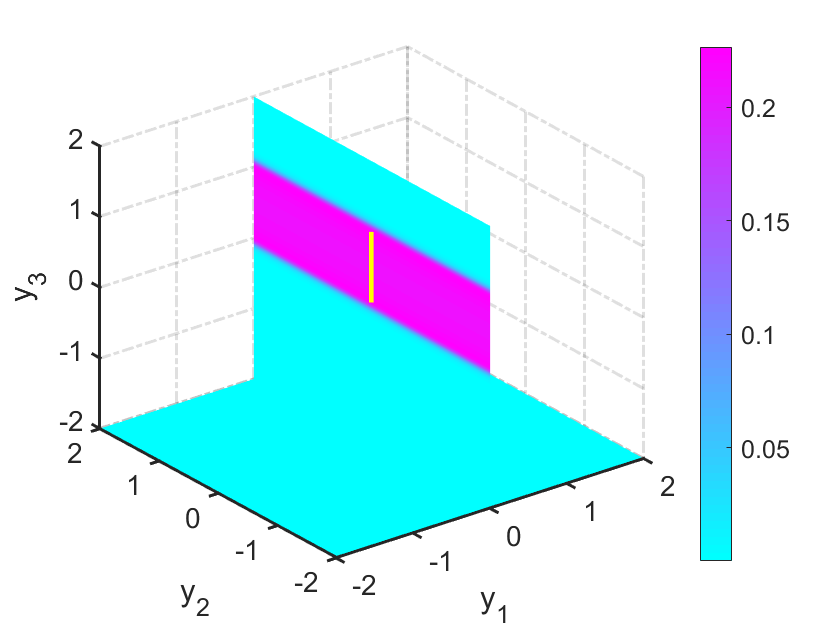}

}
\subfigure[$\phi=\pi, \theta=2\pi/8$ ]{
\includegraphics[scale=0.22]{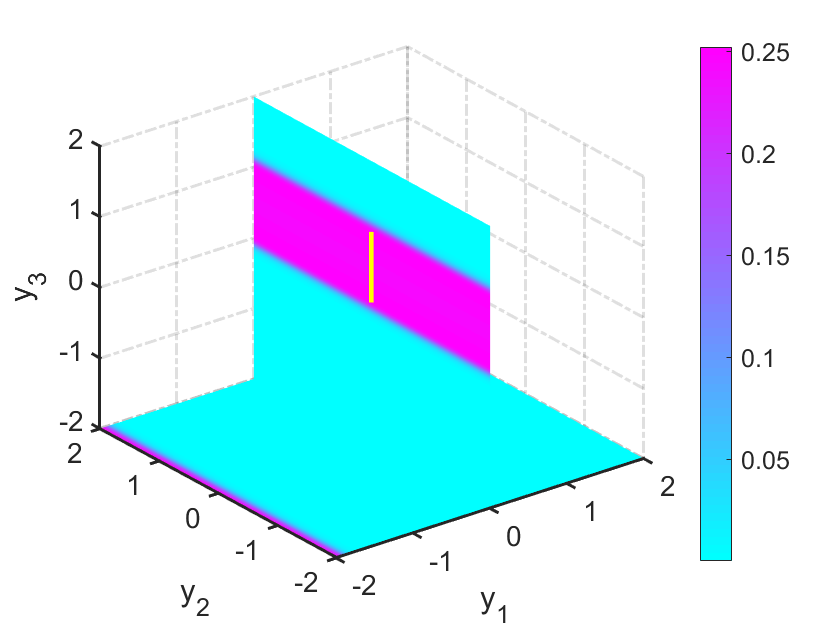}

}
\subfigure[$\phi=\pi, \theta=3\pi/8$]{
\includegraphics[scale=0.22]{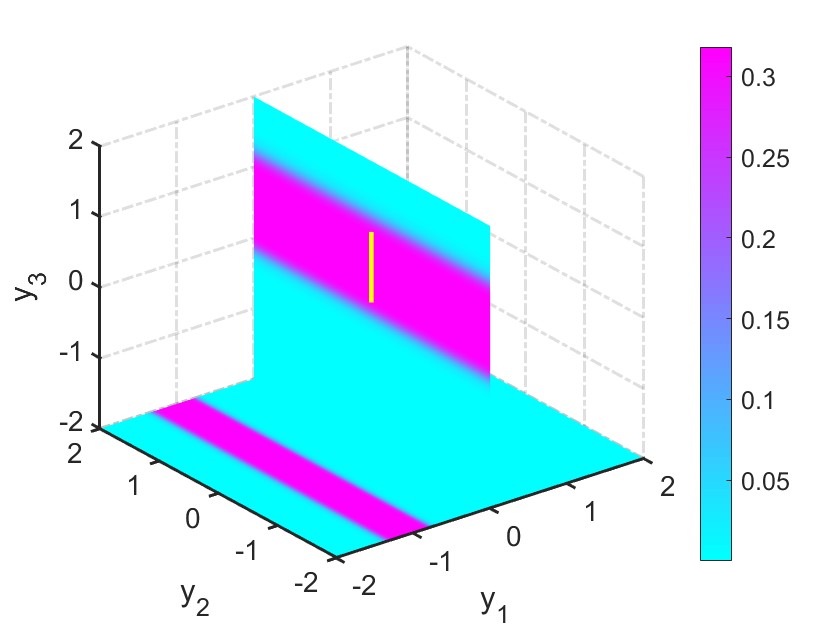}

}
\caption{Reconstruction from a single observable direction $\hat{x}=(\sin\theta\cos\varphi, \sin\theta\sin\varphi,\cos\theta)$ with $\theta\in[0,\pi/2]$ and $\phi \in [0,2\pi)$ for a straight line segment  $a(t)=(0,0,t)$ with $t\in[0,1]$ in $\R^3$. Here we  take slices at $x_1=0$ and $x_3=-2$. } \label{fig:3dline1}
\end{figure}

\begin{figure}[htb]
\centering
\subfigure[$\phi=5\pi/4, \theta=5\pi/8$]{

\includegraphics[scale=0.22]{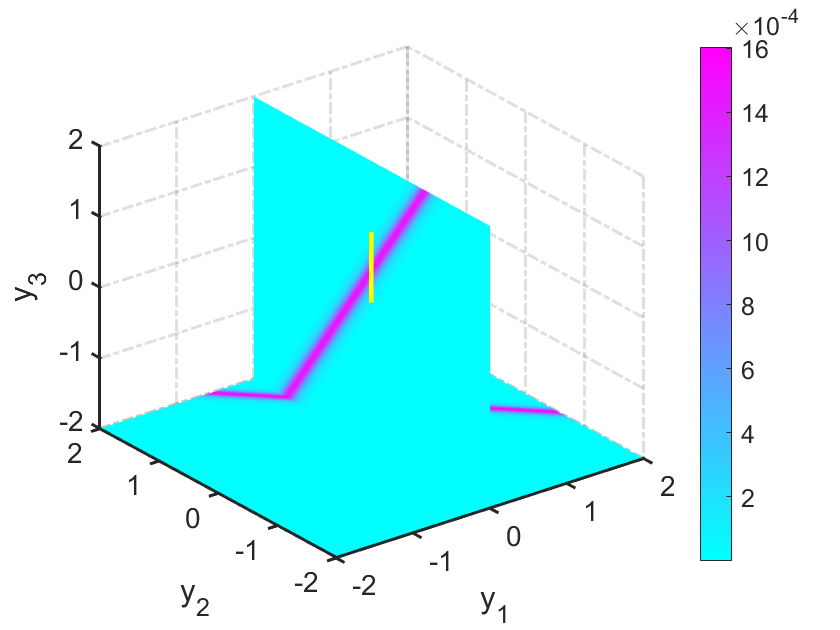}

}
\subfigure[$\phi=5\pi/4, \theta=6\pi/8$ ]{
\includegraphics[scale=0.22]{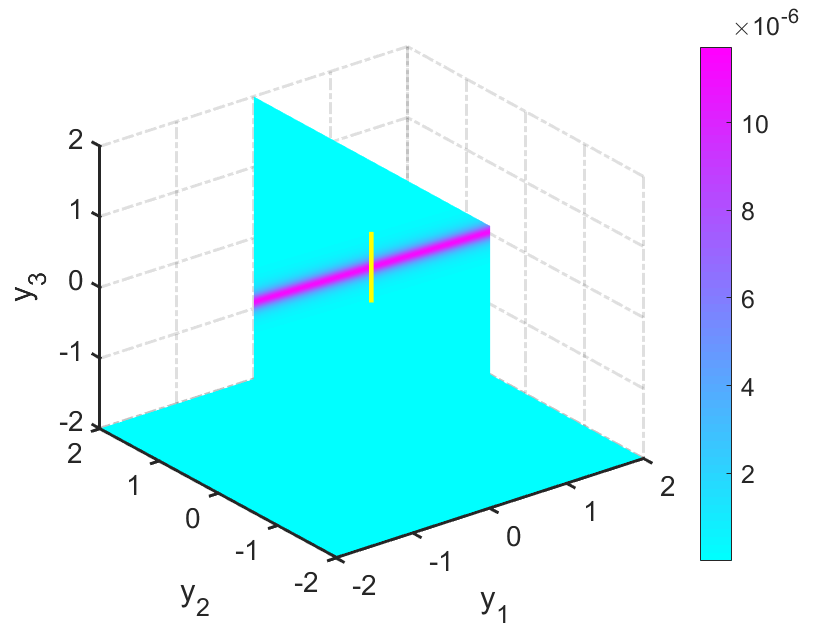}

}
\subfigure[$\phi=5\pi/4, \theta=7\pi/8$]{
\includegraphics[scale=0.22]{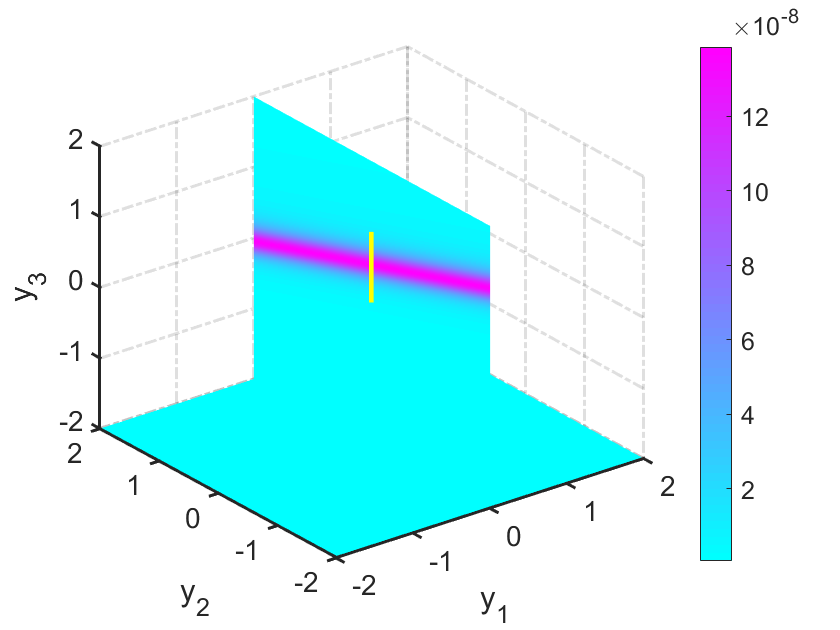}

}

\subfigure[$\phi=6\pi/4, \theta=5\pi/8$]{

\includegraphics[scale=0.22]{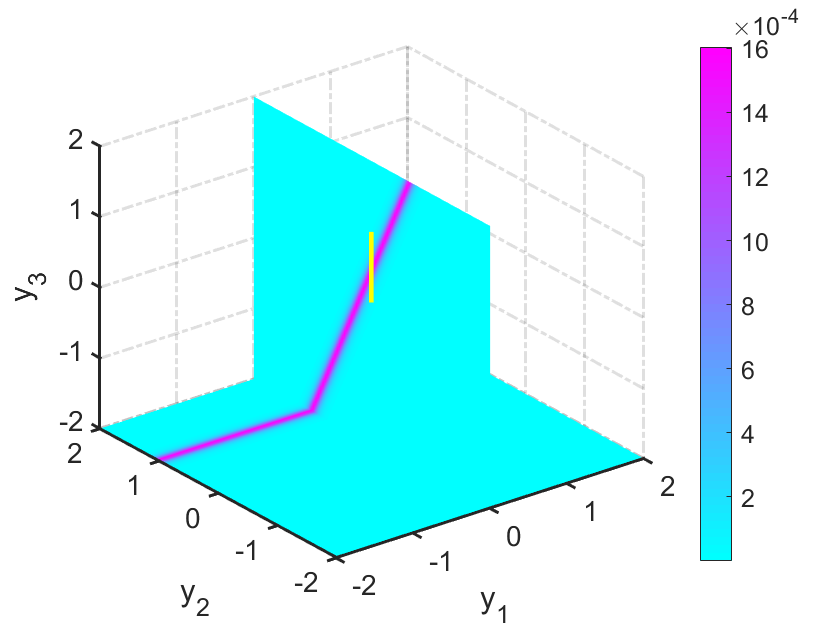}

}
\subfigure[$\phi=6\pi/4, \theta=6\pi/8$ ]{
\includegraphics[scale=0.22]{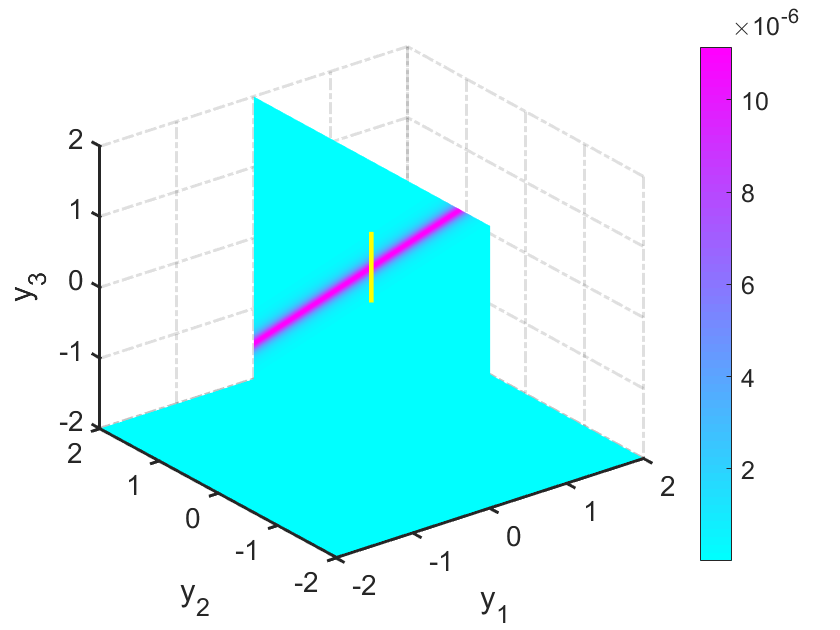}

}
\subfigure[$\phi=6\pi/4, \theta=7\pi/8$]{
\includegraphics[scale=0.22]{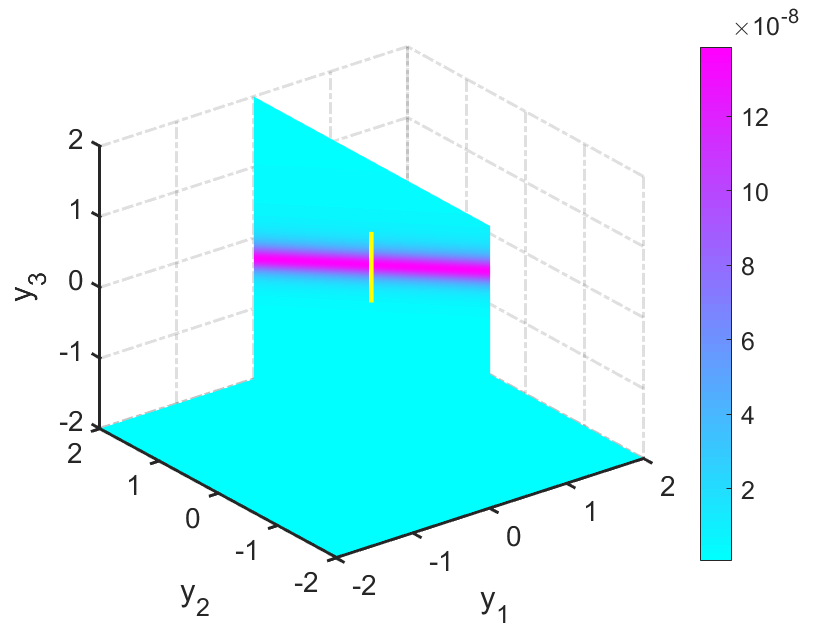}

}
\subfigure[$\phi=7\pi/4, \theta=5\pi/8$]{

\includegraphics[scale=0.22]{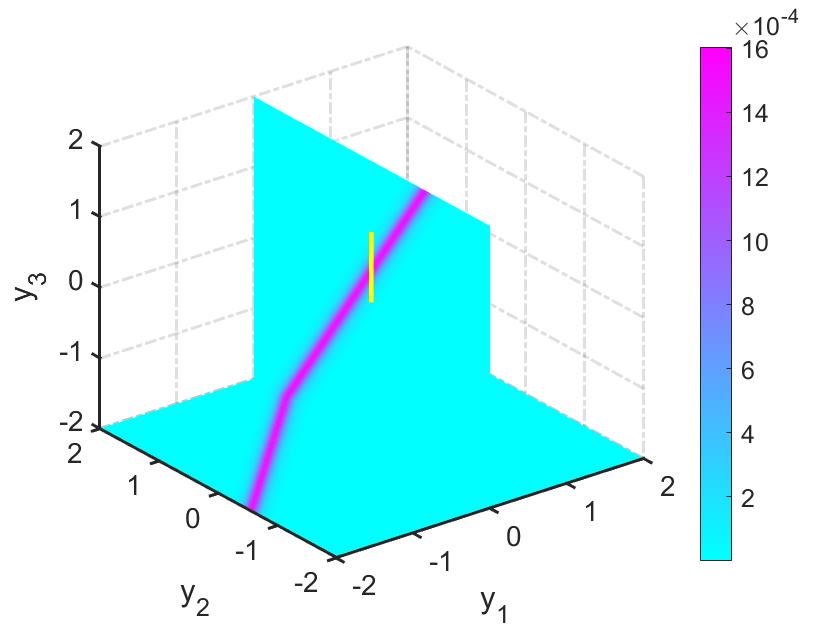}

}
\subfigure[$\phi=7\pi/4, \theta=6\pi/8$ ]{
\includegraphics[scale=0.22]{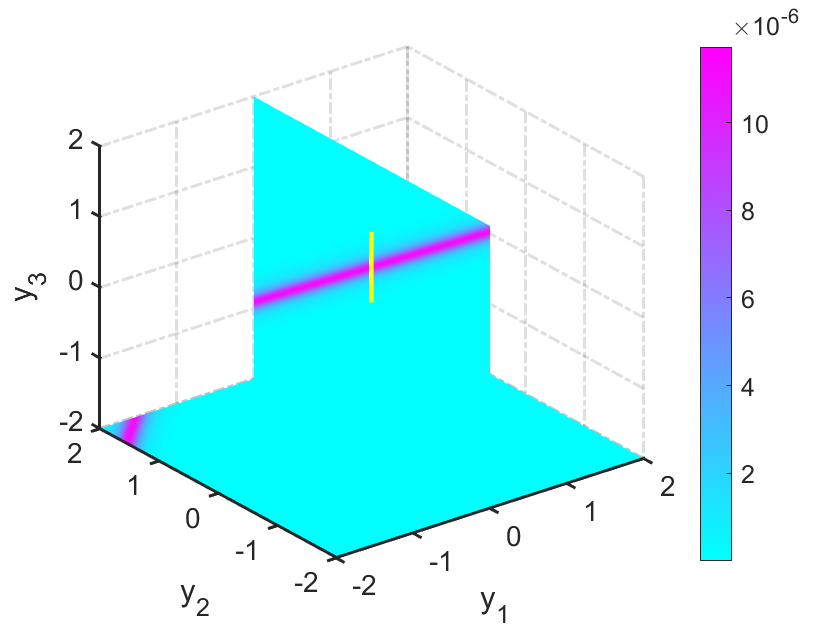}

}
\subfigure[$\phi=7\pi/4, \theta=7\pi/8$]{
\includegraphics[scale=0.22]{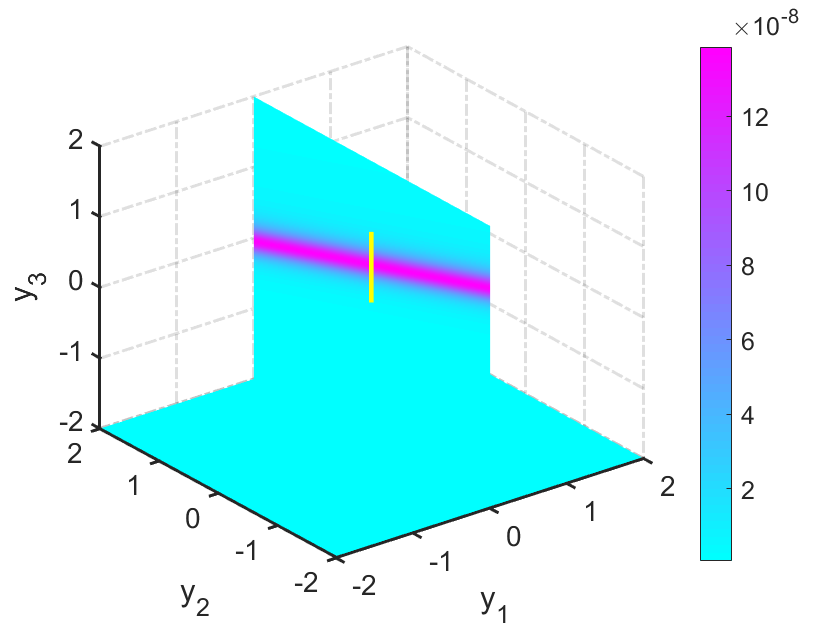}

}
\caption{Reconstruction from a single non-observable direction $\hat{x}= (\cos\theta, \sin\theta)$ with $\theta\in(\pi/2,\pi)$ and $\phi \in [0,2\pi)$ for a straight line segment  $a(t)=(0,0,t)$ with $t\in[0,1]$ in $\R^3$. Here we take take slices at $x=0$ and $z=-2$.  } \label{fig:3dline2}
\end{figure}

Consider a straight line segment in  $\R^3$ parameterized by $a(t) = (0,0,t),\,t\in[0,1]$ and write the observation direction as $\hat{x} = (\sin\theta\cos \varphi,\sin \theta \sin \varphi, \cos \theta),\,\theta\in[0,\pi],\,\varphi\in[0,2\pi)$. Then, 
$$ h(t) =t+\hat{x}\cdot a(t)= t(1+\cos \theta),\,\quad h^{\prime}(t) = 1+\hat{x}\cdot a^{\prime}(t) = 1+ \cos \theta.$$
It follows that $h^{\prime}(t)>0$ for all $t\in[t_{\min},t_{\max}]$. Hence $\hat{x}$ is a non-observable direction only if $\cos \theta<0$, that is $\theta \in (\pi/2,\pi)$, and $\hat{x}$ is an observable direction if $\theta \in [0,\pi/2]$. 
In Fig.\ref{fig:hyperplane}, we illustrate two planes perpendicular to the observable direction, between which  the trajectory of the moving source is located. Fig.\ref{fig:3dline1} presents slices of
the smallest hyperspace at $x_1=0$ and $x_3=-2$ reconstructed from the data of different observable directions.  We conclude that the trajectory of the moving source lies perfectly between the two planes that are perpendicular to the observation direction. It demonstrates effectiveness of our algorithm for imaging a straight line segment in $\R^3$.
In Fig.\ref{fig:3dline2}, we plot the indicator functions with different non-observable directions. The values of the indicator function are much smaller than $10^{-3}$.  

\begin{remark} Let us discuss the width $l^{(\hat{x})}$ of the strip $K_{\Gamma}^{(\hat{x})}$. If $\hat{x}$ is observable and $h'(t)$ remains positive, we know $l^{(\hat{x})}=\sup(\hat{x}\cdot \Gamma)-\inf(\hat{x}\cdot \Gamma)$; If $\hat{x}$ is observable and $h'(t)<0$ in $(t_{\min}, t_{\max})$, then $l^{(\hat{x})}<\sup(\hat{x}\cdot \Gamma)-\inf(\hat{x}\cdot \Gamma)$; If the direction $\hat{x}$ in the latter case is getting closer to some non-observable direction, our numerical tests show that $l^{(\hat{x})}$ tends to be thinner and thinner.
\end{remark}

\subsection{Multiple observation directions}

In this subsection, we continue  the two dimensional Examples 1, 2  and 3 but with multi-frequency far-field data measured at sparse directions.  We should truncate the indicator function (\ref{W})  by
\be \label{W1}
W(y):=\left[\sum _{j=1}^M\sum_{n=1}^N\frac{\left| \phi^{(\hat{x}_{j})}_{y}\cdot \overline{\psi_n^{(\hat{x}_{j})} }\right|^2}{ |\lambda_n^{(\hat{x}_{j})}|}\right]^{-1}, \quad y\in \R^2.
\en
where $M>0$ denotes the number of sparse observation directions equally lying on $\s^1$, the test function $\phi_y^{(\hat{x}_j)}$ is again given
by (\ref{testn}) and $\left\{(\lambda_n^{(\hat{x}_{j})} , \psi_n^{(\hat{x}_{j})}): n=1,\cdots,N\right\} $ denote an eigensystem of the operator $(\mathcal F^{(\hat{x}_j)})_\#$.
It is worthy noting that $\hat{x}_j$ ($j=1,2,\cdots, M$) may contain both observable and non-observable direction.
We set a threshold $M'>0$ to remove the contributions of the terms likes
$$\tilde w_j=\sum_{n=1}^N\frac{\left| \phi^{(\hat{x}_j)}_{y} \cdot \overline{\psi_n^{(\hat{x})} }\right|^2}{ |\lambda_n^{(\hat{x})}|},\quad j=1,2,...W$$
to the sum in (\ref{W1}). More precisely, if $\min(\tilde w_j (y))>M'$, the direction $\hat{x}_j$ can be considered as a non-observable direction by the second assertion of Theorem \ref{Th:factorization}.
In our numerical examples, the threshold value is set as $M=3.5\times10^3.$

\begin{figure}[htb]
\centering
\subfigure[$M=2$]{

\includegraphics[scale=0.22]{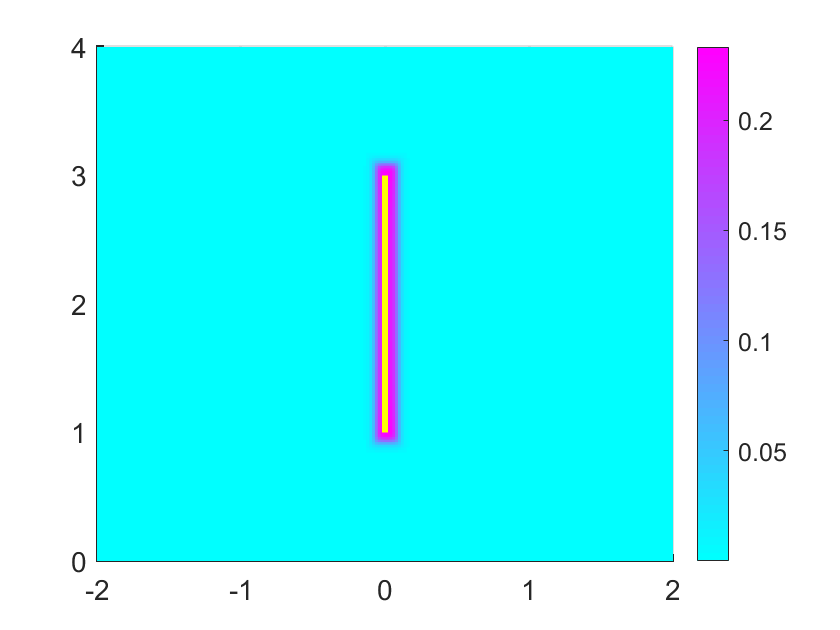}

}
\subfigure[$M=4$ ]{
\includegraphics[scale=0.22]{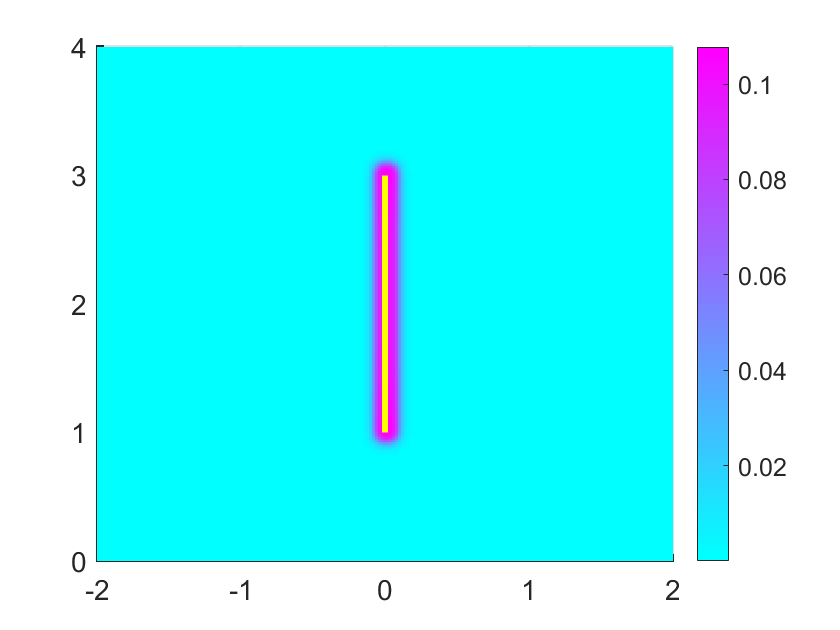}

}
\subfigure[$M=8$ ]{
\includegraphics[scale=0.22]{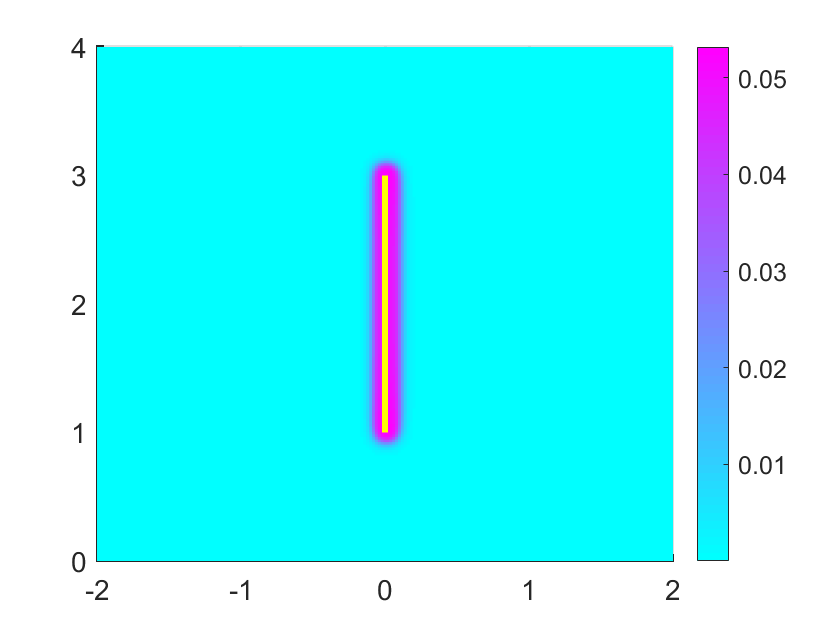}

}
\caption{Reconstruction from multiple observation  direction $\hat{x}=(\cos \theta, \sin\theta)$ with $\theta\in[0,2\pi)$ for a straight line segment $a(t)=(0,t)$ with $t\in[1,3]$. Here $M$ denotes the number of the directions. (a) $\theta=0,\pi/2$; \; (b) and (c) $\theta=(j-1)*2\pi/M$, $j=1,\cdots, M$.
 } \label{line-multi}
\end{figure}

We present in Fig.\ref{line-multi} a visualization of the reconstructed trajectory for orbit functions $a(t)=(0,t)$ with $t\in[1,3]$ with multiple observation directions. For $M=2,4,8$, there exists at one direction perpendicular to the trajectory and one parallel to the trajectory,  the intersections of the strips $K_\Gamma^{(\hat{x}_j)}$ always reflect the trajectory of the moving source. Since $h'(t)>0$  for all observable directions in Example 1,  the trajectory can be perfectly reconstructed from the data taken on sparse observation directions.

\begin{figure}[htb]
\centering
\subfigure[$M=2$]{

\includegraphics[scale=0.22]{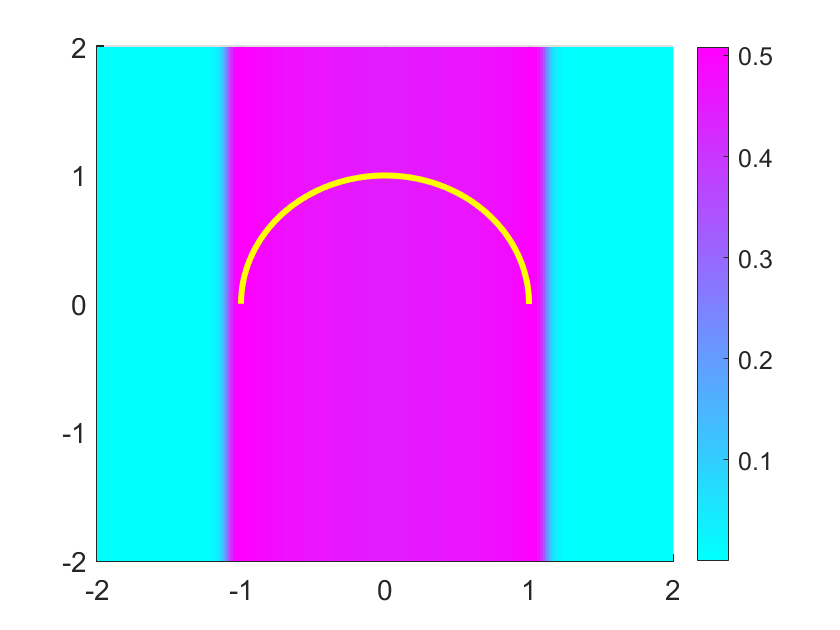}

}
\subfigure[$M=4$ ]{
\includegraphics[scale=0.22]{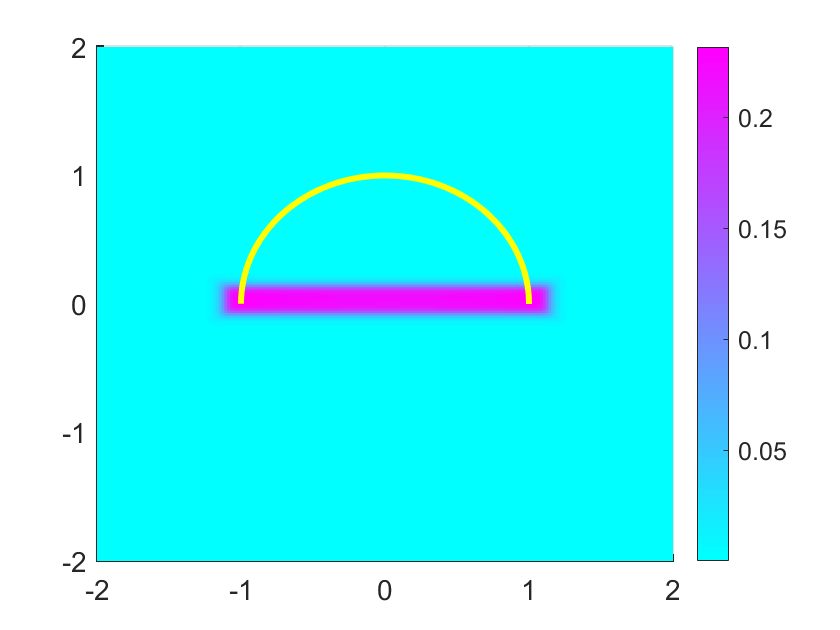}

}
\subfigure[$M=6$ ]{
\includegraphics[scale=0.22]{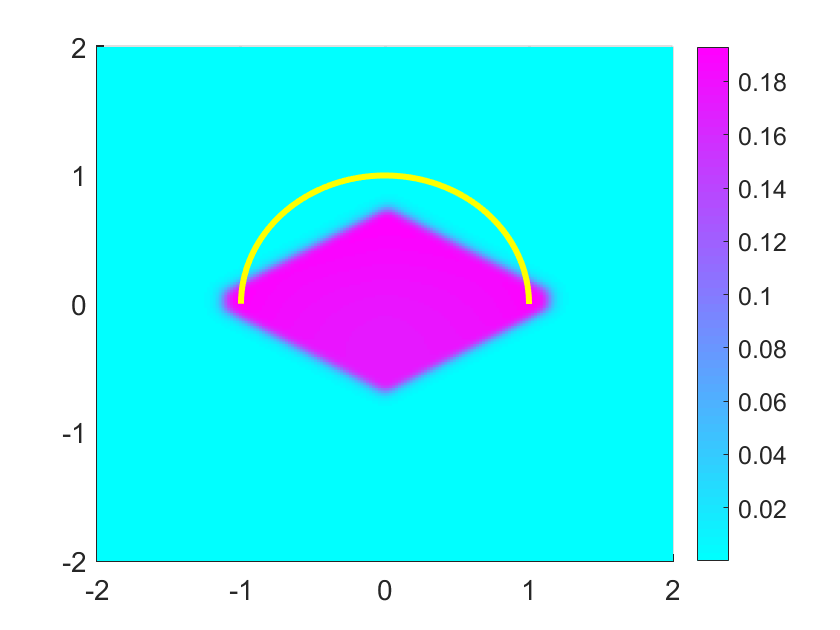}

}
\subfigure[$M=8$]{

\includegraphics[scale=0.22]{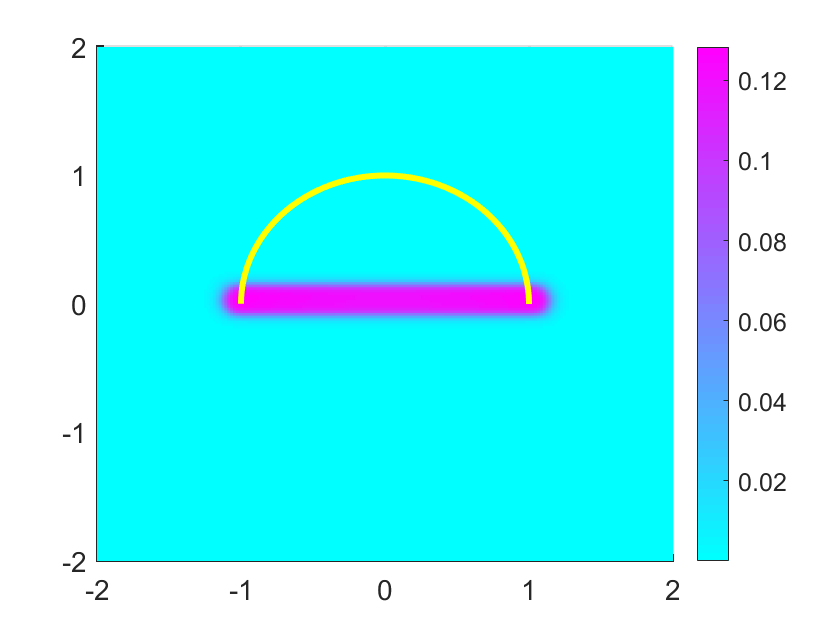}

}
\subfigure[$M=10$ ]{
\includegraphics[scale=0.22]{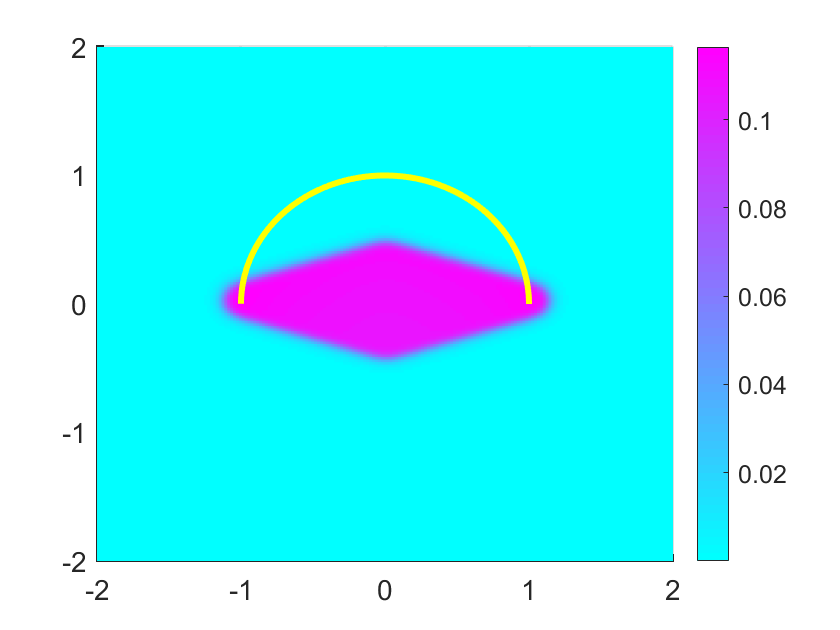}

}
\subfigure[$M=12$ ]{
\includegraphics[scale=0.22]{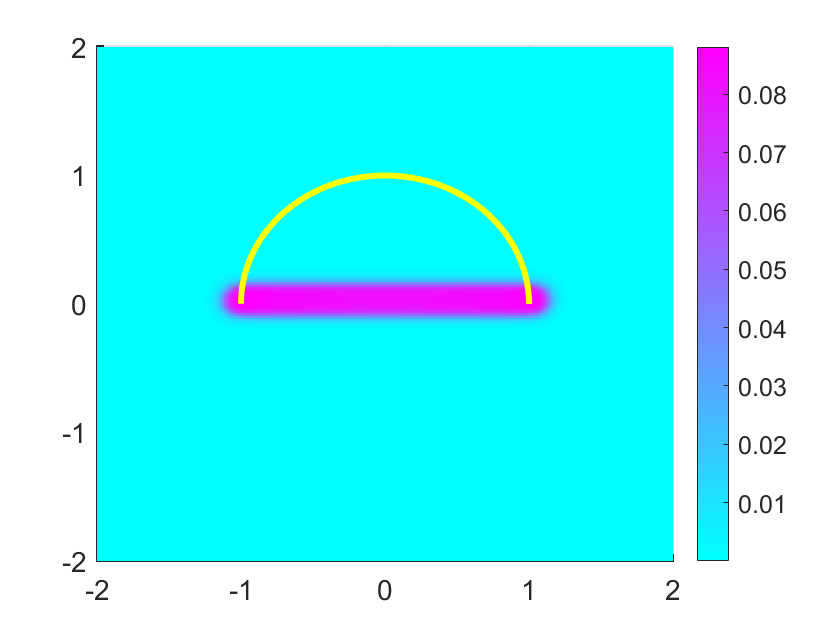}

}
\caption{Reconstruction from  multiple observation directions $\hat{x}=(\cos\theta, \sin\theta)$ with $\theta \in (0, 2\pi)$  for an arc  $a(t)=(\cos t, \sin t)$ with $t\in[0,\pi]$. Here $M$ is the number of the directions.  $\theta=(j-1)*2\pi/M$, $j=1,2,...M$. } \label{circle-multi}
\end{figure}

\begin{figure}[htb]
\centering
\subfigure[$M=2$]{

\includegraphics[scale=0.22]{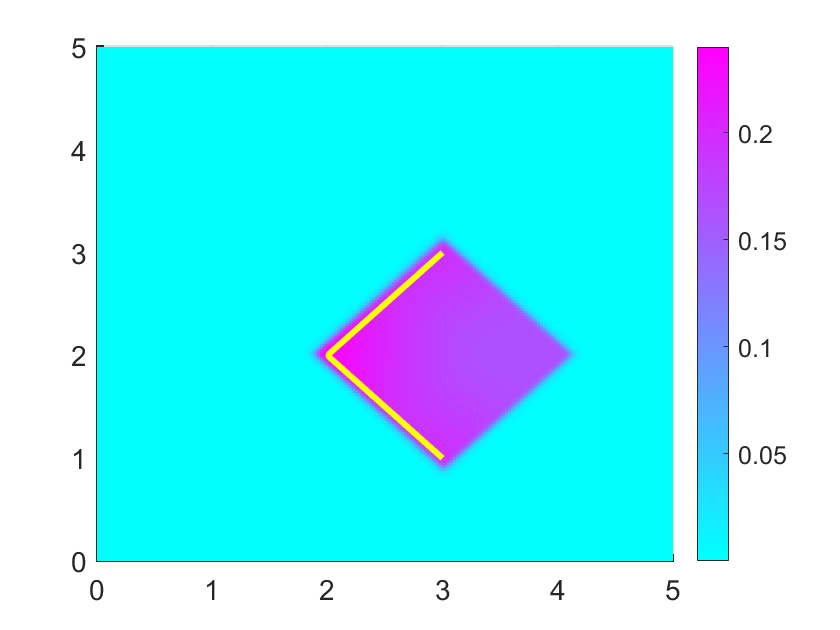}

}
\subfigure[$M=4$ ]{
\includegraphics[scale=0.22]{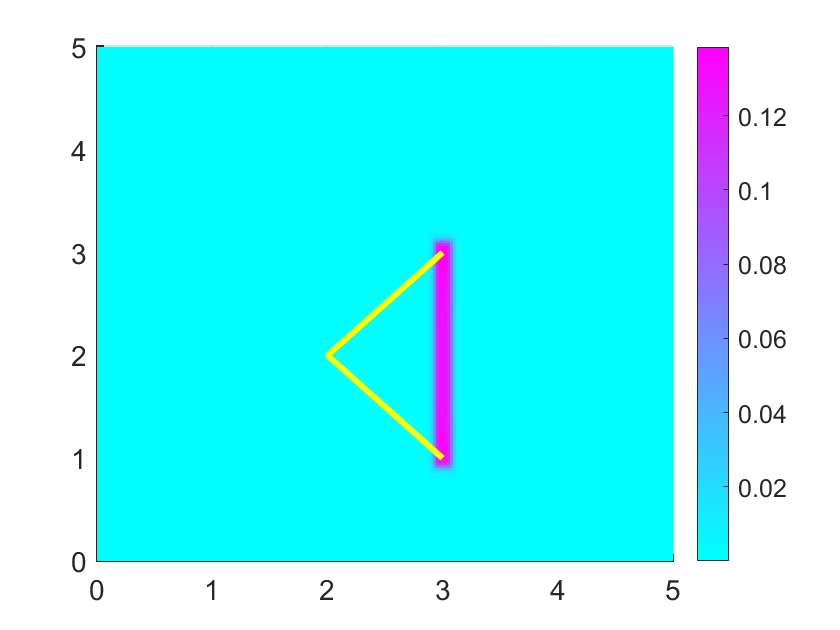}

}
\subfigure[$M=7$ ]{
\includegraphics[scale=0.22]{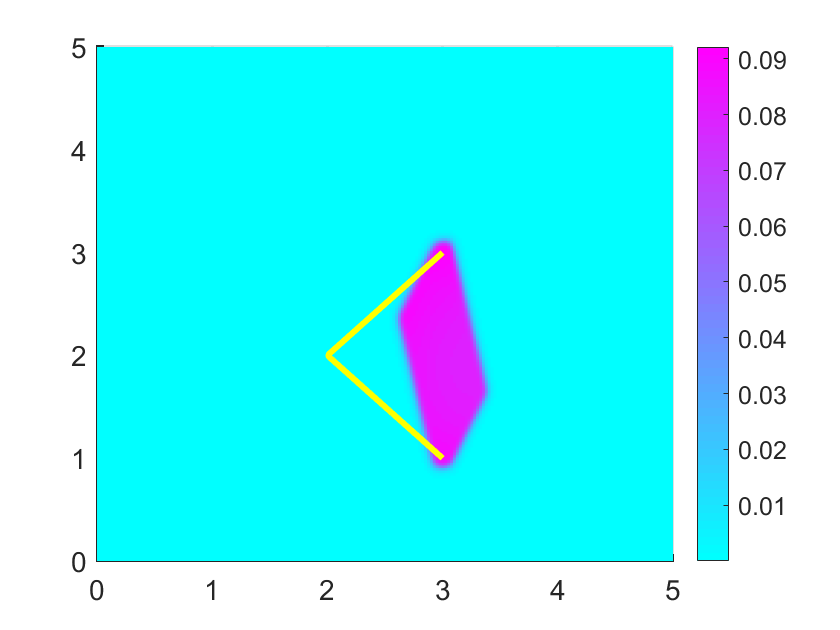}

}
\caption{Reconstruction from multiple observation directions $\hat{x}=(\cos\theta, \sin\theta)$ with $\theta\in[0,2\pi)$  for the piecewise linear curve $a(t)=(-t+3,-t+3)$ with $t\in[0,1]$ and $a(t)=(t+1,-t+3)$ with $t\in[1,2]$ in $\R^2$. Here $M$ is the number of observation directions. (a) $\theta=5\pi/4,7\pi/4$;\; (b) $\theta=(j-1)*2\pi/M, j=1,...,M$;\; (c) $\theta=2j*2\pi/15, j=1,...,M$.
 } \label{fig:xie-line2}
\end{figure}

However, in the case of the line segment in Example 3 or the arc in Example 4, we can only get partial information on the trajectory. From Figs.\ref{circle-multi} and \ref{fig:xie-line2}, one can only get the starting and ending points of the trajectory, although the data of
 multiple directions are put into use. This is due to the existence of
 $\hat x_j$ satisfying $K_{\Gamma}^{(\hat x_j)}\subset \{y\in \R^2: \sup(\hat x \cdot \Gamma)\leq \hat x \cdot y\leq \inf (\hat x \cdot \Gamma) \} $. For such observation directions, the width of the reconstructed strip $K_{\Gamma}^{(\hat x_j)}$ is very small. Hence, the intersection of $K_{\Gamma}^{(\hat x_j)}$ always appears like a line segment connecting the starting and the ending points of the trajectory.

%

\subsection{Reconstructions from noisy data}
 We test the sensitivity of the algorithm with respect to the noisy data. Consider the Case 1 in Example 1 for recovering a line segment. The far-field data are polluted by Gaussian noise in the form of
 $$
 w^{\infty}_\delta(\hat x,k)\coloneqq \real [w^{\infty}(\hat x,k)]\;\big(1+\delta\, \gamma_1\big) + {\rm Im}. [w^{\infty}(\hat x,k)]\;\big(1+\delta\, \gamma_2\big)
 $$
 where $\delta>0$ denotes the noise level and $\gamma_j\in[-1,1]$ $(j=1,2)$ are Gaussian random variables.

\begin{figure}[htb]
\centering
\subfigure[$\theta=0$]{
\includegraphics[scale=0.35]{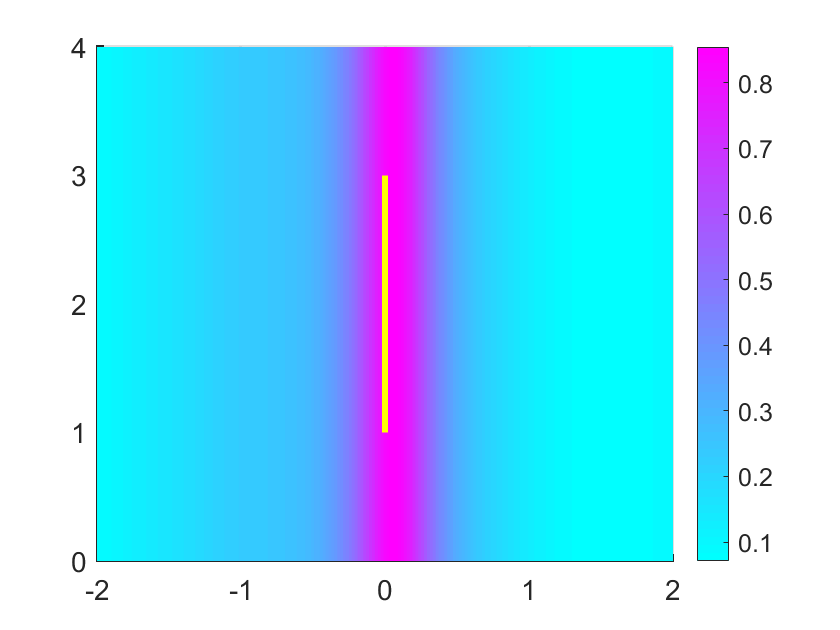}}
\subfigure[$\theta=\pi/4 $ ]{
\includegraphics[scale=0.35]{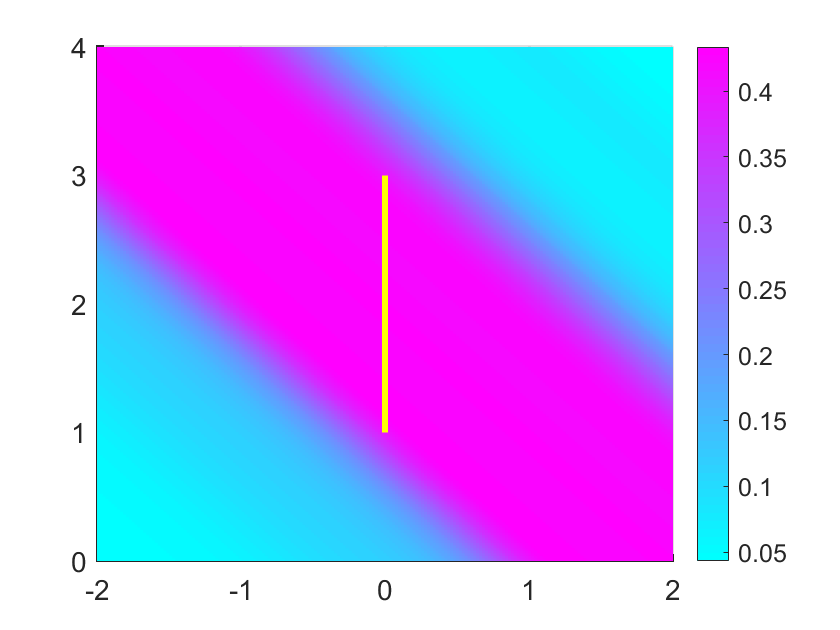}}
\subfigure[$\theta=\pi/2$]{
\includegraphics[scale=0.35]{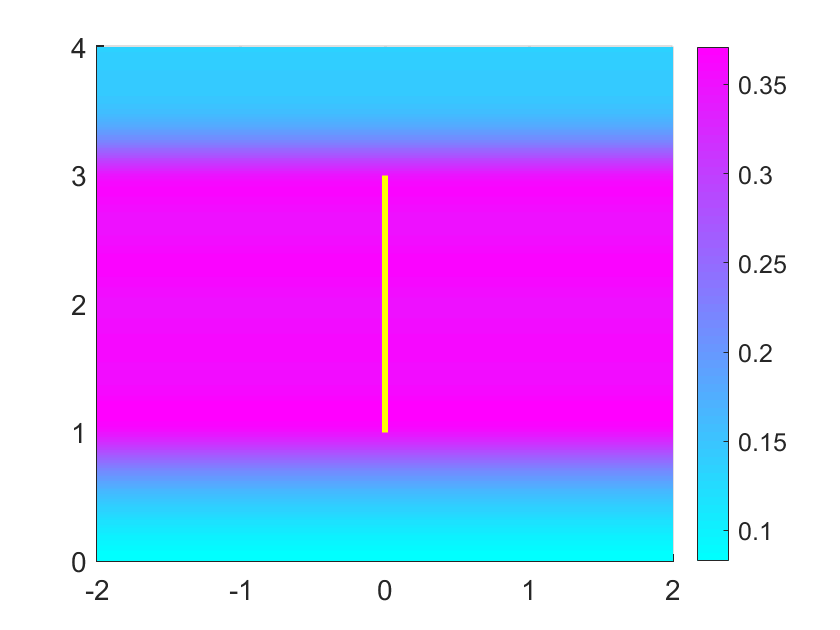}
}\subfigure[$M=4$]{
\includegraphics[scale=0.35]{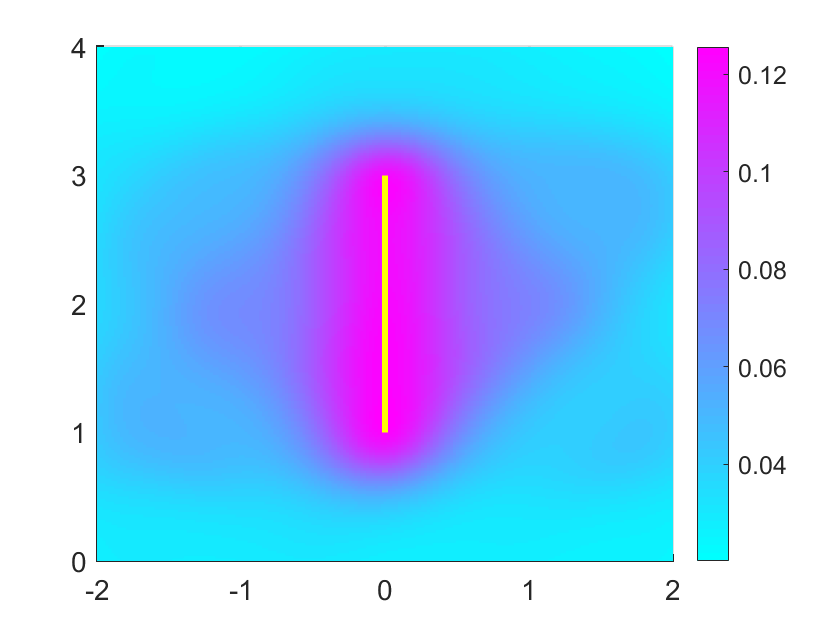}
}



%
\caption{Reconstruction of  a straight line segment $a(t)=(0,t)$, $t\in[1,3]$ from noisy data measured at a single observable direction $\hat{x}=(\cos \theta, \sin \theta)$ in (a),(b) and (c). In (d),  the polluted far-field data from $M=4$ directions are used. The noise level is set as $\delta=1$\%. } \label{fig:noise}
\end{figure}

We set $\delta=1\%$ and plot the indicator functions in Fig.\ref{fig:noise} using one and sparse observation directions. It turns out that the proposed scheme is rather sensitive to noise. Even at the noise level 1\%, one can only get a rough location of the trajectory of the moving source using the data measured at sparse directions. This shows that our inverse problems are severely ill-posed. However, a quantitive characterization of the ill-posed nature remains unclear to us.

\section*{Acknowledgements}

G. Hu is partially supported by the National Natural Science Foundation of China (No. 12071236) and the Fundamental Research Funds for Central Universities in China (No. 63213025).


\begin{thebibliography}{99}

\bibitem{AHLS} A. Alzaalig, G. Hu, X. Liu and J. Sun, Fast acoustic source imaging using multi-frequency sparse data, Inverse Problems, 36 (2020): 025009.


\bibitem{CGMS2020} B. Chen, Y. Guo, F. Ma and Y. Sun, Numerical schemes to reconstruct three-dimensional time-dependent point sources of acoustic waves, Inverse Problems, 36 (2020): 075009.


\bibitem{CB2008} M. Cheney and B. Borden,   Imaging moving targets from scattered waves, Inverse Problems 24 (2008): 035005.


\bibitem{Co79}J. Cooper, Scattering of plane waves by a moving obstacle. Arch. Ration. Mech. Anal. 71 (1979): 113-149.


\bibitem{CoStra} J. Cooper and W. Strauss, Scattering of waves by periodically moving bodies. J. Funct. Anal. 47 (1982): 180-229.






\bibitem{FGPT2017} J. Fournier, J. Garnier, G. Papanicolaou and C. Tsogka,  Matched-filter and correlation-based imaging for fast moving objects using a sparse network of receivers, SIAM J. Imag. Sci., 10 (2017): 2165-2216.






\bibitem{GF2015} J. Garnier and M. Fink, Super-resolution in time-reversal focusing on a moving source, Wave Motion, 53 (2015): 80-93.




\bibitem{GS2017} R. Griesmaire and C. Schmiedecke, A Factorization method for multifrequency inverse source problem with sparse far field measurements, SIAM J. Imag. Sci., 10 (2017): 2119-2139.

\bibitem{GGH2022} R. Griesmarier, H. Guo, G. Hu, Inverse wave-number-dependent source problems for the Helmholtz equation, in preparing.


\bibitem{GHZ22} H. Guo, G. Hu and M. Zhao,
Direct sampling method to inverse wave-number-dependent source problems (part I): determination of the support of a stationary source,
arXiv:2212.04806.


\bibitem{HKLZ2019} G. Hu, Y. Kian, P. Li and Y. Zhao, Inverse moving source problems in electrodynamics, Inverse Problems, 35 (2019): 075001.

\bibitem{HKZ2020} G. Hu, Y. Kian and Y. Zhao, Uniqueness to some inverse source problems for the wave equation in unbounded domains, Acta Mathematicae Applicatae Sinica, English Series, 36 (2020): 134-150.

\bibitem{HLY20} G. Hu, Y. Liu  and M. Yamamoto, Inverse moving source problem for fractional diffusion(-wave) equations: Determination of orbits, Inverse Problems and Related Topics ed J Cheng, S Lu and M Yamamoto (Singapore: Springer) pp. 81-100, 2020.


\bibitem{Isa1989} V. Isakov, Inverse Source Problems, AMS, Providence, RI, 1989.

\bibitem{T2022}H. A. Jebawy, A. Elbadia and F. Triki, Inverse moving point source problem for the wave equation, Inverse Problems, 38 (2022): 125003.


\bibitem{KG08} A. Kirsch and N. Grinberg, The Factorization Method for Inverse Problems, Oxford University Press, Oxford, UK, 2008.

\bibitem{LGS2021} Y. Liu, Y. Guo, and J. Sun, A deterministic-statistical approach to reconstruct moving sources using sparse partial data, Inverse Problems, 37 (2021): 065005.


\bibitem{LGY21} Y. Liu, G. Hu and M. Yamamoto, Inverse moving source problem for time-fractional evolution equations: determination of profiles, Inverse Problems, 37 (2021): 084001.


\bibitem{Liu2021}Y. Liu, Numerical schemes for reconstructing profiles of moving sources in (time-fractional) evolution equations, RIMS Kokyuroku 2174 (2021): 73-87.




\bibitem{NIO2012} E. Nakaguchi, H. Inui and K. Ohnaka. An algebraic reconstruction of a moving point source for a scalar wave equation, Inverse Problems, 28 (2012): 065018.

\bibitem{Ohe2011} T. Ohe, H. Inui and K. Ohnaka,  Real-time reconstruction of time-varying point sources in a three-dimensional scalar wave equation. Inverse Problems, 27 (2011): 115011.



\bibitem{S1991} P. D. Stefanov, Inverse scattering problem for moving obstacles, Math. Z. 207 (1991): 461-480.



\bibitem{SK} J. Sylvester and J. Kelly, A scattering support for broadband sparse far-field measurements, Inverse Problems, 21 (2005): 759-771.




\bibitem{T2020} O. Takashi. Real-time reconstruction of moving point/dipole wave sources from boundary measurements, Inverse Probl. Sci. Eng., 28 (2020): 1057-1102.









\bibitem{WKT2022}S. Wang, Mirza Karamehmedovic, Faouzi Triki, Localization of moving sources: uniqueness, stability and Bayesian inference, 
arXiv:2204.04465.




\end{thebibliography}
\end{document}